\documentclass[11pt]{article}

\usepackage[a4paper, margin=1in]{geometry}
\usepackage{graphicx,times,amsmath,amsthm,amsfonts,amssymb,mathtools,color,latexsym,amsxtra,layout,bibentry,hyperref}
\usepackage{subfig}
\usepackage{graphicx}
\usepackage[shortlabels]{enumitem}

\usepackage{array}
\newcolumntype{L}[1]{>{\raggedright\let\newline\\\arraybackslash\hspace{0pt
}}m{#1}}
\newcolumntype{C}[1]{>{\centering\let\newline\\\arraybackslash\hspace{0pt}}
m{#1}}
\newcolumntype{R}[1]{>{\raggedleft\let\newline\\\arraybackslash\hspace{0pt}
}m{#1}} 

\usepackage{times,amsmath,stackengine,amsthm,amsfonts,amssymb,mathtools,color,latexsym,amsxtra,layout}
\usepackage{graphicx}
\usepackage{algorithmicx, algpseudocode}
\usepackage{algorithm}
\usepackage{multirow}
\usepackage{comment, centernot}

\newcommand{\Vtinner}[1]{\left \langle #1  \right \rangle_{H^1_{\alpha}}}

\newcommand{\ltwoinner}[1]{\left \langle #1  \right \rangle_{2}}
\newcommand{\norm}[1]{\left \| #1 \right \|_2}
\newcommand{\Mnorm}[1]{\left \| #1 \right \|_M}

\newcommand{\Vnorm}[1]{\left \| #1 \right \|_{V}}
\newcommand{\Vtnorm}[1]{\left \| #1 \right \|_{H^1_{\alpha}}}

\newcommand{\ltwonorm}[1]{\left \| #1 \right \|_{2}}

\newcommand{\linftynorm}[1]{\left \| #1 \right \|_{\infty}}
\newcommand{\lhnorm}[1]{\left \| #1 \right \|_{h}}

\newcommand{\oneinftynorm}[1]{\left \| #1 \right \|_{1,\infty}}
\newcommand{\honenorm}[1]{\left \| #1 \right \|_{H^1}}
\newcommand{\Hnorm}[1]{\left \| #1 \right \|_{H}}

\theoremstyle{plain}
\newtheorem{theorem}{Theorem}[section]
\newtheorem{lemma}[theorem]{Lemma}

\newtheorem{o-thm}[theorem]{Theorem}

\newtheorem{corollary}[theorem]{Corollary}

\theoremstyle{definition}

\makeatletter							% Title
\def\printtitle{%						
    {\centering \@title\par}}
\makeatother									

\makeatletter							% Author
\def\printauthor{%					
    {\centering \large \@author}}				
\makeatother							

\title{Direct and Inverse Problem for Gas Diffusion in Polar Firn}
\author{Sophie Moufawad \thanks{American University of Beirut (AUB), Beirut, Lebanon.  (sm101@aub.edu.lb) } \and Nabil Nassif \thanks{American University of Beirut (AUB), Beirut, Lebanon.  (nn12@aub.edu.lb) } \and Faouzi Triki  \thanks{Grenoble-Alpes University, Grenoble, France.  (faouzi.triki@univ-grenoble-alpes.fr) }}
		
\begin{document}

\graphicspath{{fig/}}

%%%%%%%%%%%%%%%%%%%%%%%%%%%%%%%%%%%%%%%%%%%%%%%%%%%%%%%

%%%%%%%%%%%%%%%%%%%%%%%%%%%%%%%%%%%%%%%%%%%%%%%%%%%%%%%

\maketitle
%%%%%%%%%%%%%%%%%%%%%%%%%%%%%%%%%%%%%%%%%%%%%%%%%%%%%%%
\begin{abstract}
\noindent Simultaneous use of partial differential equations in conjunction with data analysis has proven to be an efficient way to obtain the main parameters of various phenomena in different areas, such as medical, biological, and ecological. In the ecological field, the study of climate change (including global warming) over the past  centuries requires estimating different gas concentrations in the atmosphere, mainly CO2. 

\noindent The mathematical model of gas trapping in deep polar ice (firns) has been derived in \cite{Witrant2012acp, Yeung, Laube, NEEM}, consisting of a parabolic partial differential equation that is almost degenerate at one boundary extreme. 
 In this paper, we consider all the coefficients to be constants, except the diffusion coefficient that is to be reconstructed. 
 We present the theoretical aspects of existence, uniqueness and simulation for such direct problem 
 and consequently formulate the inverse problem that attempts at recovering the diffusion coefficients using given generated data 
\end{abstract}

\noindent \textbf{Funding}:  This work was supported by the Alwaleed Center for American Studies and Research (CASAR) in the Faculty of Arts and Sciences at AUB; and by the AUB University Research Board grant number 104261 (Project 26742).\\


%\tableofcontents

%%%%%%%%%%%%%%%%%%%%%%%%%%%%%%%%%%%%%%%%%%%%%%%%%%%%%%%
\section{Introduction}
Antarctic and Greenland Polar snow and ice constitute a unique archive of past climates and atmospheres. Based on a good understanding of the mechanisms controlling gas trapping in deep polar ice, and therefore of the processes of densification and pore closure in Firns (typically over the first hundred meters of the polar cap), several models have been derived as a result of the collaborations between the \href{https://www.ige-grenoble.fr/ICE3-238}{ICE3 team} of the IGE and \href{http://www.gipsa-lab.fr/accueil.php}{GIPSA Lab} (24 publications \cite{Witrant2012acp} including 3 in Nature \cite{Yeung, Laube, NEEM}).
  
  Considering the mass conservation equations, the concentration $ \rho_\alpha ^{\rm o} $ of a gas $ \alpha$ in open pores satisfies an initial-value, time-dependent advection-diffusion partial differential equation on a one-space dimension segment $[0,z_F]$  with Dirichlet boundary condition at $0$ and a mixed one at $z_F$, for  $\; z\in (0, z_{\rm F}), \\ t>0$:
\begin{equation}\label{eq:trace_gas_dynamics}
\left\{
\begin{array}{l}
\displaystyle{\frac{\partial}{\partial t} [\rho_\alpha^{\rm o} f] + \frac{\partial}{\partial z} [\rho_\alpha^{\rm o} f ({v}+{w}_{\rm air})] + \rho_\alpha^{\rm o} (\tau+\lambda) =\frac{\partial}{\partial z} \left[D_\alpha \left(\frac{\partial \rho_{\rm \alpha}^{\rm o}}{\partial z} - {\rho}_{\alpha}^{\rm o} \frac{M_{\alpha}g}{RT} \right)\right]}, \vspace{1mm}\\
\rho_\alpha^{\rm o}(0,t) = \rho_\alpha^{\rm atm}(t), \; t>0,\vspace{2mm} \\
\displaystyle{ D_\alpha(z_F)\left(\frac{\partial {\rho}_{\alpha}^{\rm o} }{\partial z}(z_{\rm F},t) -  \frac{M_{\alpha}g}{RT} {\rho}_{\alpha}^{\rm o}(z_{\rm F},t)\right) = 0},\\
\rho_\alpha^{\rm o}(z,0) = 0
\end{array}
\right.
\end{equation}
with $\rho_\alpha^{\rm atm}(0)=0$. 
\\
Moreover, $ D_\alpha(z) $ is the effective diffusion coefficient of the gas $\alpha$ in the Firn ($m^2/yr$) and is given by 
\begin{eqnarray} \label{TT}
D_\alpha(z) = \left\{ \begin{array}{llcc}
D_{\rm{eddy}}(z)+r_\alpha c_f D_{\textrm{CO2, air}}(z)\;\; \textrm{ if } z \leq z_{\textrm{eddy}},\\
r_\alpha D_{\textrm{CO2, air}}(z),\;\; \textrm{ if } z > z_{\textrm{eddy}},
\end{array}
\right.
\end{eqnarray}
with $ z_{\textrm{eddy}}$, $ r_\alpha $, and $ c_f $ are known constants, and $D_{\rm{eddy}}(z), D_{\textrm{CO2, air}}(z)$ diffusion coefficients. 
The remaining terms are considered constants in this paper,  and
 summarized in Table \ref{tab:results}.

\begin{table}[H]\label{tab}
\setlength{\tabcolsep}{3pt}
\caption{The description of the model's parameters. }\label{tab:results}
  \centering
\begin{tabular}{|l|l|}\hline
$ z_{\rm F} $ & the depth of the Firn (m) \\ \hline
$ f $ & the average volume fraction in the open pores $\in (0,1)$ \\\hline
 $ v $ & the average descending speed in the Firn (m/yr)\\\hline
  $ w_{air} $ & the average speed of the air (m/yr) 
  \\\hline
  $ \tau $ & the  mass exchange rate between open and closed pores ($/yr$) \\\hline
 $ \lambda $ & the rate of radioactive decay ($/yr$) $\in [0.5, 0.999] $\\\hline
 $ M_\alpha $ & the molar mass of the gas ($kg/mol$) $\in [0.004, 0.133]$; \;\;\; $M_{CO_2} \approx 0.044$ \\\hline
 $ g $ & the gravitational acceleration ($m/s^2$) $\approx 9.80665$\\ \hline
 $ R $ & the universal constant of ideal gases ($J/mol/K$) = $8.314 $ \\\hline
  $ T $ & the mean temperature of the Firn ($K$) $\approx -31+273.15 \approx 242 K$ \\\hline
  $ \rho_\alpha^{\rm atm} $ & the concentration of gas in the atmosphere ($mol/m^3$ of void space) \\\hline
\end{tabular}
\end{table}

The main goal of this paper is to study the theoretical aspects of the underlying mathematical model, which is an almost singular,
parabolic partial differential equations. 
We start in section \ref{sec:2} by deriving the semi-variational form of \eqref{eq:trace_gas_dynamics}. Then, we prove in section \ref{sec:3} the existence and uniqueness of a solution to \eqref{eq:trace_gas_dynamics} by applying Lions' Theorem (\cite{Brezis}, page 341). Moreover, after rescaling \eqref{eq:trace_gas_dynamics} to the unit square, an Euler-Implicit in time and Finite Element Space discretization is proposed is section \ref{sec:4}, that leads to a robust Direct problem algorithm that is tested in section \eqref{sec:testings}.

Then, the inverse problem is formulated by defining the objective function in section \ref{sec:inv}, where its gradient is computed using directional derivatives (section \ref{sec:gradV}). Testings are performed on the efficiency of the computed gradient using MATLAB's FminUnc and FminCon functions and Nonlinear Conjugate Gradient method (section \ref{sec:test}). Finally, concluding remarks are given in section \ref{sec:conc}.    
\section{Direct Problem}
{
We start by stating the assumptions on the diffusion coefficient $D_\alpha$ and the the concentration of gas in the atmosphere $\rho_\alpha^{atm}$.  
\begin{enumerate}
    \item The Diffusion Coefficient $D_\alpha(z)$ satisfies the following properties:
    \begin{itemize}
\item Strictly positive on $[0,z_F)$,
    \begin{equation}\label{eq:Dpos}
        D_\alpha(z) > 0, \quad  \forall z\in [0,z_F)
    \end{equation}
    \item  Lipschitz continuous on $[0,z_F)$, specifically for all $\delta >0$, there exists $L_\delta$ such that:
    \begin{equation}\label{eq:Lipschitz}
    |D_\alpha(z)-D_\alpha(y)|\le L_\delta |z-y|,\quad\forall z,\,y\in[0,z_F-\delta].
   \end{equation}
    \end{itemize}
    \item Moreover, the boundary condition at $z_F$ can be satisfied in one of the following two ways:
    \begin{itemize}
    \item Full Degeneration $D_\alpha(z_F)=0$. 
   In this case, since $z_F$ is a singular point for $\dfrac{1}{D_\alpha(z)}$, we assume uniform convergence of its integral on $[0,z_F]$, specifically:
        \begin{equation}\label{eq:intinvD}
            I(D_\alpha)=\int_0^{z_F}{\dfrac{1}{D_\alpha(z)}}<\infty.
        \end{equation}
    \item Quasi-degeneration $D_\alpha(z_F)=\epsilon_F << 1$, in which case 
 $\rho_{\alpha}^{\rm o}$ satisfies the Robin condition  $$\frac{\partial {\rho}_{\alpha}^{\rm o} }{\partial z}(z_{\rm F},t) -  \frac{M_{\alpha}g}{RT} {\rho}_{\alpha}^{\rm o}(z_{\rm F},t) = 0.$$
In this case,  equation \eqref{eq:intinvD} is automatically satisfied as $z_F$ is no longer a singular point to $\dfrac{1}{D_\alpha(z)}$.
\end{itemize}
In this paper, we will adopt full degeneration on $D_\alpha(z)$, for theoretical considerations. \\While for numerical computations, we consider both options.
\item The gas concentration $\rho_\alpha^{atm}$ is assumed to be Lipschitz continuous on $[0,T_e]$.
\end{enumerate}

}
\subsection{Semi-Variational Formulation, Existence and Uniqueness}\label{sec:2}
\noindent Given that  $ 1/D_\alpha^{1/2}\in L^2(0,z_F)$,
i.e.
\begin{equation}\label{qalpha}
I(D_\alpha)= \ltwonorm{ 1/D_\alpha^{1/2}} = \left(\int\limits_0^{z_F} \dfrac{1}{D_\alpha(z)} \;dz \right) ^{1/2} < \infty,
 \end{equation}
 let $H^1(0,z_f)$ be the usual Sobolev subspace. \\However, for the Firn problem, since 
 $D_{\alpha}$ degenerates at $z_F$, we
 use the following subspaces of $H^1(0,z_f)$: 
 $$H^1_\alpha(0,z_F)=\{v\in H^1(0,z_F)\,|\,\Vtnorm{v}<\infty \}$$
\noindent with the following inner product and norm 
\begin{eqnarray}
\Vtinner{v,w} &=& \ltwoinner{D_{\alpha} v_z, w_z}+\ltwoinner{v,w}\\
\Vtnorm{v}^2 &=& \ltwonorm{D_\alpha^{1/2} v_z}^2 + \ltwonorm{v}^2.
\end{eqnarray}
 Note that the injection of $H^1_\alpha$ in $H^1$ is continuous with:
 $$ \Vtnorm{v}^2 \leq q_{\alpha,\infty}\honenorm{v}^2, \mbox{ where }q_{\alpha,\infty}=\max \left\{1, \linftynorm{D_\alpha}\right\}.$$
  Accordingly, we define \begin{equation} H^1_{\alpha,d}(0,z_F) = \{v\in  H^1_{\alpha}(0,z_F) \;|\; v(0) = 0 \}.\end{equation} 
\begin{lemma}\label{hyp:Dalpha}  
$H^1_{\alpha,d}$
 is a closed subspace of $H^1_{\alpha}$ and therefore itself a Hilbert space.
 \end{lemma}
 \begin{proof}
 Let $\{v_n\} \in H^1_{\alpha,d}$ be a converging sequence with $v$ its limit point, and let $\{v_n'\} \in L^2(0,z_F)$ be a uniformly converging sequence.  We need to show that  $v\in H^1_{\alpha,d}$, i.e. $v\in  L^2(0,z_F), v'\in  L^2(0,z_F)$ and $v(0) = 0$. Since $\{v_n\} \in H^1_{\alpha,d}$, then $v_n\in  L^2(0,z_F), v_n'\in  L^2(0,z_F)$ and $v_n(0) = 0$ for all $n$.\\ Moreover, $\lim\limits_{n\rightarrow \infty} v_n = v$ and $v' = (\lim\limits_{n\rightarrow \infty} v_n)' = \lim\limits_{n\rightarrow \infty} v_n'$. 
 Thus $v\in L^2(0,z_F)$ and $v'\in L^2(0,z_F)$. \\It remains to prove that $v(0) = 0$ where $\lim\limits_{n\rightarrow \infty} \Vtnorm{v-v_n} = 0$.
 
 \begin{eqnarray}
 v(z) - v(0) &=&  \int_0^z v'(s) \,ds\nonumber\\
 v_n(z) - v_n(0) &=&  \int_0^z v_n'(s)\, ds\nonumber\\
 -v(0)&=& v_n(z) - v(z)+ \int_0^z v'(s)-v_n'(s) \,ds\nonumber\\
 |v(0)|&\leq & |v_n(z) - v(z)|+ \int_0^z |v'(s)-v_n'(s)| \,ds \nonumber
 \end{eqnarray}
  \begin{eqnarray}
 \mbox{Note that} \int_0^z |v'(s)-v_n'(s)| \,ds &=& \int_0^{z_F}\dfrac{D_\alpha^{1/2}}{D_\alpha^{1/2}} |v'(s)-v_n'(s)| \,ds = \ltwoinner{D_\alpha^{1/2}|v'(z)-v_n'(z)|, 1/D_\alpha^{1/2}}\nonumber\\
 &\leq &I(D_\alpha) \ltwonorm{D_\alpha^{1/2}|v'(z)-v_n'(z)|} \;\leq \; I(D_\alpha)\Vtnorm{v(z)-v_n(z)} \nonumber\\
\therefore |v(0)|  &\leq & |v_n(z) - v(z)|+ I(D_\alpha)\Vtnorm{v(z)-v_n(z)} \label{ineq1}
 \end{eqnarray}
 Integrate \eqref{ineq1} with respect to $z$ from $0$ to $z_F$, then using Cauchy-Schwarz inequality
  \begin{eqnarray}
 z_F\,|v(0)|&\leq & \int_0^{z_F} |v_n(z) - v(z)| dz+ z_FI(D_\alpha)\,\Vtnorm{v(z)-v_n(z)} \nonumber\\
 &\leq& \sqrt{z_F}\ltwonorm{v_n - v}+ z_F\,I(D_\alpha)\,\Vtnorm{v-v_n} 
 \nonumber\\\vspace{3mm}
 \mbox{Hence: }\label{trace1} |v(0)-v_n(0)| &\leq & (1/ \sqrt{z_F}+ I(D_\alpha))\Vtnorm{v(z)-v_n(z)} \\
  \therefore |v(0)| & = & \lim\limits_{n\rightarrow \infty} |v(0)| \;\leq\;   \left(1/\sqrt{z_F}+I(D_\alpha)\right) \lim\limits_{n\rightarrow \infty} \Vtnorm{v-v_n}\;=\;0  \nonumber
 \end{eqnarray}
Thus $v(0)=0$.
 \end{proof}
\noindent In the sequel, we prove a more general estimate in $H^1_\alpha$. Specifically, one obtains the following result.
\begin{lemma}\label{inclusion} Given that
$H^1_\alpha\subset C[0,z_F]$, then $\linftynorm{v}\le (1/\sqrt{z_F}+2I(D_\alpha)) \Vtnorm{v}, \;\;\; {\forall v\in H^1_\alpha}. $
\end{lemma}
 \begin{proof}
 Using  the identity $v(z)=v(0)+\int_0^{z}{v'(x)\,dx}, $ for $z\in [0,z_F]$, then $$|v(z)|\leq |v(0)|+\int_0^{z}{|v'(x)|\,dx}, $$ and therefore on the basis of arguments that lead to \eqref{ineq1} and \eqref{trace1}, one has:
\begin{eqnarray}
\int_0^{z}{|v'(x)|\,dx}&\leq & I(D_\alpha)  \ltwonorm{D_\alpha^{1/2}v'} \leq I(D_\alpha)  \Vtnorm{v} \nonumber\\
|v(0)|& \leq &|v(z)| + \int_0^{z}{|v'(x)|\,dx} \;\leq\; |v(z)| + I(D_\alpha)  \Vtnorm{v}\nonumber\\
\ltwoinner{|v(0)|,\frac{1}{z_F}}\; = \; |v(0)| &\leq &\ltwoinner{|v|,\dfrac{1}{z_F}} + I(D_\alpha)  \Vtnorm{v} \; \leq \; 1/\sqrt{z_F}\ltwonorm{v}+ I(D_\alpha)  \Vtnorm{v} \nonumber\\& \leq &(1/\sqrt{z_F}+I(D_\alpha)) \Vtnorm{v}\nonumber
\end{eqnarray}
Hence,\vspace{-6mm}
\begin{eqnarray}
|v(z)|&\le& |v(0)|+\int_0^{z}{|v'(x)|\,dx}\le \left(1/\sqrt{z_F}+I(D_\alpha)\right) \Vtnorm{v} +I(D_\alpha)  \Vtnorm{v} \nonumber\\ & \le& (1/\sqrt{z_F}+2I(D_\alpha)) \Vtnorm{v}.\nonumber\vspace{-7mm} 
\end{eqnarray}\vspace{-7mm}
 \end{proof}
\subsubsection{Derivation of the Semi-Variational Formulation}
In what follows we denote $ \rho_\alpha^{\rm o} $  by $\rho$. \vspace{2mm}\\
Let $\phi \in H^1_{\alpha,d}(0,z_F)$; multiplying the pde in 
\eqref{eq:trace_gas_dynamics} and using integration by parts with respect to $z$, in addition to the initial and boundary conditions, then equation \eqref{eq1} is reduced to \eqref{eq3}
\begin{eqnarray}
\ltwoinner{[\rho f]_t,\phi} +  \ltwoinner{[f\rho \mathcal{F}]_z,\phi} + \ltwoinner{\rho \mathcal{G} ,\phi} &=& \ltwoinner{\left[D_\alpha \left(\rho_z - {\rho} \mathcal{M}_\alpha \right)\right]_z ,\phi}\label{eq1}\\
f\ltwoinner{\rho_t,\phi} + f\mathcal{F}[\phi\rho ]_0^{z_F} - f\mathcal{F}\ltwoinner{\rho ,\phi_z} + \ltwoinner{\rho \mathcal{G} ,\phi} &=& \left[\phi D_\alpha \left(\rho_z - {\rho} \mathcal{M}_\alpha \right)\right]_0^{z_F}\nonumber\\
&& - 
\ltwoinner{D_\alpha \left(\rho_z - {\rho}\mathcal{M}_\alpha \right),\phi_z}\label{eq2}\\
f\ltwoinner{\rho_t,\phi}  + f\mathcal{F}[\phi\rho ]({z_F})- f\mathcal{F}\ltwoinner{\rho ,\phi_z} + \ltwoinner{\rho \mathcal{G} ,\phi} &=& - 
\ltwoinner{D_\alpha \left(\rho_z - {\rho}\mathcal{M}_\alpha \right),\phi_z}\qquad\label{eq3}
\end{eqnarray}
where $f>0$, $\mathcal{M}_\alpha = \dfrac{M_{\alpha}g}{RT}>0$, $\mathcal{G} =\tau+\lambda > 0$ and $\mathcal{F} = {v}+{w}_{\rm air} > 0$ are constants. Let the bilinear form 
\begin{equation}\label{bilin}
 \mathcal{A}({\rho} ,\phi) = \dfrac{\mathcal{G}}{f} \ltwoinner{\rho  ,\phi}  + 
\dfrac{1}{f}\ltwoinner{D_\alpha \rho_z ,\phi_z} +\mathcal{F} \,\phi({z_F})\,\rho({z_F},t) - \mathcal{F}\ltwoinner{\rho ,\phi_z} - \dfrac{\mathcal{M}_\alpha}{f}\ltwoinner{\rho D_\alpha,\phi_z}  \quad
\end{equation}
 then, \eqref{eq3} becomes
\begin{equation}\label{eq:semi-var}
\ltwoinner{\rho_t,\phi} +  \mathcal{A}({\rho} ,\phi) = 0
\end{equation}
\noindent Seek $\rho: [0, T]\times [0, z_F] \rightarrow \mathbb{R}$ such that: 
\begin{equation}\label{FirnV}
 \begin{cases}
\rho(.,t)  \in H^1_{\alpha,d} +\{\rho^{atm}_\alpha(t)\},& \forall t>0, \\
\ltwoinner{\rho_t,\phi} +  \mathcal{A}({\rho} ,\phi) = 0,& \forall \phi \in H^1_{\alpha,d} , \;\; \forall t>0, \\
 \rho(z,0) = 0. &
 \end{cases}
 \end{equation}
 
 \subsubsection{Existence and Uniqueness to the Semi-Variational Formulation}\label{sec:3}
 To deal with the issue of existence and uniqueness of \eqref{FirnV}, 
 we use Lions theorem (\cite{Brezis}, page 341), which is stated below. 
 Then, we apply it to our problem.
 \begin{theorem}\label{thrm:Brez}
 Let $V$ and $H$ be 2 Hilbert spaces satisfying:
  \begin{equation} \label{VV} V \subset H \subset V^* \mbox{ (the dual of V)},\end{equation}
with the injection from $V$ to $H$ is dense and continuous.\\
 Assuming a bilinear form $a(\cdot, \cdot) : V\times V \rightarrow \mathbb{R}$ satisfies 
 \begin{equation}\label{cond:LM}
 \begin{cases}
 |a(v, w)| \leq M \Vnorm{v} \Vnorm{w} &\\
 |a(v, v)| \geq c \Vnorm{v}^2 - \hat{c} \Hnorm{v}^2
 \end{cases} 
 \end{equation}
 then for $u_0\in H$ and $F(t)\in L^2(0,T;V^*)$, the initial value problem 
\begin{equation}\label{eq:time-depend}
\begin{cases}
\ltwoinner{u_t,v} +a(u(t),v) = <F(t),v>&\\
u(0) = u_0
\end{cases} 
\end{equation}
 admits a unique solution $u$, satisfying:
 \begin{equation}\label{sol:time-depend}
 u\in L^2(0,T;V)\cap C([0,T];H),\,\,\,\,\,\frac{du}{dt}\in L^2(0,T;V^*).
 \end{equation}
 \end{theorem}
\subsubsection*{Application of Lions Theorem to \eqref{FirnV}}\label{sec:app}
 To define the Hilbert spaces $H$ and $V$, we first make a change of variable:
 \begin{equation}\label{eq:Change_variables}
\mbox{Let: }\tilde{\rho}(.,t)=\rho(.,t)- \rho_\alpha^{\rm atm}(t)
 \end{equation}
 Then \eqref{FirnV} becomes:
$$
 \begin{cases}
 \ltwoinner{(\tilde{\rho}+\rho_\alpha^{\rm atm}(t))_t,\phi} +  \mathcal{A}((\tilde{\rho}+\rho_\alpha^{\rm atm}(t)) ,\phi) = 0& \\
 \tilde{\rho}(0) = 0&
 \end{cases}
$$
 i.e.,
 \begin{equation}\label{FirnV-Tr}
 \begin{cases}
 \ltwoinner{\tilde{\rho},\phi} +  \mathcal{A}(\tilde{\rho} ,\phi) = -\ltwoinner{(\rho_\alpha^{\rm atm}(t))_t,\phi}- \mathcal{A}(\rho_\alpha^{\rm atm}(t) ,\phi)& \\
 \tilde{\rho}(0) = 0&
 \end{cases}
 \end{equation}
 with:
 \begin{equation}\label{def:A}
 \begin{cases}
  \mathcal{A}({\tilde{\rho}} ,\phi) = \dfrac{1}{f}\left(\mathcal{G} \ltwoinner{\tilde{\rho}  ,\phi}  + \ltwoinner{D_\alpha \tilde{\rho}_z -\mathcal{M}_\alpha\tilde{\rho} D_\alpha,\phi_z}\right) +\mathcal{F}\left(\phi({z_F})\,\tilde{\rho}({z_F},t) - \ltwoinner{\tilde{\rho} ,\phi_z} \right)  \quad&\\
 \mathcal{A}(\rho_\alpha^{\rm atm}(t) ,\phi)=\rho_\alpha^{\rm atm}(t) \left( \dfrac{\mathcal{G}}{f} \ltwoinner{1  ,\phi}   +\mathcal{F} \,\phi({z_F}) - \mathcal{F}\ltwoinner{1,\phi_z} - \dfrac{\mathcal{M}_\alpha}{f}\ltwoinner{ D_\alpha,\phi_z}\right)  \quad&
\end{cases}
 \end{equation}
 Then, to be in line with Theorem \ref{thrm:Brez}, we let:
 \begin{itemize} 
 \item $u=\tilde\rho$
 \item $F(t,\phi)=-\ltwoinner{(\rho_\alpha^{\rm atm}(t))_t,\phi}- \mathcal{A}(\rho_\alpha^{\rm atm}(t) ,\phi)$, which for every $t$ is a linear form in $\phi$.
 \item $u_0= 0$
\end{itemize}
Then \eqref{FirnV} can be stated as follows:
\begin{equation}\label{eq:Firn_time}
\begin{cases}
\ltwoinner{u_t,\phi} + \mathcal{A}(u(t),\phi) = F(t,\phi)&\\
u(0) = u_0
\end{cases} 
\end{equation}
 Specifically, we let:
$$H=L^2(0,z_F) \mbox{ and }V=H^1_{\alpha,d}(0,z_F). $$
Naturally,  for the above definitions of $H$ and $V$, we have 
\begin{equation}\label{eq:sub}
    V\subset H\subset V^*\;\; (V^*, \mbox{ the dual of V})
\end{equation}
with continuous injection from $V$ into $H$.
\subsubsection*{Validation of \eqref{cond:LM} for the Firn Problem \eqref{eq:Firn_time}} 

 To prove existence and uniqueness to \eqref{eq:Firn_time}, we proceed with a sequence of preliminary results.

\begin{theorem}\label{thrm:bilinear}
    The bilinear form $\mathcal{{A}}$ satisfies the following: 
\begin{enumerate}
\item $\forall v,\phi\in H^1_\alpha:\,\, $%Bi-continuity of $ \mathcal{A}(.,.)$:
 $| \mathcal{A}(v,\phi)|\le C\Vtnorm{v}.\Vtnorm{\phi},$
 \item $\forall v\in  H^1_{\alpha,d}: \,\,$%Weak coercivity of $ \mathcal{A}(.,.)$ on $H^1_{\alpha,d}$:
 $\mathcal{A}(v,v)\ge C_0\Vtnorm{v}^2-C_1\ltwonorm{v}^2,$
\end{enumerate}
where $C$, $C_0$ and $C_1$ are { positive} constants independent of $v$ and $w$.
\end{theorem}
\begin{proof}
Let $$G_f=\dfrac{\mathcal{G}}{f};\,\,f_1= \dfrac{1}{f};\,\,M_{\alpha,f}=\dfrac{\mathcal{M}_\alpha}{f}\vspace{-4mm}$$
 \begin{eqnarray}
\mbox{Then,\;\;} \mathcal{A}({v} ,\phi) &=& G_f\ltwoinner{v ,\phi}  + f_1\ltwoinner{D_\alpha v_z ,\phi_z} +\mathcal{F} \,\phi({z_F})\,v({z_F}) - \mathcal{F} \ltwoinner{v ,\phi_z} - M_{\alpha,f}\ltwoinner{v D_\alpha,\phi_z}  \qquad \nonumber
\end{eqnarray}

\begin{enumerate}
\item We start by checking the bi-continuity of $ \mathcal{A}(.,.)$. 
Given that: \begin{enumerate}
\item $\ltwoinner{v ,\phi} \leq \ltwonorm{v}.\ltwonorm{\phi} \leq  \Vtnorm{v}.\Vtnorm{\phi}$
\item   $\ltwoinner{D_\alpha v,\phi_z} \leq \ltwonorm{D^{1/2}_\alpha v_z}.\ltwonorm{D^{1/2}_\alpha \phi_z}\le \Vtnorm{v}.\Vtnorm{\phi}$
\item  $\phi({z_F})\,v({z_F})\le \linftynorm{\phi}. \linftynorm{v} \leq (1/\sqrt{z_F}+2I(D_\alpha))^2\Vtnorm{v}.\Vtnorm{\phi} $ \;\;(Lemma \ref{inclusion})
\item $|\ltwoinner{v,\phi_z}|\le \linftynorm{v}\left|\ltwoinner{\frac{1}{D^{1/2}_\alpha},D^{1/2}_\alpha\phi_z}\right|\,\le (1/\sqrt{z_F}+2I(D_\alpha))\Vtnorm{v} . \ltwonorm{1/D^{1/2}_\alpha} . \ltwonorm{D^{1/2}_\alpha \phi_z}$

$\hspace{17mm}\le { (1/\sqrt{z_F}+2I(D_\alpha))}I(D_\alpha) \Vtnorm{v}.\Vtnorm{\phi}$ \;\;(Lemma \ref{inclusion})
\item $ |\ltwoinner{v D_\alpha,\phi_z}|\le \linftynorm{D^{1/2}_\alpha}. \ltwonorm{v}. \ltwonorm{D^{1/2}_\alpha\phi_z}\le \linftynorm{D^{1/2}_\alpha}.\Vtnorm{v}.\Vtnorm{\phi}.$
\end{enumerate}
Then,
\begin{equation}\label{id:bicont}
| \mathcal{A}({v} ,\phi) |\le C\Vtnorm{v}.\Vtnorm{\phi}, 
\end{equation}
with $C=G_f+f_1+\mathcal{F}(1/\sqrt{z_F}+2I(D_\alpha))^2+\mathcal{F}{ (1/\sqrt{z_F}+2I(D_\alpha))}I(D_\alpha)+M_{\alpha,f}\linftynorm{D^{1/2}_\alpha}$.\\
This proves the bi-continuity of $ \mathcal{A}(.,.)$.
\item We turn now to the coercivity of $ \mathcal{A}(.,.)$ on $H^1_{\alpha,d}$. Let $v\in H^1_{\alpha,d}$.
\begin{eqnarray}
 \mathcal{A}({v} ,v) &=& G_f\ltwoinner{v  ,v}  + f_1\ltwoinner{D_\alpha v_z ,v_z} +\mathcal{F}\,v({z_F})^2 - \mathcal{F}\ltwoinner{v ,v_z} - M_{\alpha,f}\ltwoinner{D_\alpha v ,v_z}  \qquad \nonumber \\
 &\ge& G_f\ltwonorm{v }^2  + f_1\ltwonorm{D_\alpha^{1/2} v_z}  - M_{\alpha,f}\ltwoinner{D_\alpha v ,v_z}  \qquad \nonumber \\
 &\ge& \min\{G_f,f_1\} \Vtnorm{v}^2   - M_{\alpha,f}\ltwoinner{D_\alpha v ,v_z}  \qquad \nonumber 
\end{eqnarray}
\noindent Moreover,  $|\ltwoinner{v ,D_{\alpha}v_z}| \leq \ltwonorm{D_{\alpha}^{1/2}v}.\ltwonorm{D_{\alpha}^{1/2}v_z}\leq  \linftynorm{D_{\alpha}}^{1/2}.\ltwonorm{D_{\alpha}^{1/2}v_z} .\ltwonorm{v}$, thus 
\begin{equation}\label{ineq:2}
-\ltwoinner{v ,D_{\alpha}v_z}\ge -\linftynorm{D_{\alpha}}^{1/2}.\ltwonorm{D_{\alpha}^{1/2}v_z} .\ltwonorm{v}.
\end{equation}

This implies that for $\Gamma =   M_{\alpha,f} \linftynorm{D_{\alpha}}^{1/2}>0$
\begin{eqnarray}
 \mathcal{A}({v} ,v) &\geq &  \min\{G_f,f_1\}\Vtnorm{v}^2  -\Gamma \ltwonorm{D_\alpha^{1/2}v_z}.\ltwonorm{v} 
\end{eqnarray}

Using the geometric inequality: $ab\le \dfrac{\epsilon}{2}a^2+\dfrac{1}{2\epsilon}b^2 $, for all $\epsilon >0$, then 
\begin{eqnarray}
\ltwonorm{D_\alpha^{1/2}v_z}.\ltwonorm{v} \leq \dfrac{\epsilon}{2} \ltwonorm{D_\alpha^{1/2}v_z}^2 + \dfrac{1}{2\epsilon} \ltwonorm{v}^2 \leq \dfrac{\epsilon}{2} \Vtnorm{v}^2 + \dfrac{1}{2\epsilon} \ltwonorm{v}^2
\end{eqnarray}
 and one obtains:
 \begin{eqnarray}\label{eq:weak_coercive}
 \mathcal{A}({v} ,v) &\geq &  \left[\min\{G_f,f_1\}-\Gamma\dfrac{\epsilon}{2}\right] \Vtnorm{v}^2  -\frac{\Gamma}{2\epsilon} \ltwonorm{v}^2 
\end{eqnarray}
Thus choosing $\epsilon = \dfrac{1}{\Gamma}\min\{G_f,f_1\} >0$ such that: 
\begin{eqnarray}
C_{0}&=&\min\{G_f,f_1\}-\dfrac{\epsilon}{2} \Gamma \;=\; \dfrac{1}{2} \min\{G_f,f_1\}
>0 \nonumber\\
\mbox{and\;\;} C_{1}&=& \frac{1}{2\epsilon} \Gamma \;=\; \dfrac{\Gamma ^2}{2\min\{G_f,f_1\}} 
>0\nonumber
\end{eqnarray}
validates the weak coercivity.\vspace{-7mm}
\end{enumerate}
\end{proof}

\noindent Last point to prove is the existence of a function $f^*{(t)}\in L^2(0,T;V)$, such that:
 $$ F(t,\phi)=-\ltwoinner{(\rho_\alpha^{\rm atm}(t))_t,\phi}- \mathcal{A}(\rho_\alpha^{\rm atm}(t) ,\phi)=\Vtinner{f^*{(t)},\phi},\,\forall \phi\in H^1_\alpha. $$
 Using the bi-continuity of $ \mathcal{A}(.,.)$
 $$| \mathcal{A}(\rho_\alpha^{\rm atm}(t) ,\phi)| \leq C \Vtnorm{\rho_\alpha^{\rm atm}(t)} \Vtnorm{\phi} = C z_F^{1/2} |\rho_\alpha^{\rm atm}(t)| \Vtnorm{\phi}$$ 
 and Cauchy-Schwarz inequality on the inner product $<.,.>$, 
$$|\ltwoinner{(\rho_\alpha^{\rm atm}(t))_t,\phi}| \leq \ltwonorm{(\rho_\alpha^{\rm atm}(t))_t}.\ltwonorm{\phi} \leq z_F^{1/2}|(\rho_\alpha^{\rm atm}(t))_t| \Vtnorm{\phi}$$ 
 one has:
\begin{equation}\label{eqF}\hspace{-14mm}
 |F(t,\phi)|\le z_F^{1/2}\;[\; |(\rho_\alpha^{\rm atm})_t|+C |\rho_\alpha^{\rm atm}|\;]\; \Vtnorm{\phi} \le \hat{C}  \Vtnorm{\phi}, \quad\forall t, \forall \phi \in H^1_\alpha
\end{equation}
where $\hat{C} =  z_F^{1/2} \max\{1,C\}\oneinftynorm{\rho_\alpha^{\rm atm}}>0$ and $ \oneinftynorm{\rho_\alpha^{\rm atm}} = \max\limits_t\; [\,|(\rho_\alpha^{\rm atm})_t|+ |\rho_\alpha^{\rm atm}|\,]$
\newpage
\begin{lemma}\label{lemma:F}
$F(t,\phi)$ is linear and continuous on $H^1_{\alpha}$, i.e. $F(t,\cdot) \in (H^1_\alpha)^* \subset V^*$ for all $t$. 
\end{lemma} 
\begin{proof}
$F(t,\phi)$ is linear in $\phi$ by the linearity of the $L_2$ inner product and the bilinear form $A(\rho ,\phi)$.\\
As for the continuity of $F(t,\phi)$ in  $H^1_{\alpha}$, let $\phi_n \in H^1_{\alpha}$ be a sequence converging to $\phi$, i.e.\\ $\lim\limits_{n\rightarrow \infty} \phi_n = \phi$, then by \eqref{eqF} $$|F(t,\phi_n) - F(t,\phi)| = |F(t,\phi_n - \phi)| \leq \hat{C}\; \Vtnorm{\phi_n-\phi}.$$ Taking the limit as $n$ goes to infinity implies $\lim\limits_{n\rightarrow \infty}\,|F(t,\phi_n) - F(t,\phi)| = 0$. \\Thus, $\lim\limits_{n\rightarrow \infty} F(t,\phi_n) = F(t,\phi)$. 
\end{proof}
\begin{theorem}\label{thrm:f}
 There exists an $f^*(t)\in V$, such that: $F(t,\phi)=\Vtinner{f^*(t),\phi},\,\forall t,\,\forall \phi\in H^1_\alpha.$
\end{theorem}
\begin{proof}
By the Riesz-Fr\`echet representation and Lemma \eqref{lemma:F}, there exists $f^*(t) \in V$ such that $\forall t$ and $\forall \phi \in H^1_{\alpha}$, $$F(t,\phi) = \Vtinner{f^*(t),\phi}.$$
Thus, by isometry,
 $\Vtnorm{f^*(t)} = ||{F(t,\phi)}||_{V^*} 
= \sup\limits_{\phi \in V}\dfrac{|F(t,\phi) |}{\Vtnorm{\phi}}\le \hat{C}$\\
 Then, $f^*(t)\in L^2(0,T;V)$ since $$\int\limits_0^T \Vtnorm{f^*(t)}^2 \,dt \leq T\,\hat{C}^2. \vspace{-7mm}$$
\end{proof}
\noindent We can now state our main existence and uniqueness result.
\begin{theorem}\label{thrm:exUn}
Assuming $D_\alpha\in C[0,z_F]$, and $ 1/D_\alpha^{1/2}\in L^2(0,z_F)$,
then the Firn semi-variation formulation \eqref{eq:Firn_time} admits a unique solution $ \rho \in L^2(0,T;H^1_{\alpha})\cap C([0,T];L^2)$, and $\dfrac{du}{dt}\in L^2(0,T;(H^1_{\alpha})^*).$
\end{theorem}
\begin{proof} By applying Lions theorem to the Firn problem \eqref{eq:Firn_time} with the subspaces 
$H=L^2(0,z_F)$  and $V=H^1_{\alpha,d}(0,z_F)$ that satisfy \eqref{eq:sub}, and using theorems \ref{thrm:bilinear} and \ref{thrm:f} we get the result of this theorem.
\end{proof}
\subsection{Discretization in Space and Time}\label{sec:4}
 We start first by rescaling \eqref{eq:trace_gas_dynamics} to the unit square, then by discretizing the obtained problem in time using Finite Difference Euler-Implicit scheme, followed by space discretization using Finite Element.
 \subsubsection{Rescaling \eqref{eq:trace_gas_dynamics}}
 The Firn extracted data is discrete on some large interval, typically 100 meters, with a distance of multiple meters between the different measurements. Thus, we rescale our PDE in space from the interval $[0,z_F]$ to the interval $[0,1]$ to test our discretization numerically. Moreover, given that we would like to go back several hundred years to approximate the solution,  we also rescale our time interval $[0,T_e]$ to the interval $[0,1]$.\\
\noindent Let $\tilde{t} = t/T_e$, $\tilde{z} = z/z_F$, $\tilde{\rho}\, (\tilde{t} ,\tilde{z} ) = \rho\, (t,z)$, and $\tilde{D}(\tilde{z} ) =D(z)$, where by the chain rule we have: $$\dfrac{\partial  \rho }{\partial t}= \dfrac{1}{T_e}\dfrac{\partial \tilde{\rho}}{\partial \tilde{t}}, \;\;\dfrac{\partial  \rho }{\partial z}= \dfrac{1}{z_F}\dfrac{\partial \tilde{\rho}}{\partial \tilde{z}},
\;\;\dfrac{\partial ^2  \rho }{\partial z^2}= \dfrac{1}{z_F^2}\dfrac{\partial ^2 \tilde{\rho}}{\partial \tilde{z}^2}, \;\; \dfrac{\partial  D }{\partial z}= \dfrac{1}{z_F}\dfrac{\partial \tilde{D}}{\partial \tilde{z}}$$
Then system \eqref{eq:trace_gas_dynamics} becomes for $\tilde{\rho}_\alpha(\tilde{t},\tilde{z})$, $\tilde{z} \in [0,1], t\in [0,1]$:
{
\begin{equation}\label{eq:trace_gas_dynamics_rescaled}
\left\{
\begin{array}{l}
\dfrac{1}{T_e} \dfrac{\partial}{\partial \tilde{t}} [\tilde{\rho}_\alpha f] + \dfrac{1}{z_F} \dfrac{\partial}{\partial \tilde{z}} [\tilde{\rho}_\alpha f ({v}+{w}_{\rm air})] + \tilde{\rho}_\alpha (\tau+\lambda) =\dfrac{1}{z_F} \dfrac{\partial}{\partial \tilde{z}} \left[\tilde{D}_\alpha \left(\dfrac{1}{z_F} \dfrac{\partial \tilde{\rho}_{\rm \alpha}}{\partial \tilde{z}} - \tilde{\rho}_{\alpha} \dfrac{M_{\alpha}g}{RT} \right)\right], \vspace{1mm}\\
\tilde{\rho}_\alpha(0,\tilde{t}) = \tilde{\rho}_\alpha^{\rm \;atm}(\tilde{t}), \;\;\;\;\;\; 0<\tilde{t} \leq 1,\vspace{2mm} \\
\tilde{D}_\alpha({1})\left(\dfrac{1}{z_F}  \dfrac{\partial \tilde{\rho}_{\alpha}}{\partial \tilde{z}}(1,\tilde{t}) -  \dfrac{M_{\alpha}g}{RT} \tilde{\rho}_{\alpha}(1,\tilde{t})\right) = 0,\; \quad 0<\tilde{t} \leq 1.\vspace{2mm}\\
\tilde{\rho}_{\alpha}(\tilde{z},0) = 0
\end{array}
\right.
\end{equation}

\noindent Condition \eqref{eq:intinvD} becomes\vspace{-9mm}
\begin{equation}
           \hspace{15mm} I(D_\alpha)=\int_0^{z_F}{\dfrac{1}{D_\alpha(z)}} \;dz = I(\tilde{D}) = z_F\int_0^{1}{\dfrac{1}{\tilde{D}_\alpha(\tilde{z})}} \;d\tilde{z} <\infty.\nonumber
        \end{equation}}
Replacing $\tilde{\rho}_\alpha$ by $\rho$,  $\tilde{t}$ by $t \in [0,1]$,  $\tilde{z}$ by $z \in [0,1]$, $\tilde{D}$ by $D$, we get the following semi-variational form
\begin{eqnarray}
{ \dfrac{1}{T_e} \ltwoinner{[\rho f]_t,\phi} +  \dfrac{1}{z_F}\ltwoinner{[f\rho \mathcal{F}]_z,\phi} + \ltwoinner{\rho \,\mathcal{G} ,\phi}} &=& {\dfrac{1}{z_F} \ltwoinner{\left[ D_\alpha \left(\dfrac{1}{z_F} \rho_z - {\rho}\mathcal{M}_\alpha  \right)\right]_z ,\phi }\label{eq1-res}}\\
\dfrac{f}{T_e}\ltwoinner{\rho_t,\phi} + \dfrac{f\mathcal{F}}{z_F}[\phi\rho ]_0^{1} - \dfrac{f\mathcal{F}}{z_F}\ltwoinner{\rho ,\phi_z} + \mathcal{G}\,\ltwoinner{\rho  ,\phi} &=& \dfrac{1}{z_F}\left[\phi D_\alpha  \left(\dfrac{1}{z_F} \rho_z - {\rho} \mathcal{M}_\alpha \right)\right]_0^{1}\nonumber\\
&& - 
\dfrac{1}{z_F}\ltwoinner{D_\alpha\left(\dfrac{1}{z_F} \rho_z - {\rho} \mathcal{M}_\alpha \right),\phi_z}\label{eq2-res}\\
\dfrac{f}{T_e}\ltwoinner{\rho_t,\phi}  + \dfrac{f\mathcal{F}}{z_F}[\phi\rho ](1)- \dfrac{f\mathcal{F}}{z_F}\ltwoinner{\rho ,\phi_z} + \mathcal{G}\,\ltwoinner{\rho  ,\phi} &=& - \dfrac{1}{z_F}
\ltwoinner{D_\alpha \left(\dfrac{1}{z_F}\rho_z - {\rho}\mathcal{M}_\alpha \right),\phi_z}\qquad\label{eq3-res}
\end{eqnarray}
where $f>0$, $\mathcal{M}_\alpha = \dfrac{M_{\alpha}g}{RT}\in  [1.9496*10^{-5},6.4826*10^{-4}]$, 
$\mathcal{G} =\tau+\lambda > 0$ and $\mathcal{F} = {v}+{w}_{\rm air} > 0$ are constants. \\Let the bilinear form $\mathcal{A}$ is given by:\\
\begin{equation}\label{bilin-res}
\mathcal{A}({\rho} ,\phi) = \dfrac{T_e\mathcal{G}}{f} \ltwoinner{\rho  ,\phi}  + 
\dfrac{T_e}{z_F^2f}\ltwoinner{D_\alpha \rho_z ,\phi_z} +\dfrac{T_e\mathcal{F}}{z_F} \,\phi(1)\,\rho(1,t) - \dfrac{T_e\mathcal{F}}{z_F}\ltwoinner{\rho ,\phi_z} - \dfrac{T_e\mathcal{M}_\alpha}{z_Ff}\ltwoinner{\rho D_\alpha,\phi_z}  \quad
\end{equation}
 then, \eqref{eq3-res} becomes
\begin{equation}\label{eq:semi-var-res}
\ltwoinner{\rho_t,\phi} + \mathcal{A}({\rho} ,\phi) = 0
\end{equation}
{where 
\begin{equation}I({D}) = z_F\int_0^{1}{\dfrac{1}{{D}_\alpha({z})}} \;d{z} <\infty\end{equation} 
}
\subsubsection{Euler-Implicit Time Discretization}
By integrating  equation (\ref{FirnV}) over the temporal interval $[t,t+\Delta t]$, with $0\leq t\leq 1-\Delta t$, one reaches the following L\textsuperscript{2} Integral Formulation:
 \begin{equation}\label{Firn-tau}
\left\{\begin{array}{ll}
\ltwoinner{\rho(z, t+\Delta t)-\rho(z,t),\phi} =-\int_t^{t+\Delta t}  \mathcal{A}(\rho(z,s), \phi(z)) ds&\\
 \rho(z,0) = 0 &\\
\end{array}\right.\vspace{-2mm}
\end{equation}
For the full discretization of  the Firn equation, 
the term $\int_t^{t+\Delta t} A(\rho(z,s), \phi(z)) ds$ is first discretized using an implicit right rectangular rule:
$$\int_t^{t+\Delta t}  \mathcal{A}(\rho(z,s), \phi(z)) ds={\Delta t} \; \mathcal{A}(\rho(z,t+\Delta t), \phi(z))$$ leading to the following fully implicit scheme in time.
 \begin{equation}\label{Firn-tau23}
\left\{\begin{array}{ll}
\ltwoinner{\rho(z, t+\Delta t)-\rho(z,t),\phi} =-{\Delta t}\;  \mathcal{A}(\rho(z,t+\Delta t), \phi(z))&\\
 \rho(z,0) = 0 &\\
\end{array}\right.\vspace{-2mm}
\end{equation}
 \subsubsection{Finite Element Space Discretization}\label{sec:disc}
 Let 
$\mathcal{N}=\{z_i  \;|\;i=1,2,...,n\}$ be the set of nodes based on the partition of $(0,z_F)$ with $$0=z_1<z_2<...<z_n= 1$$  
and $\mathcal{E} = \{ E_j = [z_j,z_{j+1}] \;|\; j=1,2,..,n-1\}$ the resulting set of elements.\\
The $\mathbb{P}_1$ finite element subspace $X_n$ of $H^1(0,1)$ is given by:
 $$X_n=\{v\in C(0,1)\,|\,v\mbox{ restricted to }E_j\in\mathbb{P}_1,\,j=1,2..,n-1\} \subset H^1_{\alpha}(0,1), $$ 
 Consistently, we define $$X_{n,d} = X_n \cap H^1_{\alpha, d}.$$
 For that purpose, we let ${B}_n=\{\varphi_i|\,i=1,2,...n\}$ be a finite element basis of functions with compact support in $(0,1)$, i.e.,: 
 \begin{equation}\label{FEdef}\forall v_n \in X_n:\,v_n(z)=\sum\limits_{i=1}^n{V_{i}\varphi_i(z)},\;\;\;\; V_i=v_n(z_i),\end{equation} where 
 $\varphi_1(z) = \begin{cases} \dfrac{z_2 - z}{z_2 - z_1}, & z_1\leq z \leq z_2\\
 0,& \mbox{otherwise} \end{cases}$, \qquad $\varphi_n(z) = \begin{cases} \dfrac{z - z_{n-1}}{z_n - z_{n-1}},& z_{n-1}\leq z \leq z_n\\
 0,& \mbox{otherwise} \end{cases}$, \\and  $\varphi_i(z) = 
 \begin{cases} 
 \dfrac{z - z_{i-1}}{z_i - z_{i-1}},& z_{i-1}\leq z \leq z_i\\
 \dfrac{z_{i+1} - z}{z_{i+1} - z_i}, & z_i\leq z \leq z_{i+1}\\
 0,& \mbox{otherwise} \end{cases} $ for $i = 2,.., n-1$.\\\\
We hence obtain the following fully implicit {Computational Model}. \vspace{2mm}\\
Given $\rho_n(t) \in X_{n,d} +\{\rho^{atm}\}$, one seeks: 
 \begin{equation}\label{Firn-tau2}
\left\{\begin{array}{ll}
\rho_n(t+\Delta t)\in  X_{n,d} +\{\rho^{atm}_\alpha(t)\}
&\\
\ltwoinner{\rho_n(t+\Delta t),\phi} +{\Delta t}  \mathcal{A}(\rho_n(t+\Delta t), \phi)= \ltwoinner{\rho_n(t),\phi} ,& \quad\forall \phi \in X_{n,d}, \;\; \forall t>0\\
 \rho_n(0) = 0&\\
\end{array}\right.\vspace{-2mm}
\end{equation} 
where 
\begin{eqnarray}
\rho_n(t) &=& \rho_{n,d}(t) +  \rho^{atm}_\alpha(t) \varphi_1(z),\nonumber\\
\rho_{n,d}(t) &=& \sum\limits_{i=2}^n{\rho(z_{i},t)\varphi_i(z)},\nonumber
\end{eqnarray}
Thus, \eqref{Firn-tau2} is equivalent to 
Given $\rho_{n,d}(t) \in X_{n,d}$, one seeks $\forall t>0$, $\forall \phi \in X_{n,d}$: 
 \begin{equation}\label{Firn-tau33}
\hspace{-3mm}\left\{\begin{array}{lcl}
\rho_{n,d}(t+\Delta t)\in  X_{n,d} 
&&\\
\ltwoinner{\rho_{n,d}(t+\Delta t),\phi} +{\Delta t}  \mathcal{A}(\rho_{n,d}(t+\Delta t),\phi)=\ltwoinner{\rho_{n,d}(t),\phi}-{\Delta t} \mathcal{A}(\rho^{atm}_\alpha(t+\Delta t)\varphi_1, \phi)&&\\
\qquad \hspace{7cm}-\ltwoinner{(\rho^{atm}_\alpha(t+\Delta t)-\rho^{atm}_\alpha(t)) \varphi_1,\phi}&&\\
 \rho_{n,d}(0) = 0.&&
\end{array}\right.\vspace{-2mm}
\end{equation}
\vspace{2mm}

\noindent{ Since the  two discrete systems \eqref{Firn-tau2} and \eqref{Firn-tau33} are equivalent, studying the existence, uniqueness, and stability of one of them is equivalent to that of the other.
\subsubsection{Existence, Uniqueness and Stability to the Discrete System \eqref{Firn-tau33}}
To prove the property of existence on system \eqref{Firn-tau33} and the consequent properties of uniqueness and stability, one introduces the map:
$$\mathcal{L}_{\Delta t}:u=\rho_{n,d}(t)\to w=\rho_{n,d}(t+\tau)=\mathcal{L}_{\Delta t}(u),$$
with $u$ and $w$ satisfying:
\begin{equation}\label{eq:Firn-taud}
\hspace{-3mm}\left\{\begin{array}{lcl}
u, w\in  X_{n,d} 
&&\\
\ltwoinner{w,\phi} +{\Delta t}  \mathcal{A}(w,\phi)= \ltwoinner{u,\phi}-{\Delta t} \mathcal{A}(\rho^{atm}_\alpha(t+\Delta t)\varphi_1, \phi)-\ltwoinner{(\rho^{atm}_\alpha(t+\Delta t)-\rho^{atm}_\alpha(t)) \varphi_1,\phi}&&\\
 \rho_{n,d}(0) = 0.&&
\end{array}\right.
\end{equation}
Define on $X_{n,d}$ the bilinear form:
\begin{equation}\label{eq:Bform}
    B(\Delta t,w, \phi)=\ltwoinner{w,\phi} +{\Delta t}  \mathcal{A}(w,\phi), \quad \forall \phi \in X_{n,d}.\end{equation}
To prove the existence of $w$ from $u$, we investigate the positive definiteness of $B(\Delta t,.,.)$ on $X_{n,d}\times X_{n,d}$. For that purpose, using the weak coercivity of $\mathcal{A}$, as given in \eqref{eq:weak_coercive}, specifically:
$$ \mathcal{A}({v} ,v) \geq   \left[\min\{G_f,f_1\}-\Gamma\dfrac{\epsilon}{2}\right] \Vtnorm{v}^2  -\frac{\Gamma}{2\epsilon} \ltwonorm{v}^2. $$
Define now $Gf_{min}=\min\{G_f,f_1\}$. Letting then $\epsilon=\dfrac{Gf_{min}}{\Gamma}$, yields:
\begin{equation} \label{eq:avv}\mathcal{A}({v} ,v) \geq   \dfrac{1}{2}\ltwonorm{v}^2(Gf_{min}-\dfrac{\Gamma^2}{Gf_{min}})=\dfrac{Gf_{min}^2-\Gamma^2}{2Gf_{min}}\ltwonorm{v}^2.\end{equation}
Letting,
$$K=\dfrac{Gf_{min}^2-\Gamma^2}{2Gf_{min}}=sign(K)|K|,$$
then we can state the following lemma 
\begin{lemma}\label{lem:pd} Under the assumptions of theorem \ref{thrm:exUn},  if $0<\Delta t < \dfrac{\epsilon}{|K|} <1$ for $\epsilon << 1$, one has:
$$B(\Delta t,w,w)\;\;\ge\;\; (1+sign(K)|K|\Delta t)\ltwonorm{w}^2\;\;\ge\;\; c_\epsilon \ltwonorm{w}^2 $$
\end{lemma}
\noindent where:\vspace{-3mm}
\begin{equation}\label{eq:ceps}
  \hspace{-3mm}  c_\epsilon = \begin{cases}
    1,& if\; K\ge 0\\
    1-\epsilon,& if\; K< 0
\end{cases}\;\;>\;\;0
\end{equation}
\begin{proof}
    Let $\phi = w$ in \eqref{eq:Bform}, then by \eqref{eq:avv} we get 
    \begin{eqnarray}
        B(\Delta t,w, w)&=&\ltwonorm{w}^2 +{\Delta t}  \mathcal{A}(w,w) \;\; \ge\;\; (1+sign(K)|K|\Delta t)\ltwonorm{w}^2 \;\; \ge\;\; c_\epsilon \ltwonorm{w}^2 
    \end{eqnarray}
    If $K\geq 0$, then $(1+sign(K)|K|\Delta t) = (1+K\Delta t)\geq 1$.\\ If  $K<0$, then $(1+sign(K)|K|\Delta t) = (1-|K|\Delta t) > (1-\epsilon) > 0$ since $\Delta t < \dfrac{\epsilon}{|K|}$ and $\epsilon<<1$.
\end{proof}
}
\begin{theorem}
Under the assumptions of Lemma \ref{lem:pd} system  \eqref{eq:Firn-taud} has a unique solution, i.e. $w = \mathcal{L}_{\Delta t}(u)$ exists and is unique.
\end{theorem}
\begin{proof}
 For $w \neq 0$,    $B(\Delta t,w,w) > 0$ is positive definite, thus system  \eqref{eq:Firn-taud} has a unique solution.       
\end{proof}
\noindent We look now into the issue of stability for the discrete system \eqref{eq:Firn-taud}, in which we  let $\phi=w$, yielding:
\begin{equation}
    B(\Delta t,w,w)=\ltwoinner{u,w}-{\Delta t} \mathcal{A}(\rho^{atm}_\alpha(t+\Delta t)\varphi_1, w)-\ltwoinner{(\rho^{atm}_\alpha(t+\Delta t)-\rho^{atm}_\alpha(t)) \varphi_1,w},   
\end{equation} 
leading, using Lemma \ref{lem:pd},  to the inequality:
\begin{equation}\label{eq:mainstab}
  c_\epsilon||w||^2_2\le \ltwoinner{u,w}-{\Delta t} \mathcal{A}(\rho^{atm}_\alpha(t+\Delta t)\varphi_1, w)-\ltwoinner{(\rho^{atm}_\alpha(t+\Delta t)-\rho^{atm}_\alpha(t)) \varphi_1,w}. 
\end{equation}
Looking now at the term $\mathcal{A}(\rho^{atm}_\alpha(t+\Delta t)\varphi_1, w)$, one has:
\begin{lemma}\label{lem:discreteA} Assume $\rho^{atm}_\alpha(.)\in C[0,1]$, $D_\alpha\in C[0,1]$. Then, there exists a constant $C_{\mathcal{A}}$, independent of $h$, such that:
  $$|\mathcal{A}(\rho^{atm}_\alpha(t+\Delta t)\varphi_1, w)|\le C_{\mathcal{A}}\,h^{1/2}.$$
\end{lemma}
\begin{proof} Note that:
\begin{eqnarray}
\mathcal{A}(\rho^{atm}_\alpha(t+\Delta t)\varphi_1 ,w) &=& \dfrac{T_e\mathcal{G}}{f} \ltwoinner{\rho^{atm}_\alpha(t+\Delta t)\varphi_1,w}+\dfrac{T_e}{z_F^2f}\ltwoinner{D_\alpha \rho^{atm}_\alpha(t+\Delta t)\varphi_1 ,w_z} \nonumber\\
&& +\dfrac{T_e\mathcal{F}}{z_F} \rho^{atm}_\alpha(t+\Delta t)\varphi_1(1)\,w(1,t+\Delta t)- \dfrac{T_e\mathcal{F}}{z_F}\ltwoinner{\rho^{atm}_\alpha(t+\Delta t)\varphi_1 ,w_z} \nonumber\\ 
&& - \dfrac{T_e\mathcal{M}_\alpha}{z_Ff}\ltwoinner{ D_\alpha\rho^{atm}_\alpha(t+\Delta t)\varphi_1,w_z}.\label{eq:discrhs}
\end{eqnarray}
As $\varphi_1(1)=0$, and writing (for simplicity) $\rho^{atm}:=\rho^{atm}_\alpha(t+\Delta t)$ then \eqref{eq:discrhs} reduces to:
$$\mathcal{A}(\rho^{atm}_\alpha(t+\Delta t)\varphi_1 ,w) =T_e\rho^{atm}\left(\dfrac{\mathcal{G}}{f} \ltwoinner{\varphi_1,w}+ \dfrac{1}{z_F^2f}\ltwoinner{D_\alpha \varphi_1 ,w_z}- \dfrac{\mathcal{F}}{z_F}\ltwoinner{\varphi_1 ,w_z} - \dfrac{\mathcal{M}_\alpha}{z_Ff}\ltwoinner{ D_\alpha\varphi_1,w_z}\right).
$$
Thus:
$$|\mathcal{A}(\rho^{atm}_\alpha(t+\Delta t)\varphi_1 ,w)|\le
||\rho^{atm}_\alpha||_\infty T_e\; C_F (|<\varphi_1,w>|+|<\varphi_1,w_z>|+|<D_\alpha\varphi_1,w_z>|),$$
where $C_F = \max{\left\{\dfrac{\mathcal{G}}{f},\dfrac{1}{z_F^2f},\dfrac{\mathcal{M}_\alpha}{z_Ff}\right\}}$. \\
To handle the right hand sides, for the first term we use  $||\varphi_1||_2=\frac{h^{1/2}}{3}$ and  the Cauchy-Schwarz inequality, to get\vspace{-3mm}
$$|<\varphi_1,w>| \;\;\leq\;\; \ltwonorm{\varphi_1} \ltwonorm{w} \;\;\leq\;\; \dfrac{h^{1/2}}{3} \ltwonorm{w} $$
For the second term, we use $\int_0^h{\varphi_1}=\frac{h}{2}$ in addition to the second mean value theorem:
$$ |<\varphi_1,w_z>|\le ||w_z||_{\infty,(0,h)}\int_0^h{\varphi_1}=\dfrac{h}{2}||w_z||_{\infty,(0,h)}. $$
Using the finite element inverse inequality (cite Ciarlet, page 141, inequality (3.3.32)):
$$||w_z||_{\infty,(0,h)}\le c\,h^{-1/2}||w||_{2,(0,h)}, $$
one obtains:
$$|<\varphi_1,w_z>|\le \dfrac{c}{2}h^{1/2}\ltwonorm{w}.$$
A similar argument leads to:
$$ |<\varphi_1,D_\alpha w_z>| \,\le \,\dfrac{c}{2}h^{1/2}\,\linftynorm{D_\alpha}\,\ltwonorm{w}.$$
Hence, one obtains:
$$|\mathcal{A}(\rho^{atm}_\alpha(t+\Delta t)\varphi_1 ,w)|\le h^{1/2}||\rho^{atm}_\alpha||_\infty T_e\;C_F\,\ltwonorm{w}\,\left(\dfrac{1}{3}+\dfrac{c}{2}(\linftynorm{D_\alpha}+1)\right).$$
Factoring $h^{1/2}||w||_2$ on the right hand side and letting:
$$C_\mathcal{A}=\linftynorm{\rho^{atm}_\alpha} T_e\;C_F\;\left(\dfrac{1}{3}+\dfrac{c}{2}(\linftynorm{D_\alpha}+1)\right), $$
yields the result of the lemma.
\end{proof}
\noindent Hence using the result of this lemma and on the basis of inequality \eqref{eq:mainstab},  one obtains, using Cauchy-Schwarz inequality:
\begin{equation}
  c_\epsilon \ltwonorm{w}^2 \le \ltwonorm{u}\ltwonorm{w}+\Delta t \; C_\mathcal{A}\;h^{1/2}\;\ltwonorm{w}+\linftynorm{\rho^{atm}_\alpha(t+\Delta t)-\rho^{atm}_\alpha(t))} \dfrac{h^{1/2}}{3}\ltwonorm{w}. 
\end{equation}
\noindent This leads us to the following stability theorem.
\begin{theorem}\label{theo:stable} Assume $\rho^{atm}_\alpha(.)\in C[0,1]$ is Lipschitz, i.e., there exists $L_{atm}$, such that:
$$|\rho^{atm}_\alpha(t+\Delta t)-\rho^{atm}_\alpha(t)|\le L_{atm}\Delta t, $$
then, one has:
$$ \ltwonorm{w}\le \dfrac{1}{c_\epsilon} (\ltwonorm{u}+\,h^{1/2}\Delta t (C_\mathcal{A}+L_{atm}/3)).$$     
\end{theorem}
\begin{proof} To reach the result of this theorem, one simplifies $||w||_2$ in \eqref{eq:mainstab}, following it with a division by $c_\epsilon$, then use of the previous lemma, in addition to the Lipschitz property of $\rho^{atm}_\alpha$.
\end{proof}
\noindent As a consequence of this stability result, let $\{\rho_{n,d}\}$ be the sequence generated by the fully discrete scheme \eqref{Firn-tau2}, then one derives the following estimate.
\begin{corollary} Let $C_{\mathcal{A},atm}=C_\mathcal{A}+L_{atm}/3$. Then, under the conditions of the stability theorem \ref{theo:stable}, one has:
\begin{equation}\label{eq:ub}\forall n>0:\, ||\rho_{n,d}||_2\le 
     \begin{cases}
    t_n\,C_{\mathcal{A},atm}\,h^{1/2},& if\; K\ge 0\\
    t_n\,\gamma\,C_{\mathcal{A},atm}\,h^{1/2},\,\gamma\,>1,& if\; K< 0
\end{cases}\;
\end{equation} 
\end{corollary}
\begin{proof}
Starting with \eqref{eq:ceps}, one has with $\rho_{0,d}=0$:
\begin{enumerate}
    \item for $K\ge 0$, $c_\epsilon=1$ and therefore:
    $$||\rho_{n,d}||_2\;\;\le\;\; C_{\mathcal{A},atm}\,h^{1/2}\,n\,\Delta t\;\;=\;\;C_{\mathcal{A},atm}\;h^{1/2}\,t_n .$$
    \item for $K<0$, $c_\epsilon=1-\epsilon$:
    \begin{eqnarray}
       \ltwonorm{\rho_{n,d}}&\le& C_{\mathcal{A},atm}\;h^{1/2}\,n\,\Delta t \,\dfrac{1}{c_\epsilon}\left(1+\dfrac{1}{c_\epsilon}+...+\left(\dfrac{1}{c_\epsilon}\right)^{n-1}\right),\nonumber\\
        &\le& C_{\mathcal{A},atm}\;h^{1/2}\;\Delta t \;\dfrac{1}{\epsilon}\;\left((1/c_\epsilon)^n-1\right).\nonumber
    \end{eqnarray}
    For $\epsilon$ sufficiently small such that $1/c_\epsilon\le (1+\gamma\epsilon)$, with $\gamma \gtrapprox 1$, one has: 
        \begin{eqnarray}\ltwonorm{\rho_{n,d}}&\le& C_{\mathcal{A},atm}\;h^{1/2}\;\dfrac{\Delta t }{\epsilon}\;\left((1/c_\epsilon)^n-1\right)\;\; \le \;\; C_{\mathcal{A},atm}\;h^{1/2}\;\dfrac{\Delta t }{\epsilon}\;\left((1+\gamma\epsilon)^n-1\right)\nonumber \\
        &\lessapprox& C_{\mathcal{A},atm}\;h^{1/2}\;\gamma\,n\,\Delta t\;\;=\;\; C_{\mathcal{A},atm}\;h^{1/2}\;\gamma\,t_n. \nonumber
        \end{eqnarray}
\end{enumerate}
From 1. and 2., follows \eqref{eq:ub}.
\end{proof}
\newpage
\subsection{Algorithm and Testing} \label{sec:testings}
Given the fully discrete system \eqref{Firn-tau33}, we start by extracting its matrix form along with the properties of the matrices and the obtained algorithm in section \ref{sec:mat}. Then, we test the algorithm's robustness, performance and accuracy in sections \ref{sec:rob} and \ref{sec:run}. 
\subsubsection{Matrix Form of the Discrete System}\label{sec:mat}
Let $\phi = \varphi_j$ for $j=2,..,n$ in \eqref{Firn-tau33} and define the vector $\Lambda(t) = [\rho(z_{2},t), \, \rho(z_{3},t), \, \cdots ,\rho(z_{n},t)]^T$ of length $n-1$, then \eqref{Firn-tau33} can be written in Matrix form  
\begin{equation}\label{eq:mat2}
\left\{\begin{array}{lcl}
\left[M+T_e\Delta t\; C_\alpha\right]\Lambda(t+\Delta t) &=&  M\Lambda(t) - T_e\Delta t \;b\\
\Lambda(0)= 0
\end{array}\right.
\end{equation}
where \begin{eqnarray}
    C_\alpha &=& \dfrac{\mathcal{G}}{f} M+\dfrac{1}{z_F^2f}S(D_\alpha)-\dfrac{\mathcal{M}_\alpha}{z_Ff}A(D_\alpha)+\dfrac{1}{z_F}Q \label{eq:C}\\
    T_e\,\Delta t \,b &=&  v_1(t)+T_e\,\Delta t\, v_3(t); \qquad i.e. \quad b = \dfrac{1}{T_e\,\Delta t} v_1(t)+ v_3(t)
\end{eqnarray}
by noting that for $j=2,..,n$:
\begin{itemize}
\item $\ltwoinner{\rho_{n,d}(t+\Delta t),\varphi_j} = \sum\limits_{i=2}^n \rho(z_{i},t+\Delta t)\ltwoinner{ \varphi_i,\varphi_j}$ is equivalent to $M \Lambda(t+\Delta t)$ where $M$ is the $(n-1)\times(n-1)$ Mass matrix whose entries are $M_{i,j} = \ltwoinner{ \varphi_{i+1},\varphi_{j+1}}$ for $i,j = 1,2,\cdots, n-1$.
\item $\ltwoinner{\rho_{n,d}(t+\Delta t),\varphi_j} = \sum\limits_{i=2}^n \rho(z_{i},t)\ltwoinner{\varphi_i,\varphi_j}$ is similarly equivalent to $M \Lambda(t)$.
\item $(\rho^{atm}_\alpha(t+\Delta t)-\rho^{atm}_\alpha(t))\ltwoinner{\varphi_1,\varphi_j}$ is equivalent to the  vector of length $n-1$, \begin{equation} v_1(t) = (\rho^{atm}_\alpha(t+\Delta t)-\rho^{atm}_\alpha(t))\ltwoinner{\varphi_1,\varphi_2} e_1 = (\rho^{atm}_\alpha(t+\Delta t)-\rho^{atm}_\alpha(t))\; \dfrac{z_2}{6} \; e_1
\end{equation} where $e_1 =[1,\, 0\, \cdots , \, 0]^T$.\vspace{1mm}
\item {\color{white}.} \vspace{-9.5mm} \hspace{-4mm} \begin{eqnarray}
\hspace{-4mm}  \mathcal{A}(\rho_{n,d}(t+\Delta t),\varphi_j) &=&\mathcal{A}\left(\sum\limits_{i=2}^n{\rho(z_{i},t+\Delta t)\varphi_i} , \varphi_j\right) \nonumber\\
&=& \dfrac{T_e\mathcal{G}}{f} \sum\limits_{i=2}^n\rho(z_{i},t+\Delta t)\ltwoinner{{\varphi_i} ,\varphi_j}  - \dfrac{T_e}{z_F}\sum\limits_{i=2}^n\mathcal{F} \rho(z_{i},t+\Delta t) \ltwoinner{{\varphi_i} ,\varphi_j'} 
 \nonumber\\
&&+ 
\dfrac{T_e}{z_F^2f} \sum\limits_{i=2}^n\rho(z_{i},t+\Delta t)\ltwoinner{D_\alpha \varphi_i' ,\varphi_j'}
+ \,  \dfrac{T_e\mathcal{F} }{z_F}\,\varphi_j(z_n)\,\rho(z_n,t+\Delta t) \nonumber\\ 
&& - \dfrac{T_e\mathcal{M}_\alpha}{z_Ff}\sum\limits_{i=2}^n\rho(z_{i},t+\Delta t)\ltwoinner{D_\alpha \varphi_i ,\varphi_j'}  
,  \nonumber\end{eqnarray}
is equivalent to
 \begin{eqnarray}&&T_e\left(\frac{\mathcal{G}}{f}M+\frac{1}{z_F^2f}S(D_\alpha)-\dfrac{1}{z_F}K-\frac{\mathcal{M}_\alpha}{z_Ff}A(D_\alpha)\right)\Lambda(t+\Delta t) + \dfrac{T_e}{z_F}v_2(t+\Delta t)\qquad\nonumber\\ 
 &\iff& T_e\left(\frac{\mathcal{G}}{f}M+\frac{1}{z_F^2f}S(D_\alpha)-\dfrac{1}{z_F}K-\frac{\mathcal{M}_\alpha}{z_Ff}A(D_\alpha)+\dfrac{1}{z_F}B\right)\Lambda(t+\Delta t)\nonumber \\
 &\iff& T_e\left(\frac{\mathcal{G}}{f}M+\frac{1}{z_F^2f}S(D_\alpha)-\frac{\mathcal{M}_\alpha}{z_Ff}A(D_\alpha)+\dfrac{1}{z_F}Q\right)\Lambda(t+\Delta t) \label{bilinA}\end{eqnarray} 
 where $v_2(t+\Delta t) = \mathcal{F} \rho({z_n},t+\Delta t) e_{n-1} = B\Lambda(t+\Delta t)$ is an $(n-1) \times 1$ vector of zeros except the last entry, $B$ is an $(n-1)\times(n-1)$ zero matrix with $B(n-1,n-1) = \mathcal{F}$ and $S(D_\alpha),K,A(D_\alpha)$ are $(n-1)\times(n-1)$ matrices whose entries for $i,j = 1,\cdots, n-1$ are respectively \begin{equation}S_{i,j} = \ltwoinner{D_\alpha \varphi_{i+1}' ,\varphi_{j+1}'},\;\;\; K_{i,j} =\mathcal{F} \ltwoinner{ \varphi_{i+1}' ,\varphi_{j+1}},\;\;\; A_{i,j} = \ltwoinner{D_\alpha \varphi_{i+1}' ,\varphi_{j+1}}.\end{equation} Moreover, $Q = B-K$.
\item {\color{white}.}\vspace{-0.85cm} 
\begin{eqnarray}
\mathcal{A}\left(\rho^{atm}_\alpha(t+\Delta t)\varphi_1 , \varphi_j\right) &=& T_e\rho^{atm}_\alpha(t+\Delta t)\left[\dfrac{\mathcal{G}}{f} \ltwoinner{\varphi_1 ,\varphi_j}  + \dfrac{1}{z_F^2f} \ltwoinner{D_\alpha \varphi_1' ,\varphi_j'}- \dfrac{\mathcal{F}}{z_F} \ltwoinner{\varphi_1 ,\varphi_j'}\right. \nonumber\\
&&\left. - \frac{\mathcal{M}_\alpha}{z_Ff}\ltwoinner{D_\alpha \varphi_1 ,\varphi_j'} \right] \nonumber
\end{eqnarray}
is equivalent to the $(n-1) \times 1$  vector $T_ev_3(t)$ where \begin{eqnarray}
 v_3(t) &=& \rho^{atm}_\alpha(t+\Delta t)\left[\dfrac{\mathcal{G}}{f} \ltwoinner{\varphi_1 ,\varphi_2}  + \dfrac{1}{z_F^2f} \ltwoinner{D_\alpha \varphi_1' ,\varphi_2'}-\dfrac{\mathcal{F}}{z_F} \ltwoinner{ \varphi_1 ,\varphi_2'}- \dfrac{\mathcal{M}_\alpha}{z_Ff}\ltwoinner{D_\alpha \varphi_1 ,\varphi_2'} \right]e_1 \nonumber\\
& =&\rho^{atm}_\alpha(t+\Delta t)\, c_1 \, e_1 \label{eq:v3}
 \end{eqnarray}
\end{itemize}

\subsubsection*{Properties of the Matrices}
Assuming a uniform mesh in space, i.e. $h = z_{i+1}-z_i, \forall i=0,\cdots n-1$, then the  matrices  $M,K, Q$ are of size $(n-1)\times (n-1)$ where $\mathcal{F}>0$.

\begin{align*}
M&= \dfrac{h}{6} \begin{pmatrix}
        2 & 1 & 0 & \hdots & 0\\
         1 &4 & 1 & & \vdots\\
         0 & \ddots & \ddots & \ddots & 0\\
         \vdots & & 1 & 4 & 1\\
         0 & \hdots & 0 &1 &2\\
         \end{pmatrix},\hspace{0.35cm}
 K=  \frac{\mathcal{F}}{2}\begin{pmatrix}
         0 & 1 & 0 & \hdots & 0  \\
         -1 & 0 & 1 &  & \vdots\\
         0 & \ddots & \ddots & \ddots & 0\\
         \vdots &  & -1 & 0 & 1\\
         0 & \hdots  & 0 & -1 & 1\\
     \end{pmatrix},\hspace{0.35cm} Q=  \frac{\mathcal{F}}{2}\begin{pmatrix}
      0 & -1 & 0 & \hdots & 0  \\
         1 & 0 & -1 &  & \vdots\\
         0 & \ddots & \ddots & \ddots & 0\\
         \vdots &  & 1 & 0 & -1\\
         0 & \hdots  & 0 & 1 & 1\\
     \end{pmatrix}, 
     \end{align*}
Similarly the matrices $S(D)$ and $A(D)$ are of size $(n-1)\times (n-1)$, where $D$ is a vector of length $n$ with entries $D_i = D_{\alpha}
(z_i), \forall i=1,\cdots, n$.  To compute these matrices we approximate their integrals in $D$ using the Mean Value Theorem or Trapezoidal rule, i.e. 
\begin{eqnarray} S_{i,j} = \ltwoinner{D_\alpha \varphi_{i+1}' ,\varphi_{j+1}'}  &\approx & \dfrac{1}{2} (D_{i+1}+D_{j+1})\ltwoinner{\varphi_{i+1}' ,\varphi_{j+1}'} \label{eq:sijapp}\\
A_{i,j} = \ltwoinner{D_\alpha \varphi_{i+1}' ,\varphi_{j+1}} &\approx & \dfrac{1}{2} (D_{i+1}+D_{j+1})\ltwoinner{\varphi_{i+1}' ,\varphi_{j+1}}\label{eq:aijapp}
\end{eqnarray}

\begin{align*}
     A(D)&\approx { \dfrac{1}{4}\begin{pmatrix}
             D_1-D_3 & D_2+D_3 & 0 & \hdots & 0\vspace{3mm}\\
             -(D_2+D_3) & D_2-D_4 & D_3+D_4 & & \vdots\vspace{3mm}\\
             0 & \ddots & \ddots & \ddots & 0\vspace{3mm}\\
              \vdots& & -(D_{n-2}+D_{n-1}) & D_{n-2}-D_{n} & D_{n-1}+D_{n}\vspace{3mm}\\
             0 & \hdots& 0 & -(D_{n-1}+D_{n}) & D_{n-1}+D_{n} \\ 
         \end{pmatrix}},
         \end{align*}
         \begin{align*}
     S(D)&\approx { \dfrac{1}{2h}\begin{pmatrix}
             {D_1+2D_2+D_3}& -(D_2+D_3) & 0 & \hdots & 0\vspace{3mm}\\
             {-(D_2+D_3)}& {D_2+2D_3+D_4} & {-(D_3+D_4)}& & \vdots\vspace{3mm}\\
             0 & \ddots & \ddots & \ddots & 0\vspace{3mm}\\
             \vdots & & {-(D_{n-2}+D_{n-1})} & {D_{n-2}+2D_{n-1}}+D_n & {-(D_{n-1}+D_n)}\vspace{3mm}\\
             0&\cdots&0&{-(D_{n-1}+D_n)}&D_{n-1}+D_{n}
         \end{pmatrix}}\\
\end{align*}
Similarly, the constant $c_1$ in \eqref{eq:v3} is approximated by 
\begin{equation}\label{eq:c1}
c_1 \approx \tilde{c}_1 = 
\frac{\mathcal{G}}{6f} z_2-\left(\frac{1}{2fz_2}+\frac{\mathcal{M_\alpha}}{4f}\right)
          (D_1+D_2)-\frac{\mathcal{F}}{2}
\end{equation}
The explicit extraction of these five matrices assuming a uniform and nonuniform mesh is detailed in \cite{saramaadthesis}.
In this section we summarize the properties of these matrices.
\begin{lemma}\label{M:spd}
    The mass matrix $M$ is a symmetric positive definite tridiagonal matrix.
\end{lemma}
\begin{proof}The symmetry of $M$ follows from the symmetry of the L2 inner product $(M_{i,j} = M_{j,i})$. The tridiagonal nature follows from the 1D finite element integrals. As for the positive definiteness, 
    let  $v\neq 0$ be a vector of length $n-1$, then \vspace{-4mm}
    \begin{eqnarray}
\dfrac{6}{h}        v^TMv &=& v_1(2v_1+v_2) +\sum\limits_{i=2}^{n-2} v_i(v_{i-1}+4v_i+v_{i+1}) + v_{n-1}(v_{n-2}+2v_{n-1})\nonumber\\
&=& \sum\limits_{i=2}^{n-1} v_iv_{i-1} + \sum\limits_{i=1}^{n-2} v_iv_{i
+1} +2v_1^2+  4 \sum\limits_{i=2}^{n-2} v_i^2+ 2v_{n-1}^2\nonumber\\
&=&2v_1^2+ 2v_1v_2+2v_2^2 + 2\sum\limits_{i=2}^{n-2} v_iv_{i
+1}+ 2 \sum\limits_{i=2}^{n-2} v_i^2 + 2 \sum\limits_{i=3}^{n-1} v_i^2 \nonumber\\
&=&  2 \sum\limits_{i=1}^{n-2}(v_i^2+v_iv_{i
+1}+v_{i+1}^2) =   \sum\limits_{i=1}^{n-2}(v_i+v_{i
+1})^2 +  \sum\limits_{i=1}^{n-2}v_i^2 +  \sum\limits_{i=1}^{n-2}v_{i+1}^2 \;\;>\;\; 0 \;\; \vspace{-25mm} 
    \end{eqnarray}
\end{proof}

\begin{lemma}\label{lemma:27}
    The tridiagonal matrix $K$ and the diagonal matrix $B$ are both positive semi-definite. Moreover, $Q = B-K$ is positive semi-definite.
\end{lemma}
\begin{proof}
 Let  $v\neq 0$ be a vector of length $n-1$, then
  \begin{eqnarray}
\dfrac{2}{\mathcal{F}}        v^TKv &=& v_1v_2 +\sum\limits_{i=2}^{n-2} v_i(-v_{i-1}+v_{i+1}) + v_{n-1}(-v_{n-2}+v_{n-1})\nonumber\\
&=& \sum\limits_{i=1}^{n-2} v_iv_{i+1} - \sum\limits_{i=2}^{n-1} v_iv_{i
-1} +v_{n-1}^2 = v_{n-1}^2  \;\;\geq\;\; 0\\
\dfrac{1}{\mathcal{F}}        v^TBv &=&  v_{n-1}^2  \;\;\geq\;\; 0\\
v^T(B-K)v &=& \mathcal{F}v_{n-1}^2  - \dfrac{\mathcal{F}}{2}v_{n-1}^2  \;=\;  \dfrac{\mathcal{F}}{2}v_{n-1}^2  \geq 0
    \end{eqnarray}
\end{proof}
\begin{lemma}
    Assuming that the vector $D>0$ is strictly decreasing ( $D_i>D_{i+2}$\, for $i = 1, \cdots, n-2$),
    then the tridiagonal matrix $A(D)$ is positive definite. Moreover, $A(D)$ is linear in $D$, specifically $A(cD) = cA(D)$ for $c\in \mathbb{R}$.
\end{lemma}
\begin{proof}
    Let  $v\neq 0$ be a vector of length $n-1$, then
  \begin{eqnarray}
4     v^TAv &=& \sum\limits_{i=2}^{n-2} -v_i(D_i+D_{i+1})v_{i-1} +(D_i-D_{i+2})v_i^2+v_i(D_{i+1}+D_{i+2})v_{i+1}\nonumber \\
&&+(D_1-D_3)v_1^2 +(D_2+D_3)v_1v_2  - v_{n-1}(D_{n-1}+D_n)v_{n-2}+(D_{n-1}+D_n)v_{n-1}^2\nonumber\\
&=& \sum\limits_{i=2}^{n-1} v_i(D_i+D_{i+1})v_{i-1} -\sum\limits_{i=2}^{n-1} v_i(D_i+D_{i+1})v_{i-1} +\sum\limits_{i=1}^{n-2} (D_i-D_{i+2})v_i^2 +(D_{n-1}+D_n)v_{n-1}^2\nonumber\\
&=&\sum\limits_{i=1}^{n-2} (D_i-D_{i+2})v_i^2 +(D_{n-1}+D_n)v_{n-1}^2\;\;>\;\; 0
\end{eqnarray}
 The linearity of the matrix $A(D)$ follows from the linearity its entries, given approximately by \eqref{eq:aijapp}.
\end{proof}
\begin{lemma}
    Assuming the positivity property of the function $D$ \eqref{eq:Dpos}, then the tridiagonal matrix $S(D)$ is symmetric positive definite. Moreover, $S(D)$ is linear in $D$, specifically $S(cD) = cS(D)$ for $c\in \mathbb{R}$.
\end{lemma}
\begin{proof}
    Let  $v\neq 0$ be a vector of length $n-1$, then
  \begin{eqnarray}
2h    v^TSv &=& \sum\limits_{i=2}^{n-2} -v_i(D_i+D_{i+1})v_{i-1} +(D_i+2D_{i+1} +D_{i+2})v_i^2 - v_i(D_{i+1}+D_{i+2})v_{i+1}\nonumber \\
&&+(D_1+2D_2+D_3)v_1^2 -v_1(D_2+D_3)v_2  - v_{n-1}(D_{n-1}+D_n)v_{n-2}+(D_{n-1}+D_n)v_{n-1}^2\nonumber\\
&=& \sum\limits_{i=1}^{n-2}(D_i+2D_{i+1} +D_{i+2})v_i^2  -2\sum\limits_{i=2}^{n-1} v_i(D_i+D_{i+1})v_{i-1}+(D_{n-1}+D_n)v_{n-1}^2\nonumber
\\
&=& \sum\limits_{i=1}^{n-2}(D_i+D_{i+1})v_i^2  -2\sum\limits_{i=2}^{n-1} v_i(D_i+D_{i+1})v_{i-1} +\sum\limits_{i=2}^{n-1}(D_{i} +D_{i+1})v_{i-1}^2+(D_{n-1}+D_n)v_{n-1}^2 \nonumber
\\
&=&\sum\limits_{i=2}^{n-1} (D_i+D_{i+1})(v_i-v_{i-1})^2 + (D_1+D_2)v_1^2\;\;>\;\;0\nonumber
\end{eqnarray}
since all the terms are positive, and at least one entry $v_j\neq 0$, implying that one of the terms is nonzero. The linearity of the matrix $S(D)$ follows from the linearity its entries, given approximately by \eqref{eq:sijapp}.
\end{proof}    
\begin{algorithm}[H]
\centering 
\caption{  The Rescaled Firn Direct Problem  }
{\renewcommand{\arraystretch}{1.3}
\begin{algorithmic}[1]
\Statex{\textbf{Input:} \;\;\; Mass matrix $M$; Matrix $C_\alpha$ as defined in \ref{eq:C}; End time $T_e$; Time step $dt$; Mesh size $h$ }
\Statex{ \qquad \qquad   Discrete initial condition vector $\Lambda_0$; $c_1$ as defined in \eqref{eq:c1}; The function $\rho^{atm}(t)$.}
 \Statex{\textbf{Output:} ${\Lambda}$: $n\times m$ matrix with the computed solution vectors $\Lambda_i$ for $t_i=0,\tau, 2\tau,\cdots,1$ \vspace{3mm}}
 \State $z = 0:h:1$\;; \;\;\; $n = length(z)$;\;\;\; $v = zeros(n-1); \;\;\; \Lambda_{1} = \Lambda_0; \;\;\; \Lambda(:,1) = \Lambda_0;$\vspace{1mm}
 \State $t= 0:dt:1 ; \;\;\; m = length(t); \;\;\; \rho =  \rho^{atm}(t);\;\;\; \Lambda(1,:) = \rho$, \;\;\;$B_\alpha = (M + T_e *dt* C_\alpha)  $\vspace{1mm}
\For {$i= 1:m-1$\vspace{2mm}} 
\State $v(1)= T_e*dt*\rho(i+1)*c_1+(\rho(i+1) - \rho(i))*z(2)/6\,;$\vspace{1mm}
\State  $rhs=M*\Lambda_{i}- v;$ \vspace{1mm}
\State  $ \Lambda_{i+1} = B_\alpha \backslash rhs;$ \qquad $\%$ Solve for $\Lambda_{i+1} $\vspace{1mm}
\State  $\Lambda(2:n,{i+1}) =  \Lambda_{i+1};  $\vspace{2mm}
\EndFor
\end{algorithmic}}
\label{alg:Firn}
\end{algorithm}

The existence and uniqueness of the solution to the discrete problem \eqref{eq:mat2} is proven in the appendix. Thus,
system \eqref{eq:mat2} is solved iteratively given $\Lambda(0)$, as detailed in Algorithm \ref{alg:Firn} which is written in  {\tt MATLAB} syntax.  At each time step a system of linear equations has to be solved using either direct methods, or iterative methods such as Krylov Subspace methods. In our implementation we solve the linear system using {\tt MATLAB}'s backslash operator, which is based on direct solvers.

We consider different end times $T_e = 1, 50, 100, 150$ and firn depths $z_F = 1, 50, 100, 150$. As for the mesh sizes, we consider $h = 1/16, 1/32, 1/64, 1/128, 1/256$ and set the time step $\tau$ to be either $h$ or $h^2$. For testing purposes, we  set the constants:
\begin{itemize}
    \item $f=0.2$,
    \item $\mathcal{M}_\alpha = \dfrac{M_{\alpha}g}{RT} = \dfrac{0.04*9.8}{8.314*260} = 1.8134*10^{-4}$,
    \item $\mathcal{G} =\tau+\lambda=10+0.03 = 10.03$,
 \item $\mathcal{F} = {v}+{w}_{\rm air} = 200+485 = 685$. \end{itemize}
and the functions $\rho^{atm}(t) = 2*(T_e*t)^{1/4}$ for $t\in [0,1]$ and $\overline{\rho}(z) = 0$.\vspace{2mm}\\
We consider two general test cases for $D_\alpha(z)$:
\begin{itemize}
    \item \textbf{ Test Case 1}: A decreasing positive $D_\alpha(z) = 200-199.98*z$ for $z\in [0,1]$, that degenerates near $z=1$, specifically at $z = \frac{200}{199.98}$.
     \item \textbf{ Test Cases 2a, 2b, 2c, 2d}: A decreasing positive $D_\alpha(z) = 200*(1-z)^p$  with $p = 0.25 \,(\mathbf{2a}), $ $0.5 \,(\mathbf{2b}), 0.75\,(\mathbf{2c}), 1\,(\mathbf{2d})$, for $z\in [0,1]$, that degenerates at $z = 1$.
\end{itemize}
 
We start first by testing the robustness of  algorithm \ref{alg:Firn} in section \ref{sec:rob}, 
and then report the runtime of the algorithm in section \ref{sec:run}. It is important to have a numerically stable fast Direct Problem, as it will be called many times from the inverse problem.

\subsubsection{Robustness of the Algorithm}\label{sec:rob}
We test the robustness of the algorithm
by decreasing the mesh size $h$ from $h = 1/2^4$ to $h = 1/2^{8}$ for all the possible combinations of depth $z_F$ and end times $T_e$ that we are considering, with the time step $dt$ set to either $h$ or $h^2$.

We run algorithm \ref{alg:Firn} for test cases 1, 2a, 2b, 2c, and 2d with inputs described above, with $dt$ set to $h^2$, for all the considered $h$ values. Then, we compute the $L$ {infinity} and $L^2$ absolute and relative errors between solutions at time $T_e$ and common space points corresponding to $h=1/2^4$. These errors are computed between solutions for 
 different $h$ values less than $1/2^8$  and the solution for $h=1/2^8$.
 
Table \ref{tab:err1} summarize the obtained results for test case 1, with $T_e = 150$, and $z_F = 1, 50, 100, 150$. 
Table \ref{tab:err12} summarize the obtained results for test case 2b, with $T_e = 150$, and $z_F = 1, 50, 100, 150$. 
The orders of relative errors vary between $10^{-5}$ and $10^{-1}$, whereas the orders of errors vary between $10^{-4}$ and $1$. Moreover, it is observed that in all the cases, the errors decrease with the decrease of $h$, which proves the stability of the algorithm.

We run algorithm \ref{alg:Firn} with same settings described above for all the test cases, but with dt set to h. 
We obtain very similar results to Tables \ref{tab:err1} and \ref{tab:err12} where the corresponding errors for $dt = h$ are slightly smaller with a difference in the sixth
fractional digit in normalized format, i.e. a relative difference of order $10^{-6}$.

\noindent In Tables \ref{tab:err1} and \ref{tab:err12} we only show the results for $T_e = 150$ since for $T_e=1,50,100$  the corresponding relative errors are almost identical up to the fifth fractional digit in normalized format.  
\noindent However, the solution evolves with time as shown in figures \ref{fig:te00}, \ref{fig:te0}, \ref{fig:te}, \ref{fig:te1}, \ref{fig:te50}, \ref{fig:te150} as the concentration increases with time. 

\begin{table}[H]
\setlength{\tabcolsep}{3pt}
\renewcommand{\arraystretch}{1.1}
    \centering
    \begin{tabular}{|c|c|c|c|c|c|}
    \hline
  & & \multicolumn{2}{c|}{$\linftynorm{.}$}& \multicolumn{2}{c|}{$\norm{.}$}\\\hline
       $z_F$& h values &  Error & Relative Error & Error & Relative Error \\ \hline
             \multirow{4}{*}{1} &1/16 & 2.84260562E-01 & 4.06128811E-02 & 7.60384727E-01 & 4.05599277E-02  \\ \cline{2-6}
      & 1/32 & 1.29055000E-01 & 1.84383487E-02 & 3.47982441E-01 & 1.85618440E-02  \\ \cline{2-6}
      &  1/64 & 5.45690913E-02 & 7.79639638E-03 & 1.48751656E-01 & 7.93461018E-03  \\ \cline{2-6}
      &  1/128 & 1.80688112E-02 & 2.58152757E-03 & 4.98705049E-02 & 2.66015874E-03 \\ \hline
        \multirow{4}{*}{50}&      1/16 & 4.99303036E+00 & 7.13364341E-01 & 5.03673206E+00 & 7.15890405E-01  \\ \cline{2-6}
       & 1/32 & 8.32517513E-01 & 1.18943460E-01 & 8.32811541E-01 & 1.18370758E-01  \\ \cline{2-6}
     &  1/64 & 3.02909378E-01 & 4.32772751E-02 & 3.03385945E-01 & 4.31214297E-02  \\ \cline{2-6}
    &    1/128 & 9.35326588E-02 & 1.33632000E-02 & 9.37367578E-02 & 1.33231716E-02 \\ \hline
     \multirow{4}{*}{100}&    1/16 & 4.55610794E+00 & 6.50940351E-01 & 4.94335771E+00 & 7.06231800E-01  \\ \cline{2-6}
       & 1/32 & 5.44911050E-01 & 7.78525433E-02 & 5.44940325E-01 & 7.78527893E-02  \\ \cline{2-6}
      &  1/64 & 2.53331718E-02 & 3.61940147E-03 & 2.53336874E-02 & 3.61929213E-03  \\ \cline{2-6}
      &  1/128 & 1.34323363E-02 & 1.91910503E-03 & 1.34323467E-02 & 1.91900949E-03 \\ \hline
      \multirow{4}{*}{150}  & 1/16 & 3.74078261E+00 & 5.34453174E-01 & 4.28033772E+00 & 6.11540230E-01  \\ \cline{2-6}
       & 1/32 & 1.32873993E+00 & 1.89839760E-01 & 1.33266088E+00 & 1.90399868E-01  \\ \cline{2-6}
      &  1/64 & 5.59041058E-03 & 7.98713260E-04 & 5.59041116E-03 & 7.98712987E-04  \\ \cline{2-6}
      &  1/128 & 3.62355038E-04 & 5.17703968E-05 & 3.62356202E-04 & 5.17705399E-05 \\ \hline
    \end{tabular}
   \caption{The Errors between solutions for test case 1 with different $h$ values and the solution for $h=1/256$, where  $dt = h^2$, and $T_e = 150$.}\label{tab:err1}
   \vspace{-2mm}
\end{table}

\begin{table}[H]
\setlength{\tabcolsep}{3pt}
\renewcommand{\arraystretch}{1.1}
    \centering
    \begin{tabular}{|c|c|c|c|c|c|}
    \hline
  & & \multicolumn{2}{c|}{$\linftynorm{.}$}& \multicolumn{2}{c|}{$\norm{.}$}\\\hline
        $z_F$&h values &  Error & Relative Error & Error & Relative Error \\ \hline
      \multirow{4}{*}{1} & 1/16 & 2.79514729E-01 & 3.99348343E-02 & 7.75773809E-01 & 3.89353916E-02  \\\cline{2-6}
      &  1/32 & 1.27959433E-01 & 1.82818228E-02 & 3.55520964E-01 & 1.78432783E-02  \\ \cline{2-6}
       & 1/64 & 5.43260050E-02 & 7.76166615E-03 & 1.51020736E-01 & 7.57959527E-03  \\ \cline{2-6}
        &1/128 & 1.80244387E-02 & 2.57518799E-03 & 5.01200352E-02 & 2.51547961E-03 \\ \hline
        \multirow{4}{*}{50}    &   1/16 & 4.95768718E+00 & 7.08314789E-01 & 4.98729280E+00 & 7.08605819E-01  \\ \cline{2-6}
        &1/32 & 8.69678028E-01 & 1.24252658E-01 & 8.70281060E-01 & 1.23651498E-01  \\ \cline{2-6}
        &1/64 & 3.14178214E-01 & 4.48872766E-02 & 3.14859826E-01 & 4.47359948E-02  \\ \cline{2-6}
        &1/128 & 9.68525402E-02 & 1.38375182E-02 & 9.71290226E-02 & 1.38003108E-02 \\ \hline
     \multirow{4}{*}{100}      &  1/16 & 4.54330062E+00 & 6.49110544E-01 & 4.88710582E+00 & 6.98190073E-01  \\ \cline{2-6}
       & 1/32 & 4.89274688E-01 & 6.99036637E-02 & 4.89286755E-01 & 6.99013214E-02  \\ \cline{2-6}
      &  1/64 & 3.02779483E-02 & 4.32587167E-03 & 3.02783870E-02 & 4.32568272E-03  \\ \cline{2-6}
      &  1/128 & 1.47907959E-02 & 2.11319091E-03 & 1.47908382E-02 & 2.11307403E-03 \\ \hline
       \multirow{4}{*}{150}  &   1/16 & 3.73515622E+00 & 5.33649320E-01 & 4.23077964E+00 & 6.04459704E-01  \\ \cline{2-6}
     &   1/32 & 1.27872508E+00 & 1.82694037E-01 & 1.28153859E+00 & 1.83095907E-01  \\ \cline{2-6}
      &  1/64 & 5.76114043E-03 & 8.23105779E-04 & 5.76114187E-03 & 8.23105528E-04  \\ \cline{2-6}
     &   1/128 & 5.51343389E-04 & 7.87715446E-05 & 5.51344734E-04 & 7.87716930E-05 \\ \hline
    \end{tabular}
   \caption{The Errors between solutions for test case 2b and  different $h$ values and the solution for $h=1/256$, for  $dt = h^2$, $z_F = 1 $ and $T_e = 150$.}\label{tab:err12}
   \vspace{-2mm}
\end{table}

\begin{minipage}{0.45\textwidth}
\begin{figure}[H]
   % \centering
\hspace{-11mm}\includegraphics[width = 1.2\textwidth]{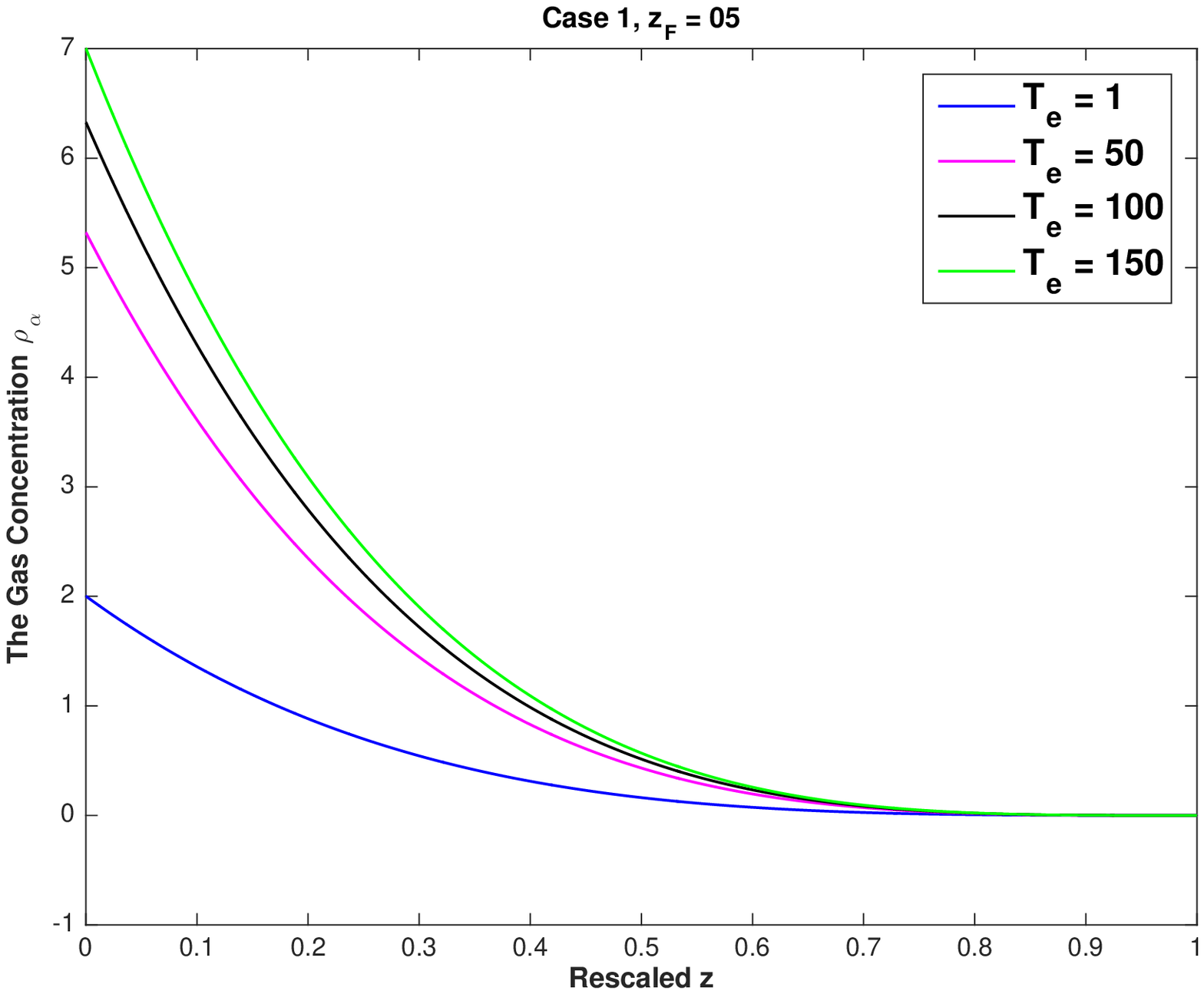}
    \caption{Solution of Rescaled Firn for different $T_e$ with $z_F = 5,$ and $ dt = h= 1/256$}
    \label{fig:te00}
\end{figure}
\end{minipage} \hspace{13mm}
\begin{minipage}{0.45\textwidth}
\begin{figure}[H]
   % \centering
   \hspace{-11mm}\includegraphics[width = 1.2\textwidth]{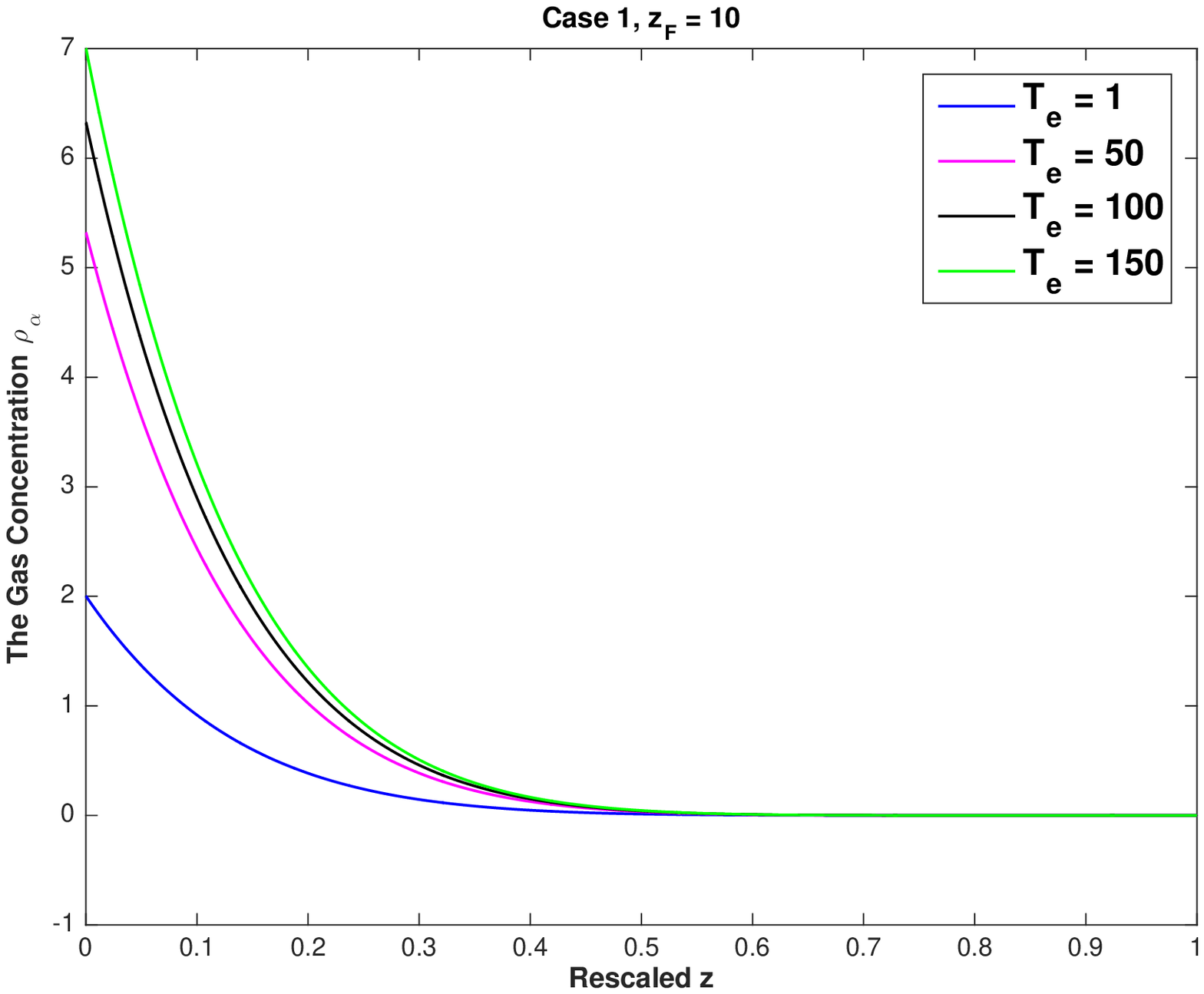}%1.png}
    \caption{Solution of Rescaled Firn for different $T_e$ with $z_F = 10,$ and $ dt = h= 1/256$}
    \label{fig:te0}
\end{figure}
\end{minipage}

\begin{minipage}{0.45\textwidth}
\begin{figure}[H]
   % \centering
\hspace{-11mm}\includegraphics[width = 1.2\textwidth]{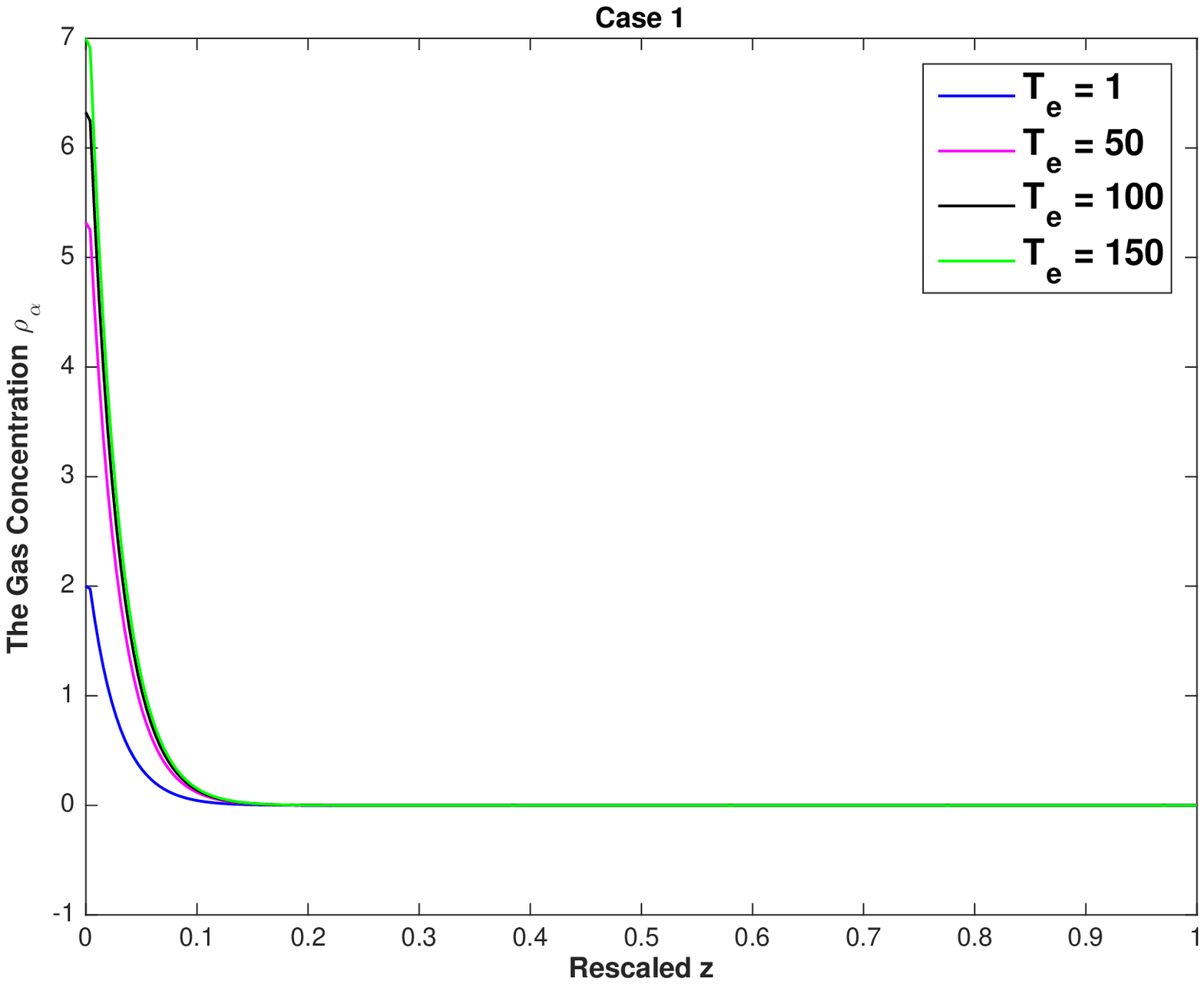}%1.png}
    \caption{Solution of Rescaled Firn for different $T_e$ with $z_F = 50,$ and $ dt = h= 1/256$}
    \label{fig:te}
\end{figure}
\end{minipage} \hspace{13mm}
\begin{minipage}{0.45\textwidth}
\begin{figure}[H]
   % \centering
   \hspace{-11mm}\includegraphics[width = 1.2\textwidth]{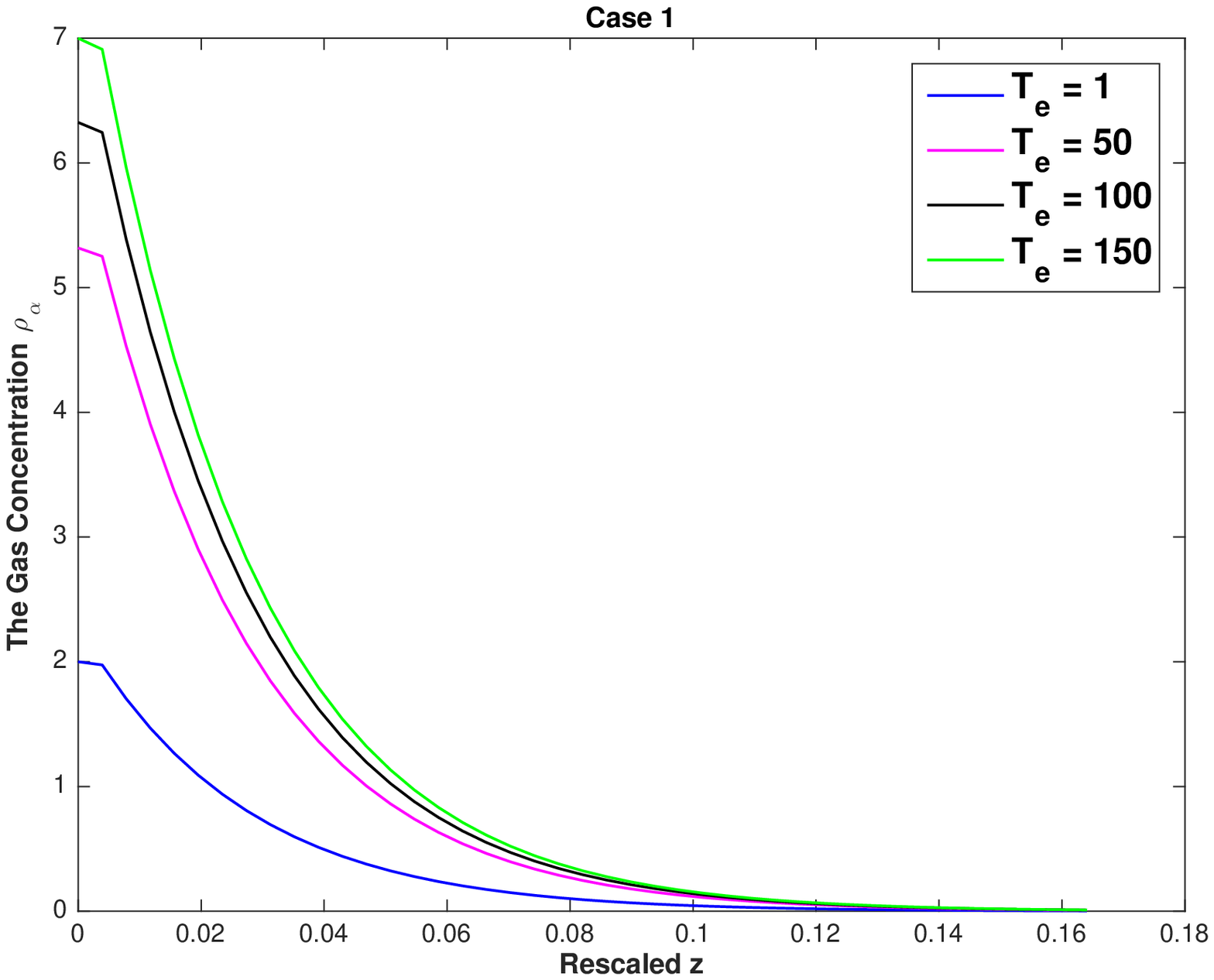}%1.png}
    \caption{Zoomed Solution over first $8m$ of $50m$ for different $T_e$ with $ dt = h= 1/256$} 
    \label{fig:te1}
\end{figure}
\end{minipage}\vspace{-2mm}

\begin{minipage}{0.45\textwidth}
\begin{figure}[H]
   % \centering
\hspace{-11mm}\includegraphics[width = 1.2\textwidth]{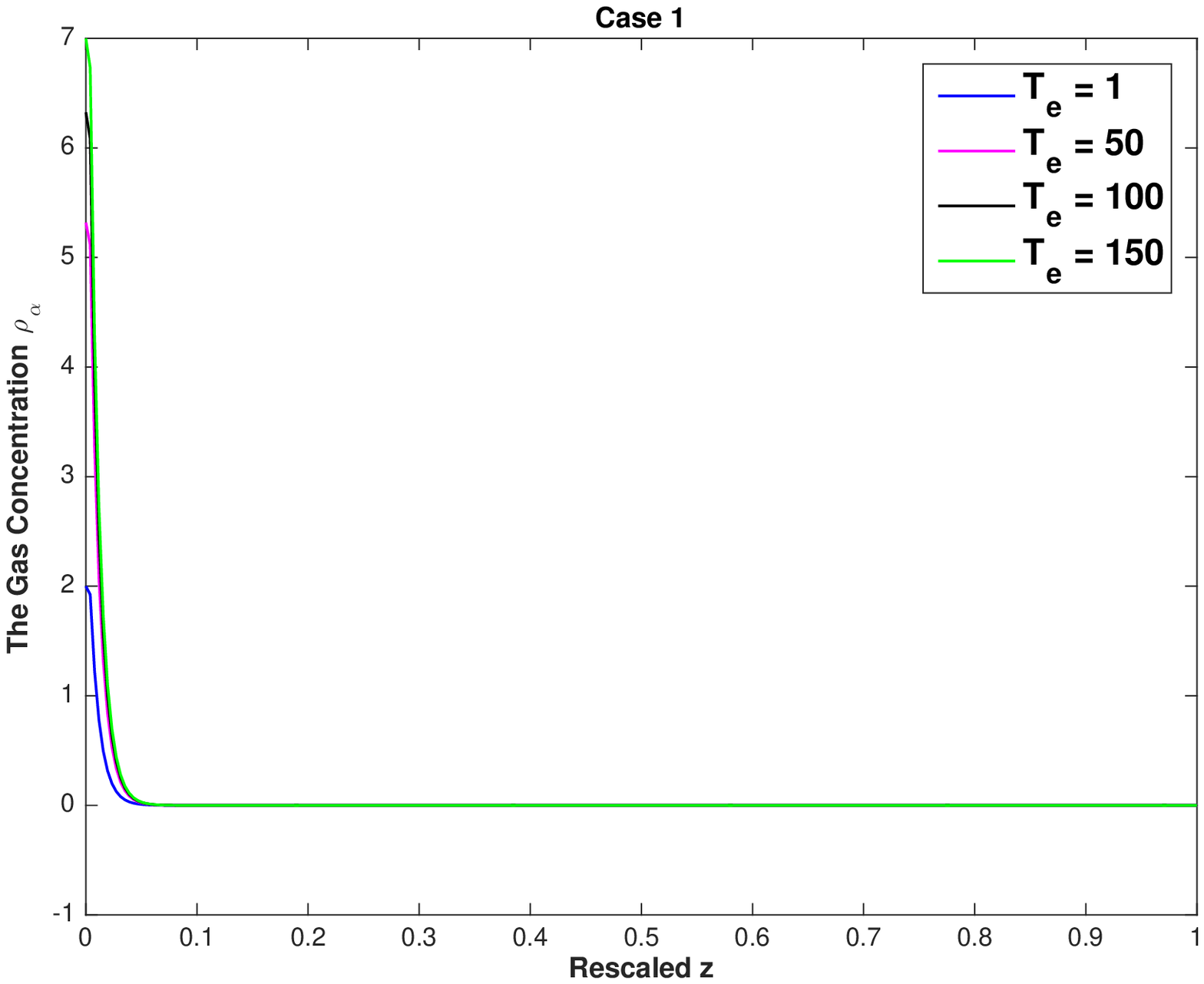}%1.png}
    \caption{Solution of Rescaled Firn for different $T_e$ with $z_F = 150,$ and $ dt = h= 1/256$}
    \label{fig:te50}
\end{figure}
\end{minipage} \hspace{13mm}
\begin{minipage}{0.45\textwidth}
\begin{figure}[H]
   % \centering
\hspace{-11mm}\includegraphics[width = 1.2\textwidth]{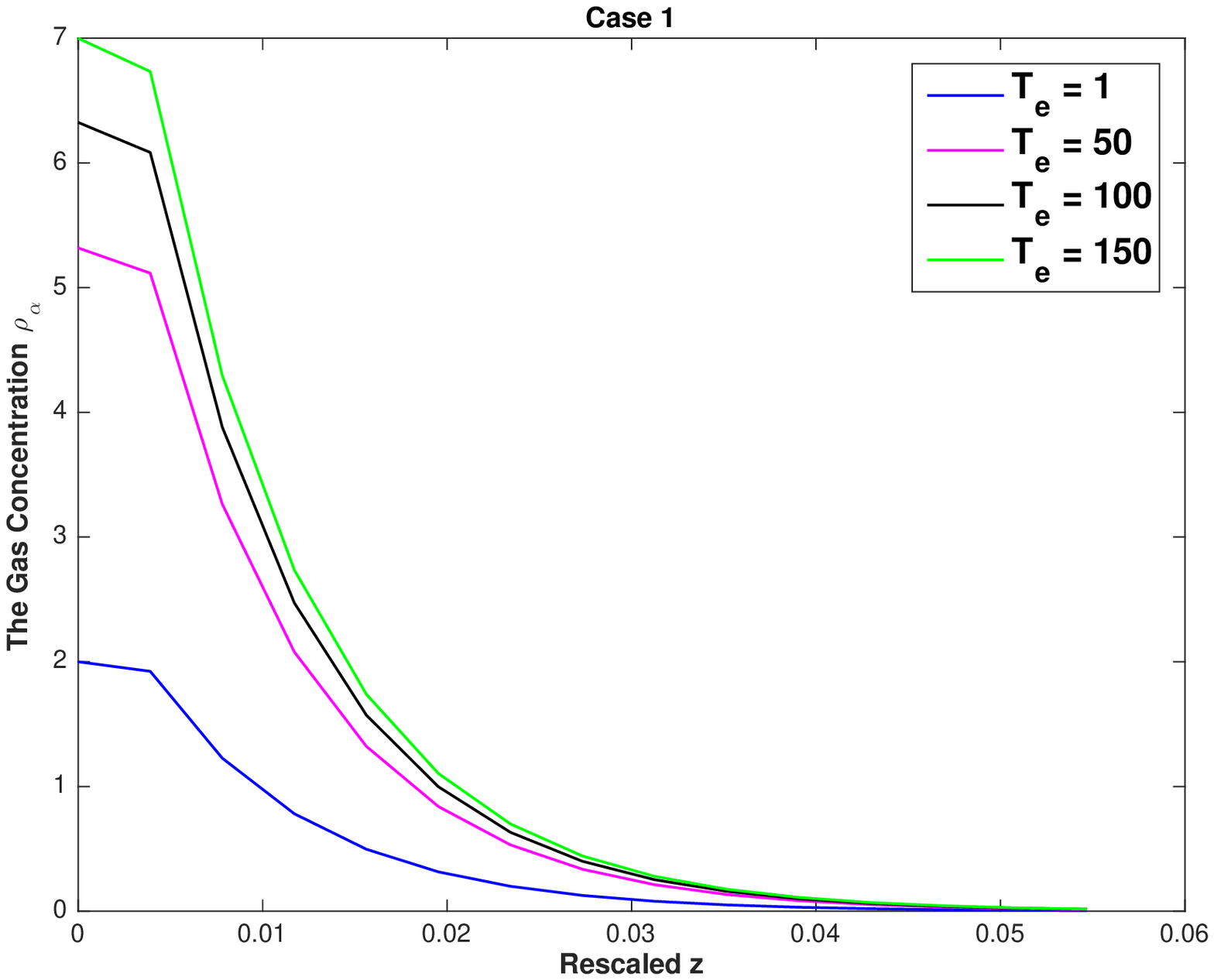}%1.png}
    \caption{Zoomed Solution over first $8m$ of $150m$ for different $T_e$ with $ dt = h= 1/256$} 
    \label{fig:te150}
\end{figure}
\end{minipage} %\vspace{-2mm}

\begin{minipage}{0.45\textwidth}
\begin{figure}[H]
   % \centering
\hspace{-11mm}
\includegraphics[width = 1.2\textwidth]{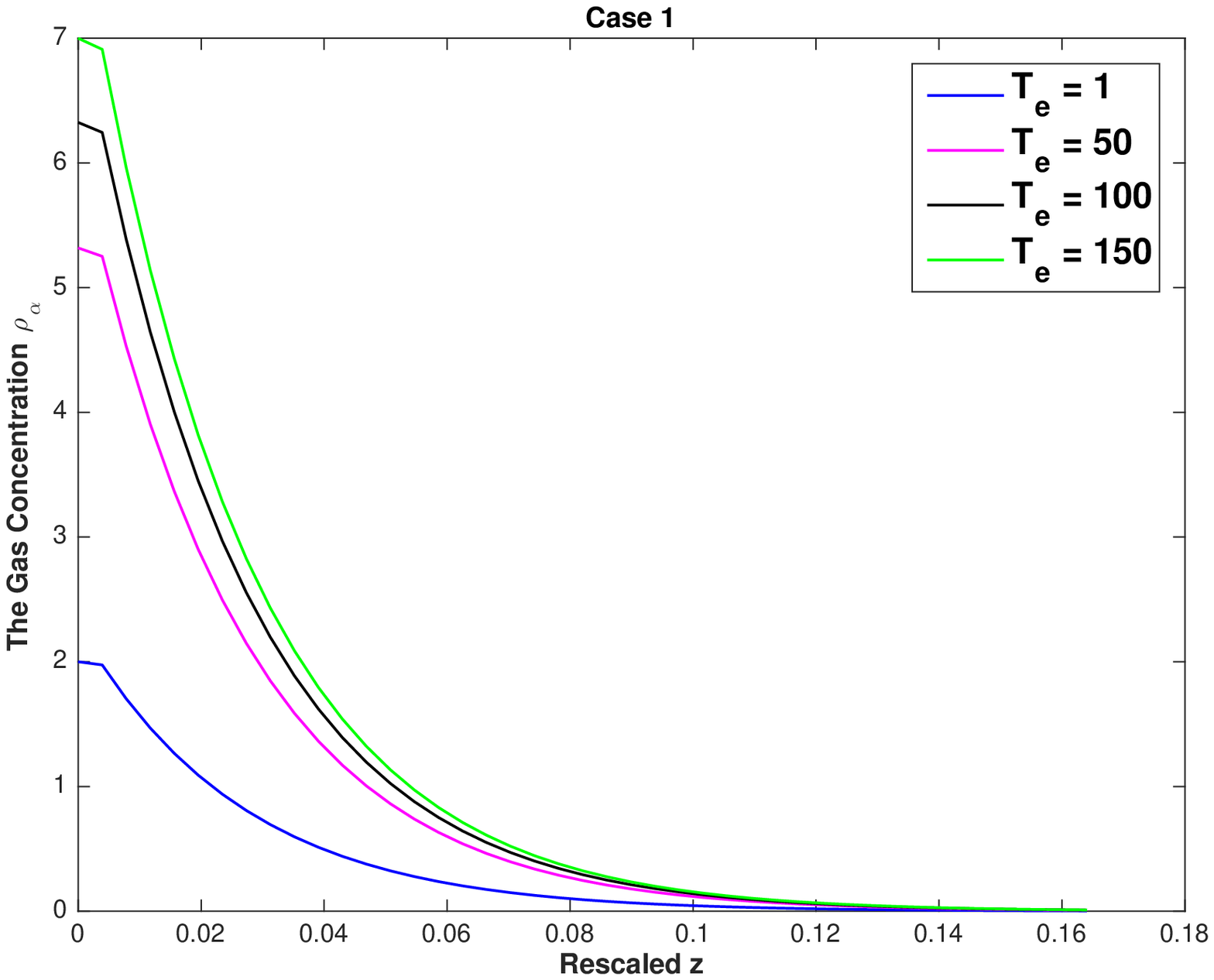}%1.png}
    \caption{Zoomed Solution over first $8m$ of $50m$ for different $T_e$ with $ dt = h= 1/256$} 
    \label{fig:te1-2}
\end{figure}
\end{minipage} \hspace{13mm}
\begin{minipage}{0.45\textwidth}
\begin{figure}[H]
\hspace{-11mm}
\includegraphics[width = 1.2\textwidth]{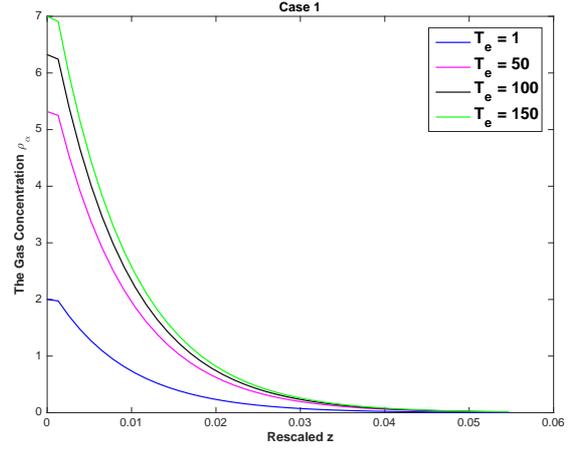}%1.png}
    \caption{Zoomed Solution over first $8m$ of $150m$ for different $T_e$ with $ dt = h= 1/3*256$} 
    \label{fig:te1-3}
\end{figure}
\end{minipage}\vspace{3mm}
\noindent When comparing the obtained concentrations for different $z_F$'s, it is clear from figures \ref{fig:te1} and \ref{fig:te150} that we obtain the same solution up to numerical errors. These errors are due to the fact that even though we are using the same rescaled $h = {1}/{256}$, the actual spacing in the original scale is $h_{z_F} = h*z_F={z_F}/{256}$. Thus, in figure \ref{fig:te1}, the actual $h_{50}
={50}/{256}$, whereas in figure \ref{fig:te150} the actual $h_{150}
={150}/{256} = 3*h_{50}$. So, if we set $h={1}/{3*256}$ for Case 1 with $z_F = 150$ we get the same solution, as shown in figures \ref{fig:te1-2} and \ref{fig:te1-3}.

\subsubsection{Performance and Accuracy}\label{sec:run}
The advantage of using the rescaled problem is that the obtained discrete problem's runtime is independent of the choice of $z_F$ and $T_e$, and solely depends on the choice of the mesh size $h$ and the time step $\Delta t = dt$. Thus, in Table \ref{tab:run1} we show the average runtime in seconds of Algorithm \ref{alg:Firn} in addition to the generation of matrix $C$ and constant $c_1$, for $h=1/16, 1/32, 1/64, 1/128, 1/256$ with $dt = h$ or $h^2$. 

\begin{table}[H]
    \centering
    \begin{tabular}{|c|c|c|c|c|}
    \hline
        ~ & \multicolumn{2}{c|}{Runtime (s)} &  \multicolumn{2}{c|}{Runtime/TimeStep (s)}   \\ \hline
        h values & dt = h$^2$ & dt = h & dt = h$^2$ & dt = h  \\ \hline
        1/16 & 7.400000E-03 & 6.637500E-03 & 2.890625E-05 & 4.148438E-04  \\ \hline
        1/32 & 6.587500E-03 & 2.687500E-04 & 6.433105E-06 & 8.398438E-06  \\ \hline
        1/64 & 3.518750E-02 & 8.125000E-04 & 8.590698E-06 & 1.269531E-05  \\ \hline
        1/128 & 3.144188E-01 & 2.987500E-03 & 1.919060E-05 & 2.333984E-05  \\ \hline
        1/256 & 4.631494E+00 & 1.933125E-02 & 7.067099E-05 & 7.551270E-05 \\ \hline
    \end{tabular}
      \caption{The average runtime and average runtime per time step  of Algorithm \ref{alg:Firn} in seconds.
      }\label{tab:run1}
\end{table}

As expected, setting $dt=h^2$ will slow down the algorithm, as compared to $dt=h$. The reason is not due to a difference in the iterations' complexity, as it is clear that the runtimes per time step are of the same order.
However, the smaller $dt$ is, the more time steps are needed to reach the end Time, implying that more linear systems are solved in total, leading to a slower runtime. 

Thus, in what follows we will be setting $dt = h$ as the algorithm approximates the solution in less runtime. Moreover, the accuracy of the solution is not affected much as shown in figures \ref{fig:16}, \ref{fig:64}, \ref{fig:128} and \ref{fig:256} ,
where we plot the Case 1 obtained $\rho_\alpha(z,T_e)$ for $z_F = T_e = 150$, and $dt = h$ or $h^2$, with $h = 1/16, 1/64, 1/128$ and $h = 1/256$ respectively. The relative error between the solution for $dt = h^2$ and the solution for $dt = h$ is of order $10^{-7}$ for $h = 1/16$, of order $10^{-8}$ for $h = 1/64$, and of order $10^{-9}$ for $h = 1/128$, and $1/256$.

\begin{minipage}{0.45\textwidth}
\begin{figure}[H]
   % \centering
\hspace{-11mm}\includegraphics[width = 1.2\textwidth]{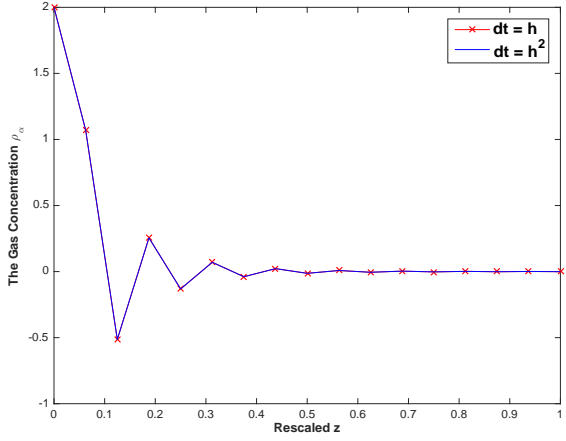}%1.png}
    \caption{The solution $\rho_\alpha(z,T_e)$ for $z_F = T_e = 150$,  $h = 1/16$ , and $dt = h$ or $h^2$, with relative error $= 2.5036*10^{-7}$.}
    \label{fig:16}
\end{figure}
\end{minipage} \hspace{13mm}
\begin{minipage}{0.45\textwidth}
\begin{figure}[H]
   % \centering
\hspace{-11mm}\includegraphics[width = 1.2\textwidth]{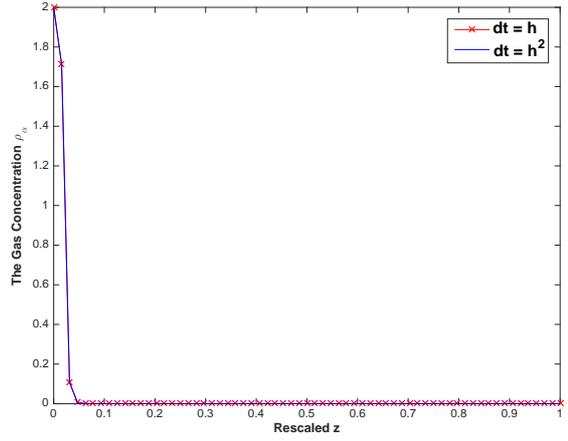}%1.png}
    \caption{The solution $\rho_\alpha(z,T_e)$ for $z_F = T_e = 150$,  $h = 1/64$ , and $dt = h$ or $h^2$, with relative error $=  1.8652*10^{-8}$.} 
    \label{fig:64}
\end{figure}
\end{minipage}

\begin{minipage}{0.45\textwidth}
\begin{figure}[H]
   % \centering
\hspace{-11mm}\includegraphics[width = 1.2\textwidth]{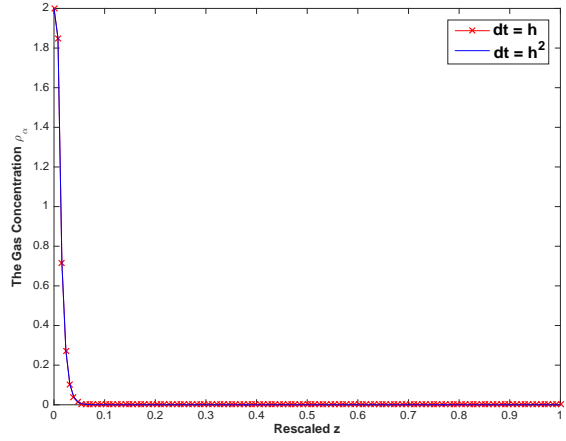}%1.png}
    \caption{The solution $\rho_\alpha(z,T_e)$ for $z_F = T_e = 150$,  $h = 1/128$ , and $dt = h$ or $h^2$, with relative error $= 6.4448*10^{-9}$.}
    \label{fig:128}
\end{figure}
\end{minipage}\hspace{13mm}
\begin{minipage}{0.45\textwidth}
\begin{figure}[H]
   % \centering
\hspace{-11mm}\includegraphics[width = 1.2\textwidth]{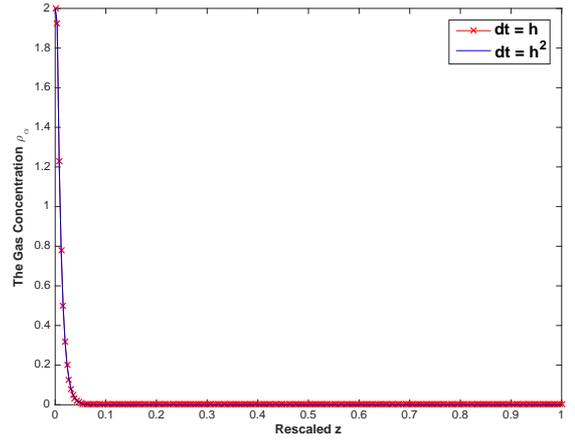}%1.png}
    \caption{The solution $\rho_\alpha(z,T_e)$ for $z_F = T_e = 150$,  $h = 1/256$ , and $dt = h$ or $h^2$, with relative error $= 2.7513*10^{-9}$.} 
    \label{fig:256}
\end{figure}
\end{minipage}\\\\

As for the mesh size $h$, it is obvious that a finer mesh leads to a smoother more accurate solution but at the expense of requiring more runtime, as $dt = h$ will be smaller. On the other hand, if the mesh size $h$, is not small enough with respect to $z_F$ then we may get a numerically wrong solution as shown in figure \ref{fig:16}. It is observed that the largest $h_{max}$ that could be chosen for this problem without obtaining a numerically oscillating solution as in figure \ref{fig:16}, is if $h_{z_F} = h*z_F < 3$. For example, for $z_F = 50$, $h$ has to be less than $3/50 = 0.06$, i.e. $h_{max} = 1/32$; for $z_F = 100$, $h$ has to be less than $3/100 = 0.03$, i.e. $h_{max} = 1/64$; and for $z_F = 150$, $h$ has to be less than $3/150 = 0.02$, i.e. $h_{max} = 1/64$. This explains the large errors observed for $h=1/16,$ and $ 1/32$ in Tables \ref{tab:err1} and \ref{tab:err12}. 
Thus, for $z_F$ ranging between $50$ and $150$, setting the uniform mesh size $h = 1/128$ is a good balance between getting a smooth accurate solution with a fast enough algorithm. 

However, noting that  the solution varies most in the neighborhood of the rescaled $z = 0$ for $z_F$ ranging between $50$ and $150$, then another alternative is using an adaptive mesh size that is finer near $z = 0$, and coarser away from $z = 0$.\\\\
We consider the following adaptive space mesh over the rescaled interval $[0,1]$ for a given $h$:\vspace{2mm}

\begin{minipage}{0.45\textwidth}
    \begin{itemize}
    \item $[0,0.0625]: h/16$
    \item $[0.0625, 0.125]: h/8$
    \item $[0.125,0.25]: h/4$
\end{itemize}
\end{minipage}
\begin{minipage}{0.45\textwidth}
\begin{itemize}
    \item $[0.25,0.5]: h/2$
    \item $[0.5,1]: h$
\end{itemize}
\end{minipage}\vspace{2mm}

\noindent Table \ref{tab:run2} shows the average runtime in seconds of Algorithm \ref{alg:Firn} in addition to the generation of matrix $C$ and constant $c_1$, using the adaptive mesh discussed above for $h=1/4, 1/8, 1/16, 1/32, 1/64$ with $dt = h$ or $h^2$. The total number of mesh points shown in parenthesis is close to that of the uniform mesh with $h=1/16,1/32,1/64,1/128,1/256$. Thus, the corresponding runtimes per time step are comparable in Tables \ref{tab:run1} and \ref{tab:run2}, as the matrices and obtained linear systems are of similar dimensions. However, the total runtimes for the nonuniform mesh are less than that of the uniform mesh since the time step is larger leading to less time iterations, where $dt = h = 1/2^i$ or $dt = h^2 = 1/2^{2i}$ for $i=2:6$ (Table \ref{tab:run2}) versus $i = 4:8$ (Table \ref{tab:run1}). \\

\noindent Similarly to the case of uniform mesh, setting $dt=h$ doesn't affect the solution as shown in figures \ref{fig:4non},\ref{fig:8non},\ref{fig:16non}, and \ref{fig:32non}, where we plot the Case 1 obtained $\rho_\alpha(z,T_e)$ for $z_F = T_e = 150$, with a nonuniform mesh for $h = 1/4, 1/8, 1/16, 1/64$ respectively and $dt = h$ or $h^2$. Moreover, the relative error between the solution for $dt = h$ and that of $dt = h^2$ is of order $10^{-8}$. Thus, using a nonuniform mesh with $h=1/16$ or $1/32$, Algorithm \ref{alg:Firn} provides a smooth and accurate solution comparable to that of a uniform mesh with $h=1/128$ or $1/256$.
\begin{table}[H]
    \centering
    \begin{tabular}{|c|c|c|c|c|}
    \hline
        ~ & \multicolumn{2}{c|}{Runtime (s)} &  \multicolumn{2}{c|}{Runtime/TimeStep (s)}   \\ \hline
        h values ($\#$ points) & dt = h$^2$ & dt = h & dt = h$^2$ & dt = h  \\ \hline
         1/4\, (13)& 1.727000E-04 & 1.168375E-04 & 1.079375E-05 & 2.920938E-05  \\ \hline
        \,1/8 \,(25)& 4.098625E-04 & 1.194625E-04 & 6.404102E-06 & 1.493281E-05  \\ \hline
        1/16 (49)& 2.098913E-03 & 2.470875E-04 & 8.198877E-06 & 1.544297E-05  \\ \hline
        1/32 (97)& 1.879503E-02 & 8.848375E-04 & 1.835452E-05 & 2.765117E-05  \\ \hline
        \;\;1/64 (193) & 2.364509E-01 & 4.129600E-03 & 5.772727E-05 & 6.452500E-05 \\ \hline
    \end{tabular}
      \caption{The average runtime and average runtime per time step in seconds  of Algorithm \ref{alg:Firn} with a nonuniform mesh.
      }\label{tab:run2}
\end{table}

\begin{minipage}{0.45\textwidth}
\begin{figure}[H]
   % \centering
\hspace{-11mm}\includegraphics[width = 1.2\textwidth]{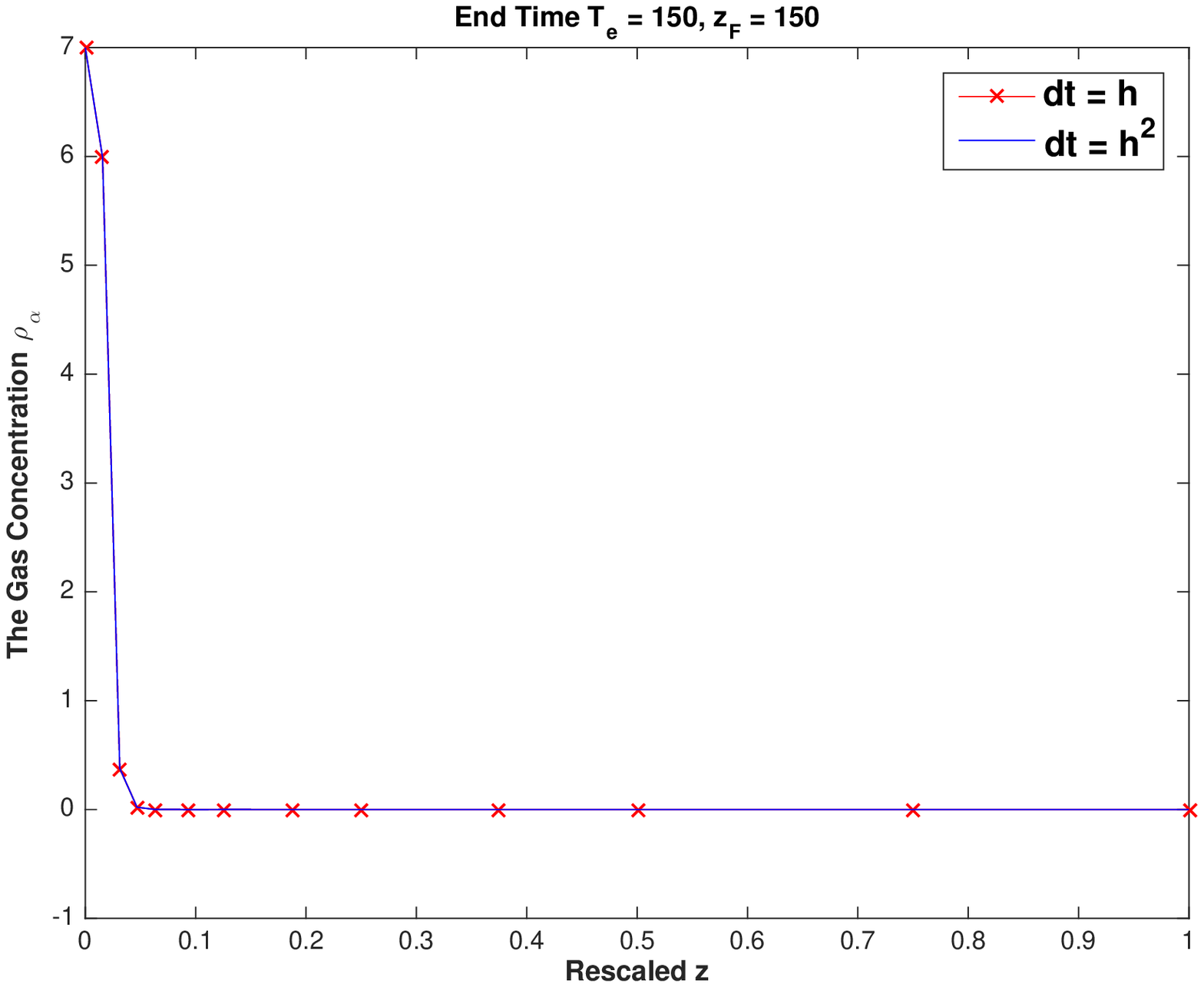}%1.png}
    \caption{The solution $\rho_\alpha(z,T_e)$ for a nonuniform mesh with  $h = 1/4$ , and $dt = h$ or $h^2$, with relative error $= 2.7550*10^{-7}$}
    \label{fig:4non}
\end{figure}
\end{minipage} \hspace{13mm}
\begin{minipage}{0.45\textwidth}
\begin{figure}[H]
   % \centering
\hspace{-11mm}\includegraphics[width = 1.2\textwidth]{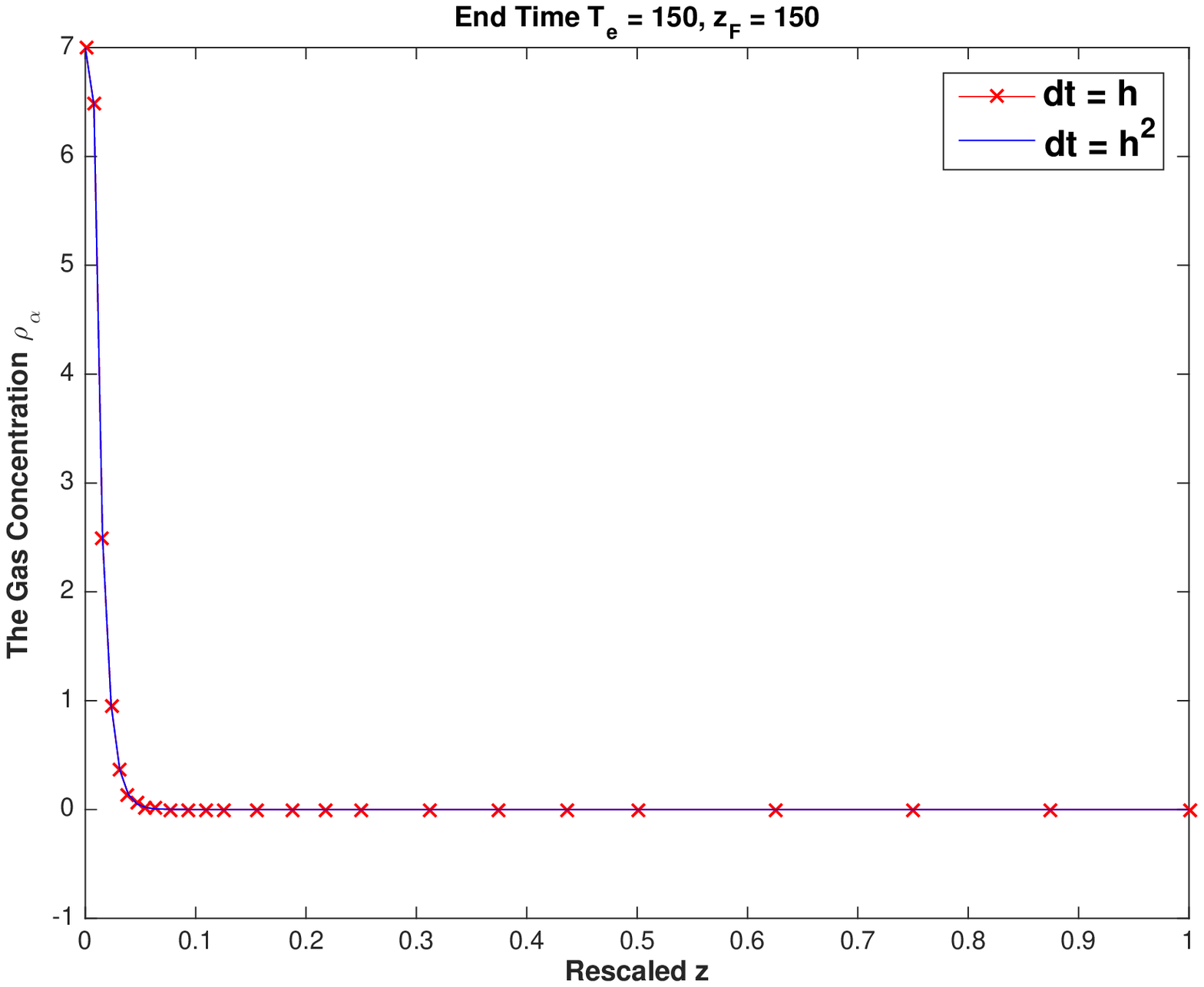}%1.png}
    \caption{The solution $\rho_\alpha(z,T_e)$ for a nonuniform mesh with  $h = 1/8$ , and $dt = h$ or $h^2$, with relative error $= 9.8661*10^{-8}$} 
    \label{fig:8non}
\end{figure}
\end{minipage}

\begin{minipage}{0.45\textwidth}
\begin{figure}[H]
   % \centering
\hspace{-11mm}\includegraphics[width = 1.2\textwidth]{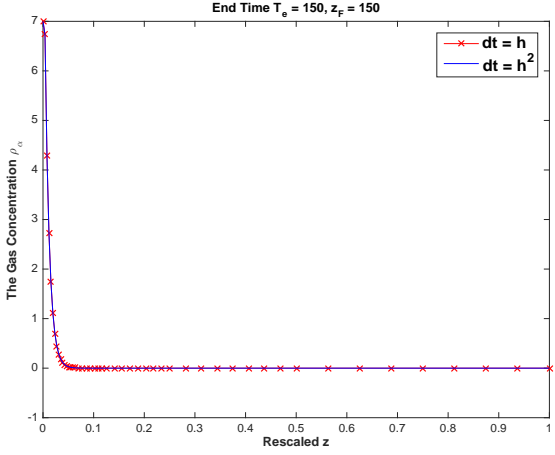}%1.png}
    \caption{The solution $\rho_\alpha(z,T_e)$ for a nonuniform mesh with  $h = 1/16$ , and $dt = h$ or $h^2$, with relative error $= 4.3011*10^{-8}$}
    \label{fig:16non}
\end{figure}
\end{minipage}\hspace{13mm}
\begin{minipage}{0.45\textwidth}
\begin{figure}[H]
   % \centering
\hspace{-11mm}\includegraphics[width = 1.2\textwidth]{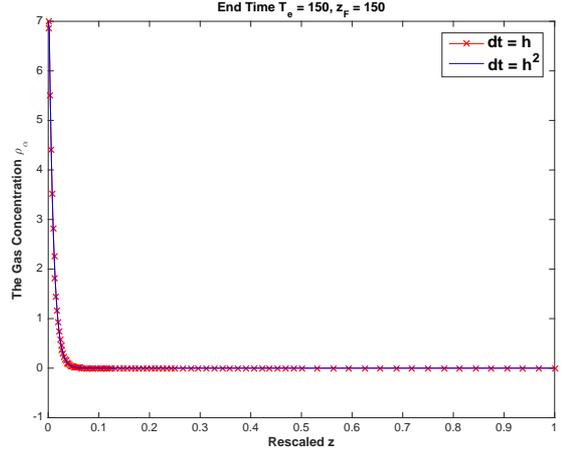}%1.png}
    \caption{The solution $\rho_\alpha(z,T_e)$ for a nonuniform mesh with  $h = 1/32$ , and $dt = h$ or $h^2$, with relative error $= 2.0969*10^{-8}$} 
    \label{fig:32non}
\end{figure}
\end{minipage}

 %%%%%%%%%%%%%%%%%%%%%%%%%%%%%%%%%%%%%%%%%%%%%%%%%%%%%%%%
\section{Inverse Problem}\label{sec:inv}
%%%%%%%%%%%%%%%%%%%%%%%%%%%%%%%%%%%%%%%%%%%%%%%%%%%%%%%%%
Given the rescaled direct problem \eqref{eq:trace_gas_dynamics_rescaled} and its matrix form \eqref{eq:mat2} where $z := \tilde{z} \in [0,1]$, and $t:=\tilde{t}\in [0,1]$, the inverse problem consists of finding the diffusion coefficients $D_\alpha(z)$ of different gases $\alpha$ , given the data measurements $ \rho_\alpha^{\rm o}$ at end time $T_e$ for $z\in (0, z _ {\rm F}),$ which is rescaled to $ \rho_\alpha^{\rm o}(z,1)$ for $z \in [0,1]$. 

Without loss of generality, we will assume that the diffusion coefficients $D_\alpha(z)$ are given by the expression \eqref{TT2} where $ r_\alpha $ are known constants for all gases $\alpha$.
\begin{eqnarray} \label{TT2}
D_\alpha(z) = r_\alpha D_{\textrm{CO2, air}}(z).
\end{eqnarray}
Thus, it is sufficient to find $D_{CO2,air}$, which can then be used to obtain all other $D_\alpha$'s by \eqref{TT2}. Consequently, the inverse problem is to find  $D_{CO2,air}$ by minimizing the objective function 
\begin{equation}\label{inv2}
V( D_{CO2,air}) =  \sum\limits_{\alpha} ||    \rho_\alpha(z) - g_{\alpha}(z)||_2^2 = \sum\limits_{\alpha} || \rho_{\alpha,c}^{\rm o}(z,1) - \rho_{\alpha,g}^{\rm o}(z,1)||_2^2 
\end{equation}
 where $g_{\alpha}(z)= \rho_{\alpha,g}^{\rm o}(z,1)$ is the given rescaled data measurement, and  $\rho_\alpha(z)  = \rho_{\alpha,c}^{\rm o}(z,1)$ is computed using the direct problem \eqref{eq:mat2} for the corresponding approximation of $D_\alpha(z)$.\\\\
 We will solve the constrained minimization problem \eqref{eq:min} and the unconstrained version
 \begin{equation}\label{eq:min}
 \begin{cases} V(D_{CO2,air}) = \min\limits_{d\,
 \geq\,0} V(  d) &\\
 subject \; to \; P( \rho_\alpha, d) = 0\;\; and \; boundary\;conditions &\\
 \end{cases}
 \end{equation}
 where the continuous differential operator of  \eqref{eq:trace_gas_dynamics_rescaled} is \begin{eqnarray}
  P( \rho_\alpha, d) &=& 
  \dfrac{\partial {\rho}_\alpha}{\partial {t}} + \dfrac{T_e\,\mathcal{F}}{z_F} \dfrac{\partial {\rho}_\alpha}{\partial {z}} +  \dfrac{T_e\,\mathcal{G}}{f} {\rho}_\alpha  - \dfrac{T_e}{f\,z_F} \dfrac{\partial}{\partial {z}} \left[r_\alpha\, d \left(\dfrac{1}{z_F} \dfrac{\partial {\rho}_{\rm \alpha}}{\partial {z}} - \mathcal{M}_\alpha {\rho}_{\alpha} \right)\right].
 \end{eqnarray}
There are different type of methods for solving constrained and unconstrained minimization problems, that could be generally categorized as gradient-free methods and gradient methods. The  gradient-free methods require only the objective function evaluation. 
Whereas gradient methods for solving a  minimization problem require at least the objective function and gradient evaluation,  such as the steepest descent, and the nonlinear Conjugate gradient method. 
Moreover, there are other methods that require the gradient and the hessian, such as SQP method, interior point method, and Newton's method for solving $\nabla V=0$. Note that a {\tt MATLAB} implementation of SQP and interior point methods within the ``fmincon" function in the optimization toolbox, allows solving the minimization problem without providing the gradient and the hessian, where the first is approximated using finite differences and the hessian using finite differences or a positive definite quasi-Newton approximation using the BFGS method. However, if the gradient expression or some approximation of it is available, then it could be passed to fmincon.\\ 

 \noindent The computation of the gradient and the hessian may be time consuming, if possible or available. Thus, we will be considering methods that will require the gradient evaluation, and the hessian is approximated, if needed.  For example, the gradient-based methods, such as steepest descent or nonlinear CG methods start with an initial guess $d_0$, compute $d_{n+1} = d_n +\delta$ where $\delta$'s computation depends on that of  
 the gradient of V. One option is to approximate the gradient using finite differences as in ``fmincon". We also consider another option where we get an approximation starting from the directional derivative of $V(d)$ along the direction $\beta$, denoted by $\mathcal{D}_\beta V$. \\

\noindent The directional derivative of $V(\rho_\alpha, d)$ with respect to $d$ is given by 
\begin{eqnarray}
\mathcal{D}_{\beta} V(d)  \;\;=\;\; \mathcal{D} V \, \cdot \beta \;\;=\;\; \nabla V \,.\, \beta 
&=& 2\sum\limits_{\alpha} \ltwoinner{ \rho_\alpha -  g_{\alpha}, \mathcal{D}_{\beta}\rho_\alpha }\nonumber \\
&=&  2\sum\limits_{\alpha} \ltwoinner{ \rho_\alpha -  g_{\alpha}, \mathcal{D}\rho_\alpha\,\cdot\, \beta} \label{eq:direc}
\end{eqnarray}
where $\mathcal{D}\rho_\alpha$ and $\mathcal{D}V$  denote the Fr\'echet derivatives of $\rho_\alpha$ and $V$ with respect to $ d$, and 
$$v_{\alpha,\beta} \;\;=\;\; \mathcal{D}\rho_\alpha\,\cdot\, \beta \;\; =\;\; \mathcal{D}_{\beta}\rho_\alpha \;\;=\;\; \lim\limits_{\epsilon \rightarrow 0} \dfrac{\rho_\alpha(d+\epsilon\beta)-\rho_\alpha(d)}{\epsilon}$$ denotes the directional derivative of $\rho_\alpha$ along $\beta(z)$. In what follows, we refer to $v_{\alpha,\beta}$ by $v_{\alpha}$.\\
\noindent To define $v_{\alpha}$, we start by differentiating the PDE $P(\rho_\alpha,d)=0$ along the direction $\beta(z)$ to obtain
\begin{equation}
{
\frac{\partial v_\alpha }{\partial t} + \dfrac{T_e\, \,\mathcal{F}}{z_F} \frac{\partial v_\alpha}{\partial z} + \dfrac{T_e\,\mathcal{G}}{f} \;v_\alpha  =\dfrac
{T_e}{f\,z_F}\frac{\partial}{\partial z} \left[r_\alpha \,d \left(\dfrac
{1}{z_F}\frac{\partial v_{\rm \alpha}}{\partial z} - \mathcal{M}_\alpha {v}_{\alpha} \right)\right]} +\dfrac
{T_e}{f\,z_F} \frac{\partial}{\partial z} \left[r_\alpha\mathcal{D}_{\beta}\{d\} \left(\dfrac
{1}{z_F}\frac{\partial \rho_{\rm \alpha}}{\partial z} - \mathcal{M}_\alpha {\rho}_{\alpha} \right)\right] 
\end{equation}
where $\mathcal{D}_{\beta}\{d\}$, the directional derivative of $d$ along $\beta$ is 
\begin{equation}
    \mathcal{D}_{\beta}\{d\} = \mathcal{D}\{d\}\cdot \beta = \beta 
\end{equation}
Thus, $v_\alpha$ is the solution of the following Initial Boundary Value Problem for $z\in (0,1)$
\begin{equation}\label{eq:va}
\hspace{-3mm}\left\{
\begin{array}{l}
\dfrac{\partial v_\alpha }{\partial t} + \dfrac{T_e \,\mathcal{F}}{z_F} \dfrac{\partial v_\alpha}{\partial z} + \dfrac{T_e\,\mathcal{G}}{f} \;v_\alpha  =\dfrac{T_e}{f\,z_F} \dfrac{\partial}{\partial z} \left[r_\alpha\, d \left(\dfrac{1}{z_F}\dfrac{\partial v_{\rm \alpha}}{\partial z} - \mathcal{M}_\alpha {v}_{\alpha} \right)\right]
+ \dfrac{T_e}{f\,z_F}\dfrac{\partial}{\partial z} \left[r_{\alpha} \; \beta  \left(\dfrac{1}{z_F}\dfrac{\partial \rho_{\rm \alpha}}{\partial z} - \mathcal{M}_\alpha {\rho}_{\alpha} \right)\right]  \vspace{1mm}\\
v_\alpha(0,t) = 0, \; t>0,\vspace{2mm} \\
\dfrac{1}{z_F}\displaystyle{ \frac{\partial {v}_{\alpha}}{\partial z}(1,t) -  \mathcal{M}_\alpha {v}_{\alpha}(1,t) = 0},
\; t>0.\vspace{2mm}\\
v_\alpha(z,0) = 0,
\end{array}
\right.
\end{equation}
In section \ref{sec:gradV} we define the procedure that gives the gradient of the objective function using the directional derivatives and optimize this procedure in terms of runtime as much as possible. Then, we test the obtained Algorithm \ref{alg:grad} using different constrained and unconstrained optimization methods in section \ref{sec:test}.

\subsection{Gradient of Objective Function }\label{sec:gradV}
 In section \ref{sec:discv} we discretize the directional derivative \eqref{eq:va} and then approximate $\nabla V$ in section \ref{sec:gradv}.
\subsubsection{Directional Derivative Approximation} \label{sec:discv}
 To solve \eqref{eq:va} we first find its equivalent semi-variational form, then use Euler-Implicit discretization in time and Finite Element discretization in space.
 
 Let $\mathcal{T} = \{ \phi \in H^1(0,1)\;|\; \phi(0) = 0 \}$ then \eqref{eq:va} is given in variational form by \eqref{eq1:va}. Using integration by parts with respect to $z$, in addition to the initial and boundary conditions, and $\phi(0) = 0$,  then equation \eqref{eq1:va} is reduced to \eqref{eq3:va}
 \begin{eqnarray}
\hspace{-5mm}0\hspace{-1mm}&=&\hspace{-1mm}
\ltwoinner{[v_\alpha]_t,\phi} +  \dfrac{T_e\, \mathcal{F}}{z_F}\ltwoinner{[v_\alpha]_z,\phi} + \dfrac{T_e\,\mathcal{G}}{f}\ltwoinner{v_\alpha ,\phi} -  \dfrac{T_e}{zf\,_F}\ltwoinner{\left[r_\alpha \,d\left(\dfrac{1}{z_F}[v_\alpha]_z - {v_\alpha} \mathcal{M}_\alpha \right)\right]_z ,\phi}\nonumber\\
&&\hspace{-1mm} - \dfrac{T_e}{f\,z_F}\ltwoinner{\left[r_\alpha \,\beta\left(\dfrac{1}{z_F}[\rho_\alpha]_z - {\rho_\alpha} \mathcal{M}_\alpha \right)\right]_z ,\phi}\label{eq1:va}\\
&=&\hspace{-1mm}
\ltwoinner{[v_\alpha]_t,\phi} + \dfrac{T_e\, \mathcal{F}}{z_F}[\phi\, v_\alpha ]_0^{1} - \dfrac{T_e\, \mathcal{F}}{z_F}\ltwoinner{v_\alpha ,\phi_z} +\dfrac{T_e\,\mathcal{G}}{f}\ltwoinner{v_\alpha  ,\phi} - \dfrac{T_e}{fz_F}\left[\phi \,d\,r_\alpha \left(\dfrac{1}{z_F}[v_\alpha]_z - {v_\alpha} \mathcal{M}_\alpha \right)\right]_0^{1}\nonumber\\
&&\hspace{-1mm} +\dfrac{T_e}{fz_F}\ltwoinner{r_\alpha\,d \left(\dfrac{1}{z_F}[v_\alpha]_z - v_\alpha\mathcal{M}_\alpha \right),\phi_z} - \dfrac{T_e}{fz_F}\left[\phi \,\beta\,r_\alpha \left(\dfrac{1}{z_F}[\rho_\alpha]_z - {\rho_\alpha} \mathcal{M}_\alpha \right)\right]_0^{1}\nonumber\\
&&\hspace{-1mm}+\dfrac{T_e}{fz_F}\ltwoinner{r_\alpha\,\beta \left(\dfrac{1}{z_F}[\rho_\alpha]_z - {\rho_\alpha}\mathcal{M}_\alpha \right),\phi_z}\nonumber\\
&=&\hspace{-1mm}
\ltwoinner{[v_\alpha]_t,\phi} + \dfrac{T_e \mathcal{F}}{z_F}[\phi\, v_\alpha](1) - \dfrac{T_e\, \mathcal{F}}{z_F}\ltwoinner{v_\alpha ,\phi_z} +\dfrac{T_e\,\mathcal{G}}{f} \ltwoinner{v_\alpha  ,\phi} \nonumber\\
&&\hspace{-1mm}+\dfrac{T_e}{fz_F}\ltwoinner{r_\alpha\,d \left(\dfrac{1}{z_F}[v_\alpha]_z - v_\alpha\mathcal{M}_\alpha \right),\phi_z}
+\dfrac{T_e}{fz_F}\ltwoinner{r_\alpha\,\beta \left(\frac{1}{z_F}[\rho_\alpha]_z - {\rho_\alpha}\mathcal{M}_\alpha \right),\phi_z}\label{eq3:va}
\end{eqnarray}
Let the bilinear form 
\begin{equation}\label{bilin-v}
 \mathcal{B}({v_\alpha} ,\phi) = \dfrac{T_e\mathcal{G}}{f} \ltwoinner{v_\alpha  ,\phi}+\dfrac{T_e}{z_F^2f}\ltwoinner{r_\alpha\,d\, [v_\alpha]_z ,\phi_z} +\dfrac{T_e \mathcal{F}}{z_F}\phi(1)\, v_\alpha(1,t) - \dfrac{T_e \mathcal{F}}{z_F}\ltwoinner{v_\alpha ,\phi_z}  -\dfrac{T_e\mathcal{M}_\alpha}{z_Ff}\ltwoinner{r_\alpha\,d \, v_\alpha ,\phi_z} \quad
\end{equation} and $U(\rho_\alpha,\phi) = \dfrac{T_e}{z_Ff}\ltwoinner{r_\alpha\,\beta \left(\frac{1}{z_F}[\rho_\alpha]_z - {\rho_\alpha}\mathcal{M}_\alpha \right),\phi_z} = \dfrac{T_e}{z_F^2f}\ltwoinner{r_\alpha\,\beta\, [\rho_\alpha]_z ,\phi_z}  -\dfrac{T_e\mathcal{M}_\alpha}{z_Ff}\ltwoinner{r_\alpha\,\beta \, \rho_\alpha ,\phi_z}$
 then, \eqref{eq3:va} becomes
\begin{equation}\label{eq:semi-var:va}
\ltwoinner{[v_\alpha]_t,\phi} +  \mathcal{B}(v_\alpha ,\phi) + U(\rho_\alpha,\phi) = 0.
\end{equation}

\noindent \textbf{Semi-Variational formulation}\\
Given $\rho_\alpha(z,t)$, seek $v_\alpha: [0, 1]\times [0, 1] \rightarrow \mathbb{R}$ such that for all $t>0$ and $v_\alpha(.,t)  \in \mathcal{T}$
\begin{equation}\label{FirnV:v}
 \begin{cases}
\ltwoinner{[v_\alpha]_t,\phi} +  \mathcal{B}(v_\alpha ,\phi) + U(\rho_\alpha,\phi) = 0& \\
 v_\alpha(z,0) = \overline{v}(z) = 0 &
 \end{cases}
 \end{equation}
 Note that $\mathcal{B}(. ,\phi) = \mathcal{A}(. ,\phi) $ for $T = T_e$ and $D_\alpha = r_\alpha d$. Thus the time and space discretization of the first two terms is similar to that of \eqref{FirnV} with the exception that the vectors with $\rho^{atm}$ are zero since $v_\alpha(0,t) = 0$.

 By integrating (\ref{FirnV:v}) over the temporal interval $[t,t+\Delta t]$, with $0\leq t\leq 1-\Delta t$, one reaches the following L\textsuperscript{2} Integral Formulation:
 \begin{equation}\label{Firn-tau:v}
\hspace{-8mm}\left\{\begin{array}{ll}
\ltwoinner{v_\alpha(z, t+\Delta t)-v_\alpha(z,t),\phi} =-\int_t^{t+\Delta t}  \mathcal{B}(v_\alpha(z,s), \phi(z)) \,ds -\int_t^{t+\Delta t}  U(\rho_\alpha(z,s), \phi(z)) \,ds&\\
 v_\alpha(z,0) = \overline{v}(z) = 0 &\\
\end{array}\right.\vspace{-2mm}
\end{equation}
For the full discretization of equation \eqref{Firn-tau:v}, 
the term $\int_t^{t+\Delta t} \mathcal{B}(v_\alpha(z,s), \phi(z)) \,ds$ is first discretized using an implicit right rectangular rule:
$$\int_t^{t+\Delta t}  \mathcal{B}(v_\alpha(z,s), \phi(z)) ds={\Delta t} \; \mathcal{B}(v_\alpha(z,t+\Delta t), \phi(z))$$ and the second term $\int_t^{t+\Delta t}  U(\rho_\alpha(z,s), \phi(z)) \,ds$ using the trapezoidal rule $$\hspace{-10mm}\int\limits_t^{t+\Delta t}  U(\rho_\alpha(z,s), \phi(z)) \,ds = \dfrac{\Delta t}{2} \; \left[U(\rho_\alpha(z,t+\Delta t), \phi(z))+U(\rho_\alpha(z,t), \phi(z))\right] = \dfrac{\Delta t}{2} \; U(\rho_\alpha(z,t+\Delta t)+\rho_\alpha(z,t), \phi(z))$$ leading to the following fully implicit scheme in time.
 \begin{equation}\label{Firn-tau23:v}
\hspace{-5mm}\left\{\begin{array}{ll}
\ltwoinner{v_\alpha(z, t+\Delta t)-v_\alpha(z,t),\phi}  =-{\Delta t}\;  \mathcal{B}(v_\alpha(z,t+\Delta t), \phi(z))-\dfrac{\Delta t}{2} \; U\left(\rho_\alpha(z,t+\Delta t)+\rho_\alpha(z,t), \phi(z)\right)&\\
 v_\alpha(z,0) = \overline{v}(z) = 0&\\
\end{array}\right.
\end{equation}
Applying Finite Element in space, and using definition \eqref{FEdef} and the fact that $v_\alpha(z_1,t)=0$ and that $\rho_\alpha(z_1,t+\Delta t)+\rho_\alpha(z_1,t)=\rho^{atm}(t+\Delta t)+\rho^{atm}(t)$, then \eqref{Firn-tau23:v} simplifies as follows
 \begin{eqnarray}
 \sum\limits_{i=2}^n\ltwoinner{{(v_\alpha (z_{i},t+\Delta t) - v_\alpha(z_{i},t))\varphi_i},\phi}
&=&-{\Delta t} \mathcal{B}\left(\sum\limits_{i=2}^n{v_\alpha(z_{i},t+\Delta t)\varphi_i} , \phi\right) 
 \nonumber\\
 &&-\dfrac{\Delta t}{2} U\left({\left(\rho^{atm}(t+\Delta t)+\rho^{atm}(t)\right)\varphi_1} , \phi\right)
 \nonumber\\
&&-\dfrac{\Delta t}{2} \,U\left(\sum\limits_{i=2}^n{\left(\rho(z_{i},t+\Delta t)+\rho(z_{i},t)\right)\varphi_i} , \phi\right)\nonumber\vspace{-15mm}
\end{eqnarray}\vspace{-6mm}
 \begin{eqnarray}
\implies{\Delta t}  \mathcal{B}\left(\sum\limits_{i=2}^n{v_\alpha(z_{i},t+\Delta t)\varphi_i} , \phi\right)&+&\sum\limits_{i=2}^n v_\alpha(z_{i},t+\Delta t)\ltwoinner{ \varphi_i,\phi} =\sum\limits_{i=2}^n v_\alpha(z_{i},t)\ltwoinner{\varphi_i,\phi}\nonumber\\
&&-\dfrac{\Delta t}{2} \,U\left({\left(\rho^{atm}(t+\Delta t)+\rho^{atm}(t)\right)\varphi_1} , \phi\right)\nonumber\\
&& -\dfrac{\Delta t}{2} \,U\left(\sum\limits_{i=2}^n{\left(\rho(z_{i},t+\Delta t)+\rho(z_{i},t)\right)\varphi_i} , \phi\right)\;\;\;\quad\label{FirnN:va}
\end{eqnarray}
Let $\phi = \varphi_j$ for $j=2,..,n$ in \eqref{FirnN:va} and define the vectors $\Lambda_\alpha(t) = [\rho_\alpha(z_{2},t), \, \rho_\alpha(z_{3},t), \, \cdots ,\rho_\alpha(z_{n},t)]^T$ and $\mathcal{V}_{\alpha,\beta}(t) = [v_\alpha(z_{2},t), \, v_\alpha(z_{3},t), \, \cdots ,v_\alpha(z_{n},t)]^T$ of length $n-1$, then \eqref{FirnN:va} can be written in Matrix form 
\begin{equation}\label{eq:mat1-res:va}
\hspace{-5mm}\left\{\begin{array}{lcl}
\left[M+T_e\,\Delta t \,C\right]\,\mathcal{V}_{\alpha,\beta}(t+\Delta t) &=&  M\,\mathcal{V}_{\alpha,\beta}(t) - T_e\,\Delta t \,J_\beta\, (\Lambda_\alpha(t+\Delta t)+\Lambda_\alpha(t)) \\
&&-  T_e\,\Delta t\,(p_\beta(t)+p_\beta(t+\Delta t)) \\
\mathcal{V}_{\alpha,\beta}(0)= \overline{\mathcal{V}} = 0
\end{array}\right.
\end{equation} 
where $C = \dfrac{\mathcal{G}}{f} M+\dfrac{1}{z_F^2 f}S(r_\alpha \,d)-\dfrac{\mathcal{M}_\alpha}{z_Ff}A(r_\alpha \,d)+\dfrac{1}{z_F}Q$, and $J_\beta = \dfrac{1}{2z_F^2 f}S(r_\alpha \,\beta)-\dfrac{\mathcal{M}_\alpha}{2z_Ff}A(r_\alpha \,\beta)$, and 
\begin{equation}\label{eq:c2}
p_\beta(t) = \rho^{\rm atm}(t)\left[ \dfrac{1}{2z_F^2f} \ltwoinner{r_\alpha\,\beta \varphi_1' ,\varphi_2'}- \dfrac{\mathcal{M}_\alpha}{2z_Ff}\ltwoinner{r_\alpha\,\beta \varphi_1 ,\varphi_2'} \right]e_1 = \rho^{\rm atm}(t) \,c_{2,\beta} \,e_1. \qquad \quad \end{equation}%.\vspace{1mm}\\
Note that the left-hand side and the first term on the right-hand side of \eqref{eq:mat1-res:va} are obtained exactly the same way as their corresponding terms in \eqref{eq:mat2}. As for the remaining 2 terms, they are obtained as follows.\\ For $j=2,3,\cdots, n$
\begin{itemize}
    \item $U\left(\sum\limits_{i=2}^n\rho(z_{i},t)\varphi_i,\varphi_j\right)  = \dfrac{T_e}{z_F^2f}\sum\limits_{i=2}^n\rho(z_{i},t)\ltwoinner{r_\alpha\,\beta\, \varphi_i' ,\varphi_j'}  -\dfrac{T_e\mathcal{M}_\alpha}{z_Ff}\rho(z_{i},t)\ltwoinner{r_\alpha\,\beta \, \varphi_i ,\varphi_j'}$\\is equivalent to $2\,T_e\,J_\beta\, \Lambda(t)$.
    \item $U\left(\rho^{\rm atm}(t)\varphi_1,\varphi_j\right)  = \dfrac{T_e}{z_F^2f}\rho^{\rm atm}(t)\ltwoinner{r_\alpha\,\beta\, \varphi_1' ,\varphi_j'}  -\dfrac{T_e\mathcal{M}_\alpha}{z_Ff}\rho^{\rm atm}(t)\ltwoinner{r_\alpha\,\beta \, \varphi_1 ,\varphi_j'}$
    \\is equivalent to the $(n-1) \times 1$  vector $2\,T_e\,p_\beta(t)$. 
\end{itemize}
Similarly to $c_1$ in \eqref{eq:c1}, the constant $c_{2,\beta}$ in \eqref{eq:c2} is approximated by \begin{equation} c_{2,\beta}\approx -
\dfrac{r_\alpha}{4\,z_f^2\,f\,z_2}(\beta(z_1)+\beta(z_2))-\dfrac{r_\alpha\,\mathcal{M}_\alpha}{8\,z_f\,f}
(\beta(z_1)+\beta(z_2))\end{equation}

\subsubsection{Computing the gradient using the Directional Derivative $v_\alpha$}\label{sec:gradv}
 Given the directional derivative of $V(\rho_\alpha, d)$ with respect to $d$ is
$$
\mathcal{D}_{\beta} V(d)  \;\;=\;\; \mathcal{D} V \, \cdot \beta \;\;=\;\; \nabla V \,.\, \beta 
\;\;=\;\; 2\sum\limits_{\alpha} \ltwoinner{ \rho_\alpha -  g_{\alpha}, \mathcal{D}_{\beta}\rho_\alpha } \;\;=\;\;  2\sum\limits_{\alpha} \ltwoinner{ \rho_\alpha -  g_{\alpha}, v_{\alpha,\beta}} 
$$
with the directional derivative of $\rho_\alpha$ along $\beta(z)$
$$v_{\alpha,\beta} \;\;=\;\; \mathcal{D}\rho_\alpha\,\cdot\, \beta \;\; =\;\; \mathcal{D}_{\beta}\rho_\alpha \;\;=\;\; \lim\limits_{\epsilon \rightarrow 0} \dfrac{\rho_\alpha(d+\epsilon\beta)-\rho_\alpha(d)}{\epsilon}$$ where $\mathcal{D}\rho_\alpha$ and $\mathcal{D}V$  denote the Fr\'echet derivatives of $\rho_\alpha$ and $V$ with respect to $ d$, then it  is possible compute the gradient entry by entry, considering $\beta$ as the canonical basis $\{e_1, e_2, \cdots, e_n\}$. \\

\noindent At a first glance this might appear computationally intense as it requires solving problem \eqref{eq:mat1-res:va} $n\alpha$ times for computing the gradient once. However, since the matrix $M+T_e\,\Delta t \,C$ is independent from $\beta$ and only the right-hand side vector $T_e\,\Delta t \,J_\beta\, (\Lambda(t+\Delta t)+\Lambda(t)) +  T_e\,\Delta t\,(p_\beta(t)+p_\beta(t+\Delta t))$ is $\beta$-dependent, then it is possible to solve the $n$ $v_{\alpha,e_j}$ systems simultaneously by solving a linear system with multiple right-hand sides.

\noindent Moreover, it is possible to reduce the time needed to compute the right-hand side vectors by noting that the matrices $S(e_j)$ and $A(e_j)$ have at most 7 nonzero entries. Thus,  $J_{e_j}\, (\Lambda(t+\Delta t)+\Lambda(t))$ can be computed without performing matrix-vector multiplication, nor generating the matrices $S(e_j)$ and $A(e_j)$ for $j=1:n$.\\ For that purpose, note that the multiplication of $A(e_j)$ with any vector $v$ costs at most 6 flops

\noindent $A(e_1)v = \dfrac{v_1}{4}e_1, \;A(e_2)v = \dfrac{1}{4}\begin{bmatrix}
    v_2\\
    v_2-v_1\\
    0\\
    0\\
    \vdots\\0
\end{bmatrix}, \;A(e_{n})v = \dfrac{1}{4}\begin{bmatrix}
    0\\
     \vdots\\0\\
   0\\
    v_{n-1}-v_{n-2}\\
    v_{n-1}-v_{n-2}
\end{bmatrix} $\\ $A(e_3)v = \dfrac{1}{4}\begin{bmatrix}
    v_2-v_1\\
    v_3-v_1\\
    v_3-v_2\\
    0\\
    \vdots\\0
\end{bmatrix},  \;A(e_4)v = \dfrac{1}{4}\begin{bmatrix}
    0\\
    v_3-v_2\\
    v_4-v_2\\
    v_4-v_3\\
    \vdots\\0
\end{bmatrix}, \cdots , \;A(e_{n-1})v = \dfrac{1}{4}\begin{bmatrix}
    0\\
     \vdots\\0\\
    v_{n-2}-v_{n-3}\\
    v_{n-1}-v_{n-3}\\
    v_{n-1}-v_{n-2}
\end{bmatrix}.$\\
However, given the repeating patterns in \eqref{eq:Ae}, the $n$ matrix-vector multiplications cost $4(n-1)$ flops  \begin{equation}\label{eq:Ae}
[A(e_1)v\, A(e_2)v \cdots A(e_j)v\cdots A(e_n)v ] =  \dfrac{1}{4}\begin{bmatrix}
    a_0&b_0&a_1&0&0&\cdots&0\\
    0&a_1&b_1&a_2&0&\cdots&0\\
\vdots&\ddots&\ddots&\ddots&\ddots&\ddots&\vdots\\
    0&\cdots&&a_{n-4}&b_{n-4}&a_{n-3}&0\\
     0&\cdots&0&0&a_{n-3}&b_{n-3}&a_{n-2}\\
      0&\cdots&0&0&0&a_{n-2}&a_{n-2}\\
\end{bmatrix}\end{equation}
where $a_j = v_{j+1}-v_j$ for $j=0:n-2$; and
  $b_j = v_{j+2}-v_j$ for $j=0:n-3$; and $v_0=0$ .\\
Similarly for $S(e_j)$, its multiplication with any vector $v$ costs at most 8 flops

\noindent $S(e_1)v = \dfrac{v_1}{2h}e_1, \;S(e_2)v = \dfrac{1}{2h}\begin{bmatrix}
    2v_1-v_2\\
    v_2-v_1\\
    0\\
    0\\
    \vdots\\0
\end{bmatrix}, \;S(e_{n})v = \dfrac{1}{4}\begin{bmatrix}
    0\\
     \vdots\\0\\
   0\\
    v_{n-2}-v_{n-1}\\
    v_{n-1}-v_{n-2}
\end{bmatrix} $\\ $S(e_3)v = \dfrac{1}{2h}\begin{bmatrix}
    v_1-v_2\\
    2v_2-v_1-v_3\\
    v_3-v_2\\
    0\\
    \vdots\\0
\end{bmatrix},  \;S(e_4)v = \dfrac{1}{2h}\begin{bmatrix}
    0\\
    v_2-v_3\\
    2v_3-v_2-v_4\\
    v_4-v_3\\
    \vdots\\0
\end{bmatrix}, \cdots , \;S(e_{n-1})v = \dfrac{1}{2h}\begin{bmatrix}
    0\\
     \vdots\\0\\
    v_{n-3}-v_{n-2}\\
    2v_{n-2}-v_{n-1}-v_{n-3}\\
    v_{n-1}-v_{n-2}
\end{bmatrix}.$\\
 Given the repeating patterns in \eqref{eq:Se}, the $n$ matrix-vector multiplications cost $2(n-1)+n-2$ flops   
 \begin{equation}\label{eq:Se} 
 \hspace{-3mm}[S(e_1)v\, S(e_2)v \cdots S(e_j)v\cdots S(e_n)v ] =  \dfrac{1}{2h}\begin{bmatrix}
    a_0&c_0&-a_1&0&0&\cdots&0\\
    0&a_1&c_1&-a_2&0&\cdots&0\\
\vdots&\ddots&\ddots&\ddots&\ddots&\ddots&\vdots\\
    0&\cdots&0&a_{n-4}&c_{n-4}&-a_{n-3}&0\\
     0&\cdots&0&0&a_{n-3}&c_{n-3}&-a_{n-2}\\
      0&\cdots&0&0&0&a_{n-2}&a_{n-2}\\
\end{bmatrix}\;\end{equation}
where $a_j = v_{j+1}-v_j$
for $j=0:n-2$; and $c_j = 2v_{j+1} - v_{j+2} -v_{j} = a_j - a_{j+1}$  for $j=0:n-3$; and $v_0=0$.\\
 Thus, all the right-hand side vectors $J_{e_j} v$ can be computed simultaneously using $O(n)$ flops in each time iteration of \eqref{eq:mat1-res:va}, to obtain the $(n-1)\times n$ block matrix $\mathcal{J}$, where the vector $v = \Lambda(t+\Delta t)+\Lambda(t)$  is given\begin{equation}\label{eq:Je} 
 \hspace{-3mm}\mathcal{J} = [J_{e_1}v\;\;\; J_{e_2}v \;\cdots \;J_{e_j}v\;\cdots \;J_{e_n}v ] = 
 \dfrac{r_\alpha}{2z_F^2 f} \eqref{eq:Se} -\dfrac{r_\alpha\,\mathcal{M}_\alpha}{2z_Ff}\eqref{eq:Ae}\,.
 \end{equation}
 Moreover, given that \begin{equation}\label{eq:c2e1}
     c_{2,e_1} = c_{2,e_2} = -\dfrac{r_\alpha}{4\,z_f^2\,f\,z_2}-\dfrac{r_\alpha\,\mathcal{M}_\alpha}{8\,z_f\,f} \mbox{ and }c_{2,e_j} = 0 \mbox{ for } j=3:n,  \end{equation}
     then $$p_{e_1}(t) = p_{e_2}(t) = \rho^{atm}(t)\;c_{2,e_2}\;e_1, \qquad p_{e_j}(t) = 0 \;\; for \;\; j=3:n$$
 Similarly, $$p_{e_1}(t+\Delta t) = p_{e_2}(t+\Delta t) = \rho^{atm}(t+\Delta t)\;c_{2,e_2}\;e_1, \qquad p_{e_j}(t+\Delta t) = 0 \;\; for \;\; j=3:n$$
 Thus, computing $\hat{p}_{\beta}(t) = {p}_{\beta}(t) + {p}_{\beta}(t+\Delta t) = c_{2,\beta} \;(\rho^{atm}(t)+\rho^{atm}(t+\Delta t))\; e_1$ for the canonical basis consists of 2 nonzero entries of the $(n-1)\times n$ matrix below
 \begin{equation}\label{eq:pe}[\hat{p}_{e_1}(t)\;\; \hat{p}_{e_2}(t)\;\; \hat{p}_{e_3}(t)\;\; \cdots \;\; \hat{p}_{e_n}(t)] =     c_{2,e_2} \;(\rho^{atm}(t)+\rho^{atm}(t+\Delta t))  \begin{bmatrix}
 1&1&0&\cdots&0\\
 0&0&0&\cdots&0\\
 \vdots&\vdots&\vdots&&\vdots\\
  0&0&0&\cdots&0\\
 \end{bmatrix}
 \end{equation}
 Thus, it is possible to find $v_{\alpha,e_j}(t)$ for $j=1:n$ and $t\in (0,1]$ simultaneously by solving
 \begin{equation}\label{eq:mat1-res:va:block}
\left\{\begin{array}{lcl}
\left[M+T_e\,\Delta t \,C_\alpha\right]\,\mathcal{V}_{\alpha}(t+\Delta t) &=&  M\,\mathcal{V}_\alpha(t) - T_e\,\Delta t \; \mathcal{R}\\
\mathcal{V}(0)=  0_{(n-1)\times n}
\end{array}\right.
\end{equation} where the $(n-1)\times n$ block $\mathcal{R}$ consists of adding \eqref{eq:pe} and \eqref{eq:Je}; i.e. adding $ c_{2,e_1}\;(\rho^{atm}(t)+\rho^{atm}(t+\Delta t))$ to the first 2 entries in the first row of matrix \eqref{eq:Je}; and the $(n-1)\times n$ block $\mathcal{V}_\alpha(t)$ is $$\mathcal{V}_\alpha(t) = \begin{bmatrix}
 v_{\alpha,e_1}(z_{2},t) &v_{\alpha,e_2}(z_{2},t) &\cdots &v_{\alpha,e_n}(z_{2},t) \\
 v_{\alpha,e_1}(z_{3},t)&v_{\alpha,e_2}(z_{3},t)&\cdots&v_{\alpha,e_n}(z_{3},t)\\
 \vdots &\vdots& &\vdots \\
 v_{\alpha,e_1}(z_{n},t)&v_{\alpha,e_2}(z_{n},t)&\cdots&v_{\alpha,e_n}(z_{n},t)\\
\end{bmatrix}$$
  At the end of the time iterations of \eqref{eq:mat1-res:va:block}, the $(n-1)\times n$ block $\mathcal{V}_\alpha(1)$ is obtained. Note that a zero first row could be added to $\mathcal{V}_\alpha(1)$  to get an $n\times n$ matrix due to the boundary condition $v_{\alpha}(0,t) = 0$. This procedure is summarized in Algorithm \ref{alg:Bvalpha}.

The above procedure is repeated for each gas $\alpha$. Then, 
$$ \nabla V ^T
\;\;=\;\; 2\sum\limits_{\alpha} ( \Lambda_\alpha(1) -  g_{\alpha})^T \mathcal{V}_\alpha(1)$$
This gradient calculation is summarized in Algorithm \ref{alg:grad}.
Note that since the matrix $B_\alpha = M+T_e\,\Delta t \, C_\alpha$ is fixed throughout the time iterations of both Algorithms \ref{alg:Firn} and \ref{alg:Bvalpha} for the same $D_\alpha$, we perform a PLU decomposition of the matrix once, and pass the obtained triangular matrices to the Algorithms, instead of the matrix $C_\alpha$. 
\begin{algorithm}[H]
\centering 
\caption{  The Block $v_\alpha$ Direct Problem  }
{\renewcommand{\arraystretch}{1.3}
\begin{algorithmic}[1]
\Statex{\textbf{Input:} \;\;\; Mass matrix $M$; Matrix $C_\alpha$ as defined in \ref{eq:C}; End time $T_e$; Time step $dt$; Mesh size $h$ }
\Statex{ \qquad \qquad  $n\times m$  matrix $\Lambda$, output of Algorithm \ref{alg:Firn}; $c_{2,e_1}$ as defined in \eqref{eq:c2e1}; The function $\rho^{atm}(t)$.}
 \Statex{\textbf{Output:} ${\mathcal{V}_\alpha}$: $(n-1)\times n$ matrix with the computed solution vectors for $\beta = e_1, e_2, \cdots, e_n$ and }
 \Statex{\qquad\quad \;\;\; $t=1$, i.e. $v_{\alpha,e_j}(z_i,1)$ for $i=2:n$ and $j=1:n$. \vspace{3mm}}
 \State $z = 0:h:1$\,; \;\;\;$n = length(z)$;\;\;\;$t= 0:dt:1 ; \;\;\; m = length(t); \;\;\; \rho =  \rho^{atm}(t);$ \vspace{1mm}
 \State $\mathcal{V}_1 = zeros(n-1,n)$;\;\;\; $B_\alpha = (M + T_e *dt* C_\alpha)  $\vspace{1mm}
\For {$i= 1:m-1$\vspace{2mm}} 
\State Compute the block matrix $\mathcal{J}$ as defined in \eqref{eq:Je}, and let $\mathcal{R} = \mathcal{J}$ ;\vspace{1mm}
\State Let $p= c_{2,e_1}*(\rho(i) + \rho(i+1))$\,; and $\mathcal{R}(1,1) = \mathcal{R}(1,1)+ p$;\;\;\; $\mathcal{R}(1,2) =\mathcal{R}(1,2)+  p;$\vspace{1mm}
\State  $RHS=M*\mathcal{V}_{i}-  T_e *dt*\mathcal{R};$ \vspace{1mm}
\State  $ \mathcal{V}_{i+1} = B_\alpha \backslash RHS;$ \qquad $\%$ Solve for $\mathcal{V}_{i+1} $\vspace{1mm}
\EndFor
\State  $\mathcal{V}_\alpha =  \mathcal{V}_m; $ \vspace{2mm}
\end{algorithmic}}
\label{alg:Bvalpha}
\end{algorithm}

\begin{algorithm}[H]
\centering 
\caption{  Gradient of objective function $V(d)$}
{\renewcommand{\arraystretch}{1.3}
\begin{algorithmic}[1]
\Statex{\textbf{Input:} \;$d$: CO2 Diffusion coefficient; cnsts: Direct Problem constants; }
\Statex{\qquad \quad\;\;$g_{\alpha,i}(1) = \rho_{\alpha,g}(z_i,1)$: given gas concentrations at time $t=1$ and positions $z_i, i=1:n$;   }
 \Statex{\textbf{Output:} $V(d)$: the objective function evaluated at $d$,} 
 \Statex{\qquad \quad \;\;${\nabla V (d)}$: $1\times n$ vector with the Gradient of the objective function $V(d)$, }\vspace{0.1mm}
 \State{Generate the matrices $M$ and $K$.}
  \For {each gas $\alpha$}%\vspace{2mm}} 
  \State Compute the constant $c_1$ as defined in \eqref{eq:c1}
  \State{Let $D_\alpha = r_\alpha d$ and generate the matrices $A(D_\alpha)$, $S(D_\alpha)$, $C_\alpha$, and $B_\alpha$. }
  \State{Compute the PLU Decomposition of $B_\alpha$, $P_\alpha B_\alpha = L_\alpha U_\alpha$} 
 \State {Call the Direct Problem Algorithm \eqref{alg:Firn} for $M, P_\alpha, L_\alpha, U_\alpha$, to get the $n\times m$ solution matrix $\Lambda_\alpha$}
 \vspace{1mm}
 \State Call the Block $v_\alpha$ Algorithm \ref{alg:Bvalpha} for  $M, P_\alpha, L_\alpha, U_\alpha$, to get the $(n-1)\times n$ matrix $\mathcal{V}_\alpha$ 
 \vspace{1mm}
 \EndFor
\State $ \nabla V (d) = \sum_\alpha (\Lambda_\alpha(2:n,m) - g_{\alpha}(2:n))^T \mathcal{V}_\alpha$; \vspace{1mm}
\State  $V (d) = \sum_\alpha (\Lambda_\alpha(:,m) - g_{\alpha})^T(\Lambda_\alpha(:,m) - g_{\alpha})$; \vspace{1mm}
\end{algorithmic}}
\label{alg:grad}
\end{algorithm}

\subsection{Testing}\label{sec:test}
We first start by briefly discussing different constrained and unconstrained minimization methods to solve our inverse problem in section \ref{sec:min}.
 Then, we test these methods using our gradient algorithm \ref{alg:grad} or with finite difference approximations, and compare the obtained results in sections .. 

In all the testings, we generate data using the direct
problem for 3 gazes with $r_\alpha = [0.5, 1, 1.5]$ and $D_{CO2}$ corresponding to test cases described in the Direct
problem part. Moreover, we consider $T_e = 1, 50, 100, 150, z_F = 1, 5, 10$,  $dt = h_g$, and $h_g = 1/65$ or $1/128$; where a data set is generated for each combination. Then, for each data set, $\rho_\alpha(z,1)$ values  for $h = 1/16, 1/32, 1/64$ are generated using linear splines, as the inverse problem is solved using these meshings ($h = 1/16, 1/32, 1/64$). 

 \subsubsection{Minimization Algorithms} \label{sec:min}  
  We start by briefly introducing the steepest descent, and  nonlinear Conjugate Gradient. 
  As mentioned earlier, 
  both methods aim at minimizing $V(d)$, by starting with an initial guess $d_0$, then iterating with%
  $$d_{n+1} = d_n +\alpha_n\delta_n,$$ where $\delta_n$ is the descent's direction and $\alpha_n$ is the step length along this direction. 
  
\noindent  The computation of $\delta_n$ depends on that of 
 the gradient of V. The steepest descent method considers  $$\delta_n = -\nabla V(d_n)$$ which is the steepest descent from $d_n$. Whereas, the nonlinear Conjugate Gradient method considers $$\begin{cases}
\delta_0 = -\nabla V(d_0) &\\
\delta_n = -\nabla V(d_n) + \beta_{n-1} \delta_{n-1}
 \end{cases}$$
There are different choices for $\beta_{n-1}$ based on the already computed $g_n =  \nabla V(d_n), \,g_{n-1} ,\, \delta_{n-1}$, and $y_{n-1} = (g_n - g_{n-1})$, such as:
\begin{itemize}
    \item $\beta_{n-1}^{HS} = \dfrac{g_n^Ty_{n-1}}{\delta_{n-1}^Ty_{n-1}}$, \quad Hestenes and Stiefel, 1952 \vspace{-2mm}
    \item $\beta_{n-1}^{FR} = \dfrac{\ltwonorm{g_n}^2}{\ltwonorm{g_{n-1}}^2}$, \quad Fletcher and Reeves,  1964 \vspace{-2mm}
    \item  $\beta_{n-1}^{PR} = \dfrac{g_n^Ty_{n-1}}{\ltwonorm{g_{n-1}}^2}$, \quad Polak and Ribi\`ere, 1969 \vspace{-2mm}
     \item $\beta_{n-1}^N = \beta_{n-1}^{HS} - \dfrac{2\ltwonorm{y_{n-1}}^2}{(\delta_{n-1}^Ty_{n-1})^2} g_n^T\delta_{n-1}$,  \quad Hager and Zhang, 2005
\end{itemize}
 \noindent Once the descent's direction is set, it remains to compute the step length $\alpha_n$. In both methods, $\alpha_n$ is obtained by solving the following line search $$\min\limits_{\alpha \ge 0} \phi_n (\alpha)= \min\limits_{\alpha \ge 0} V(d_n +\alpha\delta_n)$$ 
where $\alpha$ should satisfy the Armijo condition and curvature condition, or what is known as the Wolfe condition. For a survey of different nonlinear CG methods, refer to \cite{NCG}.

In our testings we use 
nonlinear Conjugate Gradient method as implemented in the Poblano's MATLAB Toolbox \cite{poblano} for solving the unconstrained version of the inverse problem, where the gradient is computed via Algorithm \ref{alg:grad}. A line search satisfying the strong Wolfe conditions is used to guarantee global convergence of the Poblano optimizers. Note that we modified the algorithm to include the "Hager and Zhang" nonlinear CG version. We have also tested nonlinear Conjugate Gradient method, where the gradient is approximated using Finite Differences as implemented in ``Adaptive Robust Numerical Differentiation" Toolbox \cite{derivest}.  
 
Matlab's Optimization Toolbox has several functions that solve unconstrained and constrained optimization problems, such as fminUn and fminCon with different methods such as SQP (sequential quadratic programming) and Interior Point methods.  As mentioned, Matlab's optimization functions require the gradient as input, and if not available a Finite Difference approximation is used. We test these methods on our problem where the gradient is computed via Algorithm \ref{alg:grad} and compare the results with the finite difference approximation version.

\subsubsection{Unconstrained Inverse Problem}\label{sec:noncg}
We solve the unconstrained inverse problem using 
nonlinear Conjugate Gradient (NCG) method as implemented in the Poblano's MATLAB Toolbox \cite{poblano}, where the stopping criteria are set to $10^{-8}$. We also solve it using MATLAB fminUnc's Quasi-Newton method with tolerance $10^{-6}$. The data sets are generated using a mesh size $h_g = 1/65$. 

Table \ref{tab:inv4P} show 
the runtime of the tested algorithms in seconds needed till convergence in (iter) iterations and the L2 relative error with the exact $D_{CO2}$. The
obtained results for NCG with the "Hager and Zhang" (HZ) and "Hestenes and Stiefel" (HS) beta updates for Cases 1 and 2d respectively are shown, where the gradient is computed via Algorithm \ref{alg:grad}. The algorithm was tested where the gradient is computed using Finite Differences as implemented in ``Adaptive Robust Numerical Differentiation" Toolbox \cite{derivest}, but the results are not shown as it is at least 100 times slower than the Algorithm \ref{alg:grad} version, with similar relative errors.

Moreover, the results for MATLAB fminUnc's Quasi-Newton method are shown, where the gradient is computed using MATLAB's built-in Finite Differences or Algorithm \ref{alg:grad}. Both versions converge in number of iterations of the same order with similar relative errors. It is clear that the version with Algorithm \ref{alg:grad} is faster, specifically 10 times faster for $h=1/32, 1/64$.

\begin{figure}[H]
\begin{tabular}{cc}
\subfloat{\includegraphics[scale=0.4]{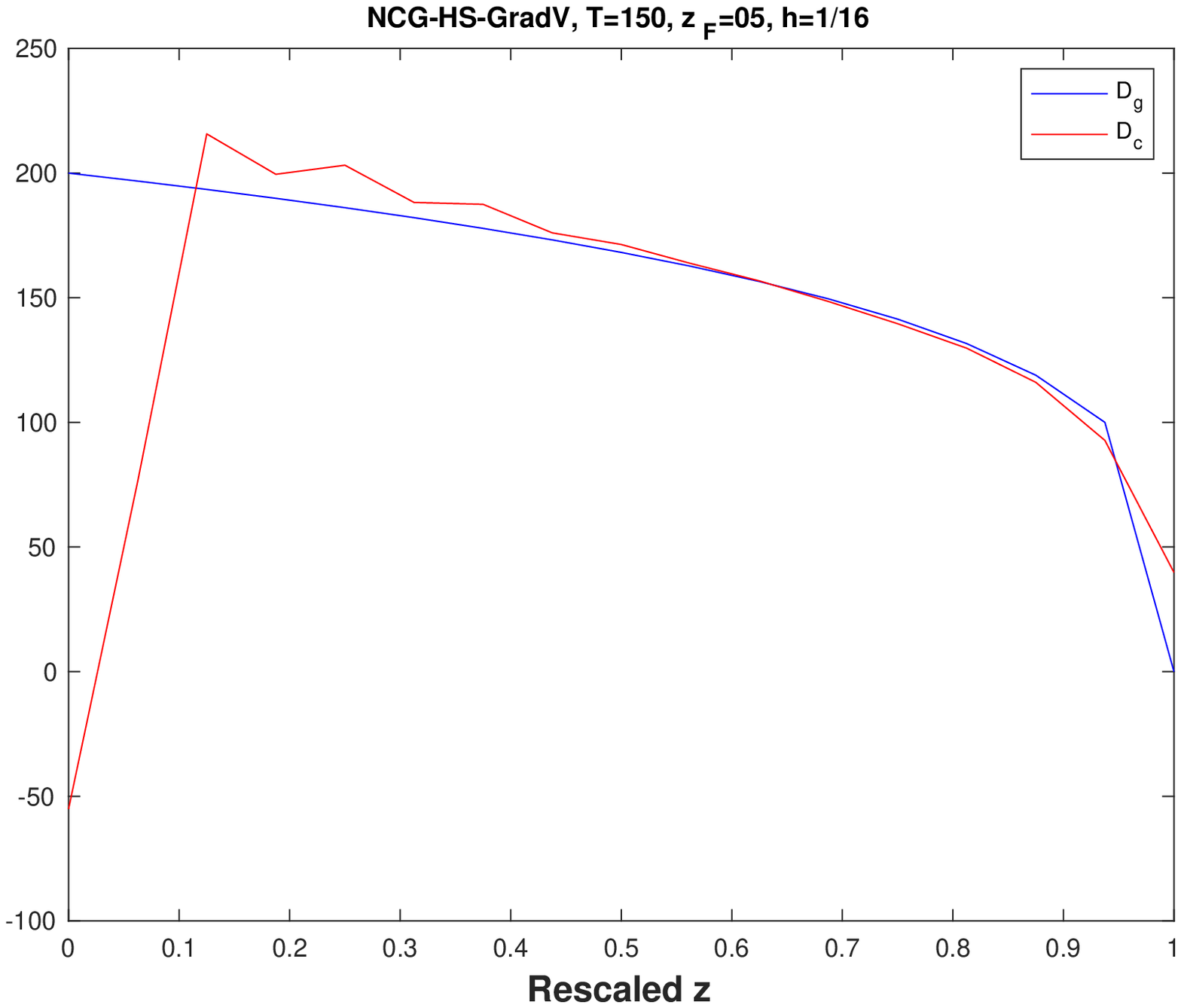}} \hspace{-10mm}&
\subfloat{\includegraphics[scale=0.4]{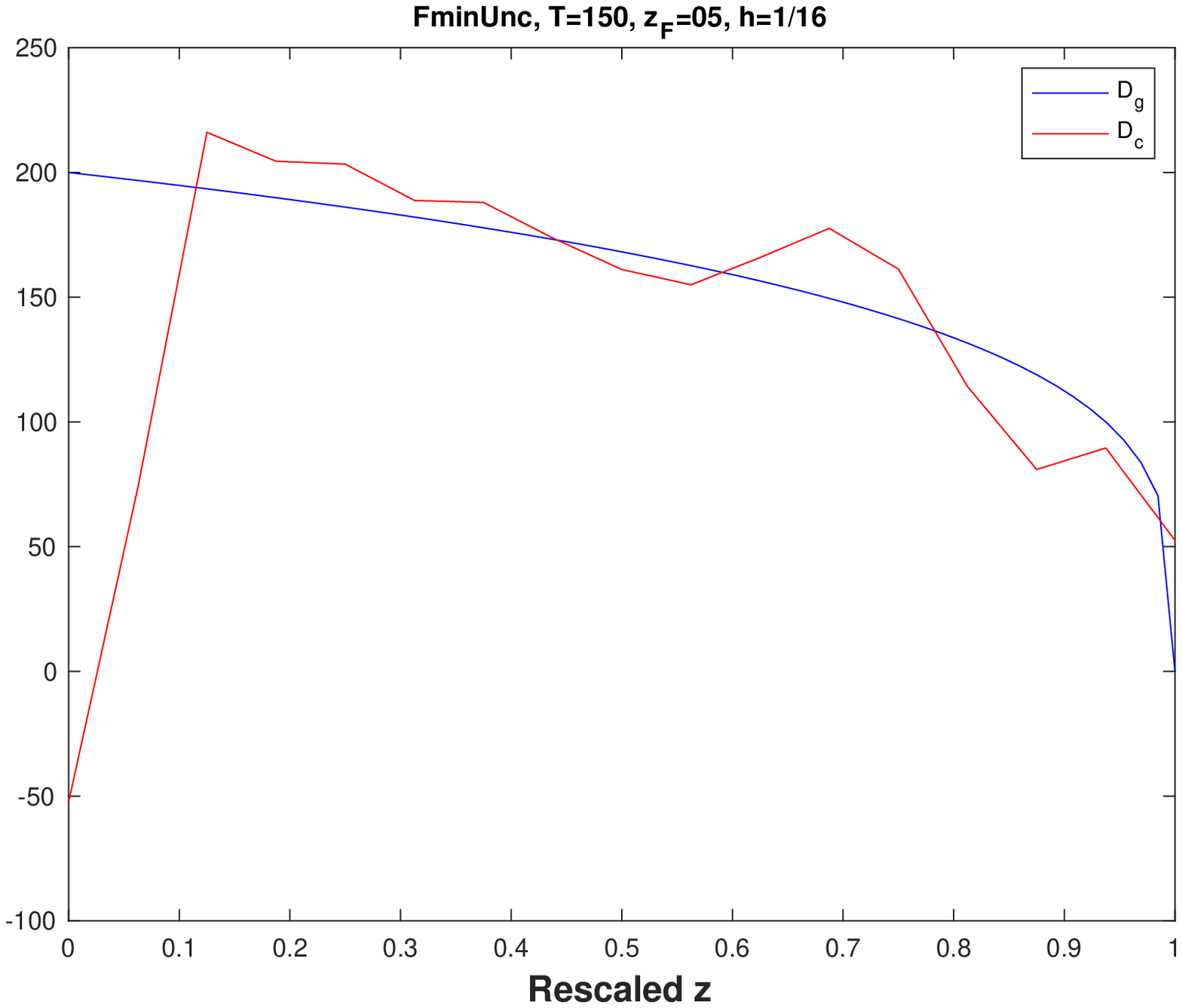}}\vspace{-5mm}\\
\subfloat{\includegraphics[scale=0.4]{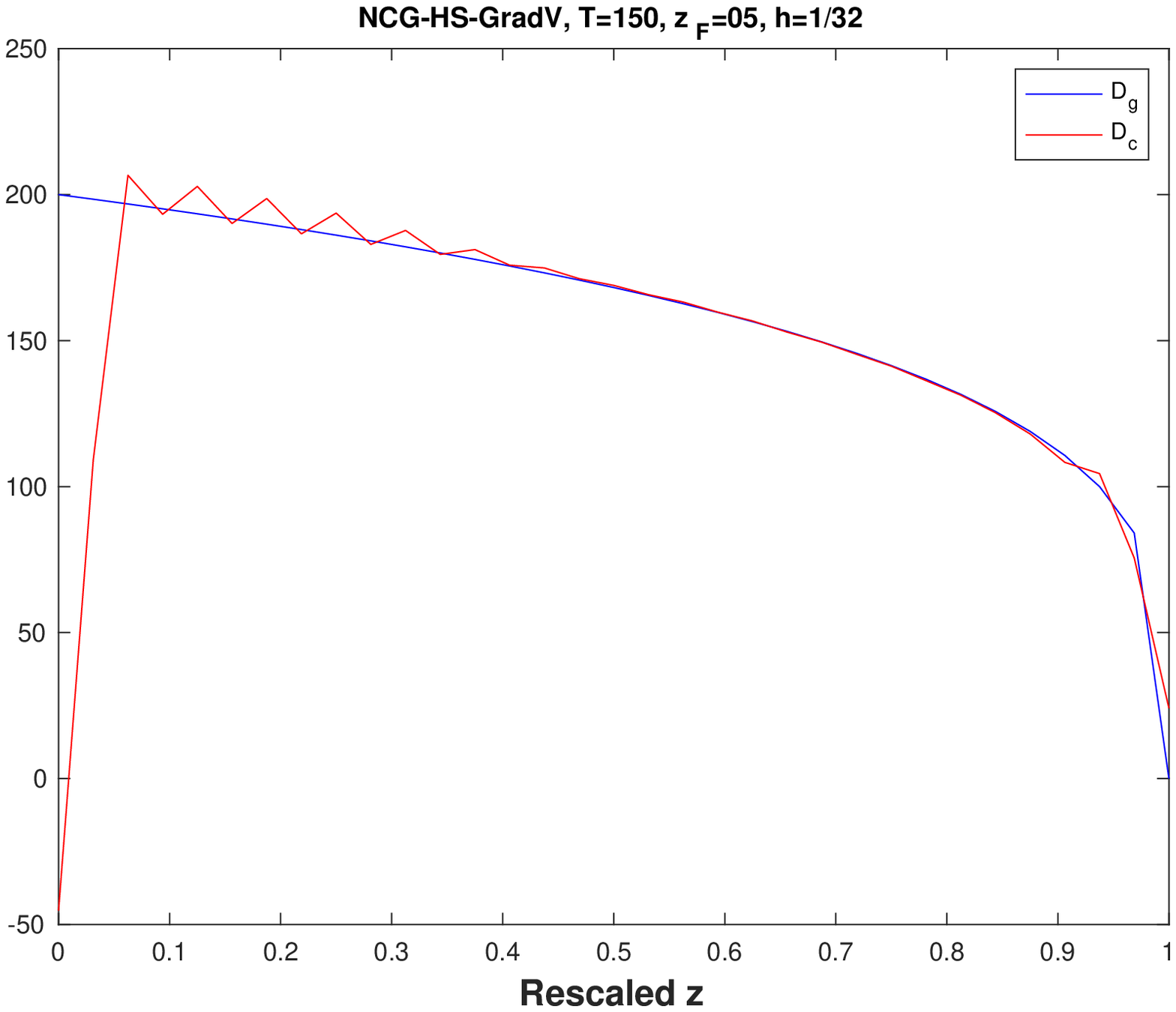}} &
\subfloat{\includegraphics[scale=0.4]{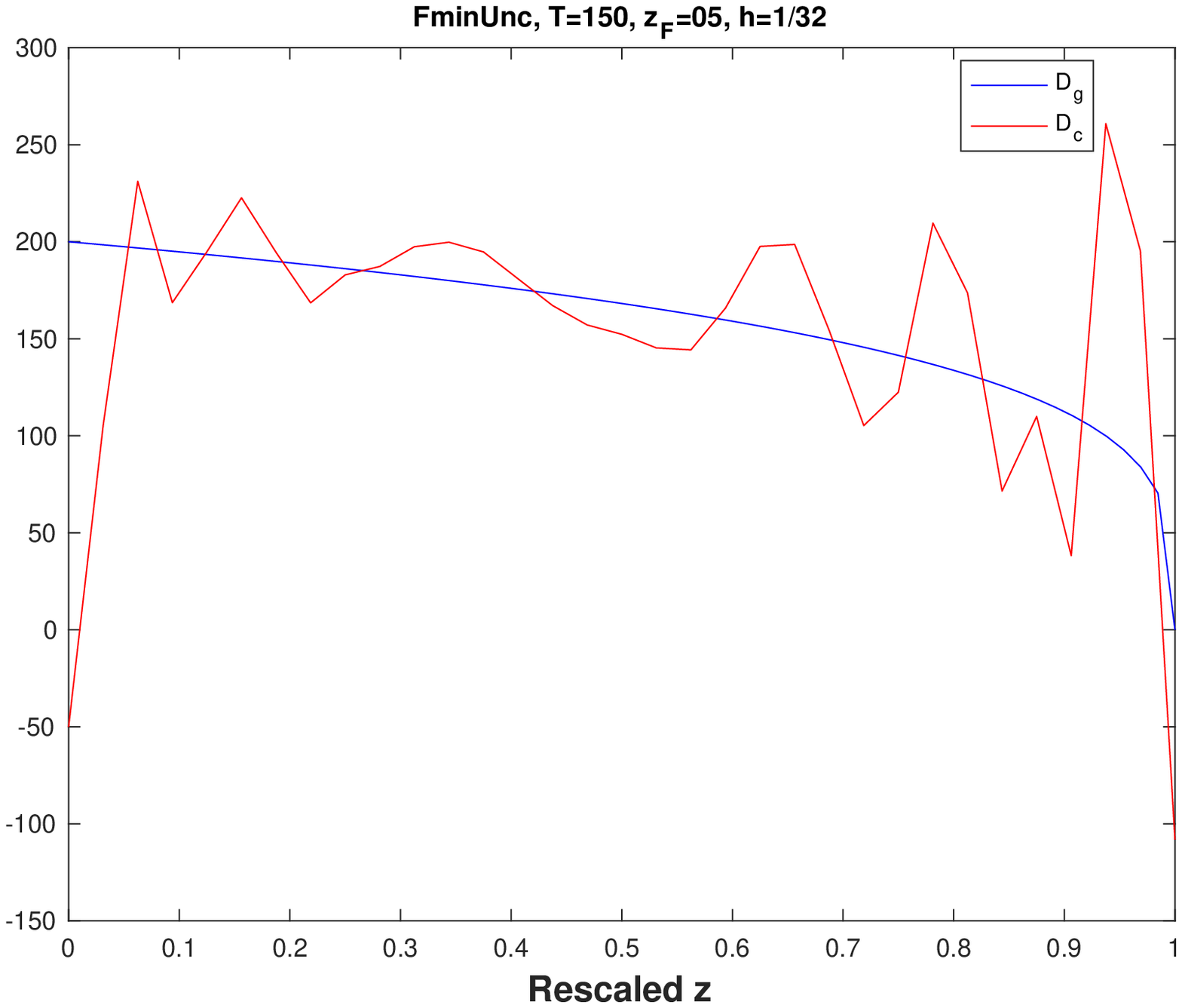}}\vspace{-5mm}\\
\end{tabular}
\centering
\caption{\it \small The $D_{CO2}$ solution for
Case 2d, using NCG (left) and FminUnc (right) methods with initial guess $d_0 = 0$,  $z_F = 5$, and $T_e = 150$, for $h =1/16,1/32$}\label{fig:fminUnc}
\end{figure}

All the methods fail to converge for $h=1/64$, with the exception of NCG, Case1, $z_F=5, T_e = 50$. Moreover, the relative error is of the order of $10^{-1}$ when the methods converge. Figure \ref{fig:fminUnc} shows a sample of obtained $D_{CO2}$ for both NCG and FminUnc. At $z=0$ we get a negative value of $D_{CO2}(0) = -50$. One could set these negative values to zero, since  $D_{CO2} \geq 0$. Another observation is that the obtained  $D_{CO2}$ is not smooth. Thus, it is possible to apply some polynomial regression to the obtained data, to get a smoother solution. 

\begin{table}[H]
%\centering
\small
 \hspace{-2.3cm} \begin{tabular}{|c|c|c||l|l|l|l|l|l||l|l|l|l|l|l|}
       \cline{4-15}
      \multicolumn{1}{c}{}&\multicolumn{1}{c}{}&\multicolumn{1}{c|}{}&\multicolumn{6}{c||}{NCG with Alg \ref{alg:grad} }&\multicolumn{6}{c|}{fminUnc}\\
     \cline{4-15}
 \multicolumn{1}{c}{}& \multicolumn{1}{c}{}& \multicolumn{1}{c|}{} & \multicolumn{3}{c|}{HZ} & \multicolumn{3}{c||}{HS} & \multicolumn{3}{c|}{Quasi-Newton with FD} & \multicolumn{3}{c|}{Quasi-Newton with Alg \ref{alg:grad}  }   \\ \hline
        \multicolumn{1}{|c}{$z_F$}&\multicolumn{1}{|c}{$T_e$}&\multicolumn{1}{|c||}{$h$} & Time/s & iter & L2 err & Time/s & iter & L2 err & Time/s & iter & L2 err& Time/s & iter & L2 err   \\ \hline
   \multirow{9}{*}{5} & \multirow{3}{*}{50} & 
          1/16 & 7.83E-1 & 80 & 3.58E-1 & 1.52E+0 & 74 & 4.30E-1 & 4.14E-1 & 92 & 4.40E-1 & 1.27E-1 & 93 & 4.40E-1  \\ \cline{3-15}
        ~ & ~ & 1/32 & 2.87E+0 & 97 & 2.38E-1 & 6.64E+0 & 126 & 2.82E-1 & 3.47E+0 & 162 & 3.54E-1 & 6.58E-1 & 187 & 4.25E-1   \\\cline{3-15}
        ~ & ~ & 1/64 & 2.32E+2 & 1418 & 4.02E-1 & 5.43E+2 & 1895 & 2.04E+0 & 5.62E+1 & 292 & 5.48E+0 & 6.03E+0 & 325 & 6.75E+0   \\ \cline{2-15}        
        ~ & \multirow{3}{*}{100} &
        1/16 & 6.75E-1 & 70 & 3.59E-1 & 1.47E+0 & 87 & 4.30E-1 & 4.18E-1 & 93 & 4.40E-1 & 1.23E-1 & 93 & 4.40E-1  \\ \cline{3-15}
        ~ & ~ & 1/32 & 3.55E+0 & 122 & 2.37E-1 & 1.40E+2 & 4741 & 9.28E+0 & 3.64E+0 & 162 & 5.69E-1 & 4.25E-1 & 136 & 4.12E-1  \\ \cline{3-15}
        ~ & ~ & 1/64 & 3.26E+2 & 2020 & 4.31E-1 & 5.43E+2 & 1338 & 1.61E+0 & 4.03E+1 & 229 & 2.19E+2 & 3.91E+0 & 216 & 2.89E+2 \\ \cline{2-15}

        ~ & \multirow{3}{*}{150} & 
        1/16 & 5.90E-1 & 61 & 3.59E-1 & 1.25E+0 & 70 & 4.30E-1 & 4.27E-1 & 92 & 4.40E-1 & 1.23E-1 & 94 & 4.40E-1  \\ \cline{3-15}
        ~ & ~ & 1/32 & 3.87E+0 & 132 & 2.37E-1 & 7.78E+0 & 137 & 2.82E-1  & 4.02E+0 & 171 & 4.27E-1 & 4.22E-1 & 136 & 4.12E-1 \\ \cline{3-15}
        ~ & ~ & 1/64 & 3.22E+2 & 2031 & 2.02E-1 & 6.97E+2 & 2342 & 1.54E+0 & 4.68E+1 & 257 & 3.72E+1 & 5.49E+0 & 288 & 6.99E+0  \\\hline
        \multirow{9}{*}{10}  & \multirow{3}{*}{50}  &
        1/16 & 6.19E-1 & 64 & 4.85E-1 & 1.41E+0 & 84 & 6.05E-1 & 4.48E-1 & 89 & 6.13E-1 & 1.25E-1 & 91 & 6.12E-1  \\ \cline{3-15}
        ~ & ~ & 1/32 & 3.30E+0 & 117 & 2.93E-1 & 6.92E+0 & 128 & 3.49E-1 & 2.25E+0 & 102 & 5.51E-1 & 3.33E-1 & 103 & 5.51E-1 \\ \cline{3-15}
        ~ & ~ & 1/64 & 3.59E+1 & 225 & 5.35E-2 & 2.88E+2 & 1636 & 3.26E+0 & 3.88E+1 & 240 & 6.28E-1 & 3.72E+0 & 240 & 8.68E-1  \\ \cline{2-15}
        ~ & \multirow{3}{*}{100}  & 
        1/16 & 6.44E-1 & 69 & 4.89E-1 & 1.33E+0 & 74 & 5.98E-1  & 4.70E-1 & 93 & 6.13E-1 & 1.27E-1 & 91 & 6.12E-1 \\ \cline{3-15}
        ~ & ~ & 1/32 & 4.05E+0 & 148 & 2.86E-1 & 7.89E+0 & 140 & 3.43E-1 & 2.27E+0 & 103 & 5.51E-1 & 3.24E-1 & 103 & 5.51E-1  \\ \cline{3-15}
        ~ & ~ & 1/64 & 2.74E+1 & 175 & 5.41E-2 & 9.23E+2 & 6446 & 1.63E+1 & 3.57E+1 & 229 & 6.75E-1 & 3.90E+0 & 249 & 6.47E-1  \\ \cline{2-15}
        ~ & \multirow{3}{*}{150}  & 
        1/16 & 5.89E-1 & 61 & 4.86E-1 & 1.40E+0 & 94 & 5.96E-1 & 4.69E-1 & 94 & 6.13E-1 & 1.24E-1 & 91 & 6.12E-1   \\\cline{3-15}
        ~ & ~ & 1/32 & 4.27E+0 & 152 & 2.86E-1 & 8.38E+0 & 151 & 3.43E-1 & 2.25E+0 & 102 & 5.51E-1 & 3.23E-1 & 103 & 5.51E-1 \\ \cline{3-15}
        ~ & ~ & 1/64 & 4.67E+2 & 2939 & 7.50E-1 & 1.41E+3 & 6573 & 1.71E+1  & 2.98E+1 & 188 & 7.46E-1 & 3.57E+0 & 231 & 8.96E-1 \\ \hline
    \end{tabular}
    \caption{Comparison of the convergence of Poblanos's NCG using HZ and HS beta updates, and MATLAB fminUnc's Quasi-Newton method with built-in Finite Differences or Algorithm \ref{alg:grad},  with initial guess $d_0 = 0$, for Case 2d with different $z_F, T_e,$ and $h$ values.
      }\label{tab:inv4P}
\end{table}

\subsubsection{Constrained Inverse Problems }
We test Matlab fmincon's SQP and interior point (IP) methods for two constrained inverse problems, the first with the constraint that $d\geq 0$ \eqref{eq:min}, and the second with the constraints that $d\geq 0$ and $d$ is decreasing ($d_i\geq d_{i+1}$). We solve these inverse problems using all the generated data sets with a mesh size $h_g = 1/65$. 

\begin{figure}[H]
\begin{tabular}{ll}
\subfloat{\includegraphics[scale=0.4]{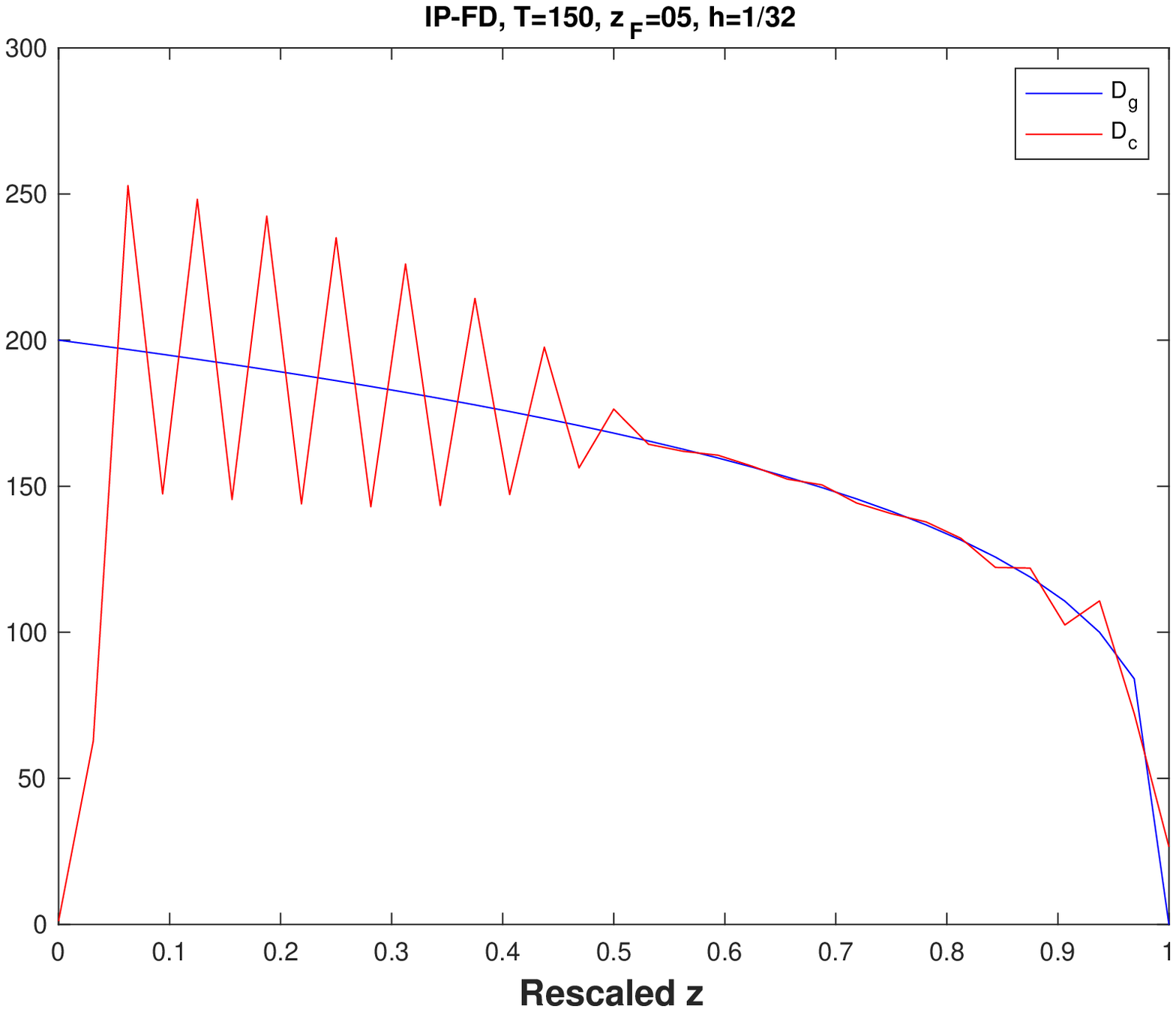}} \hspace{-10mm}&
\subfloat{\includegraphics[scale=0.4]{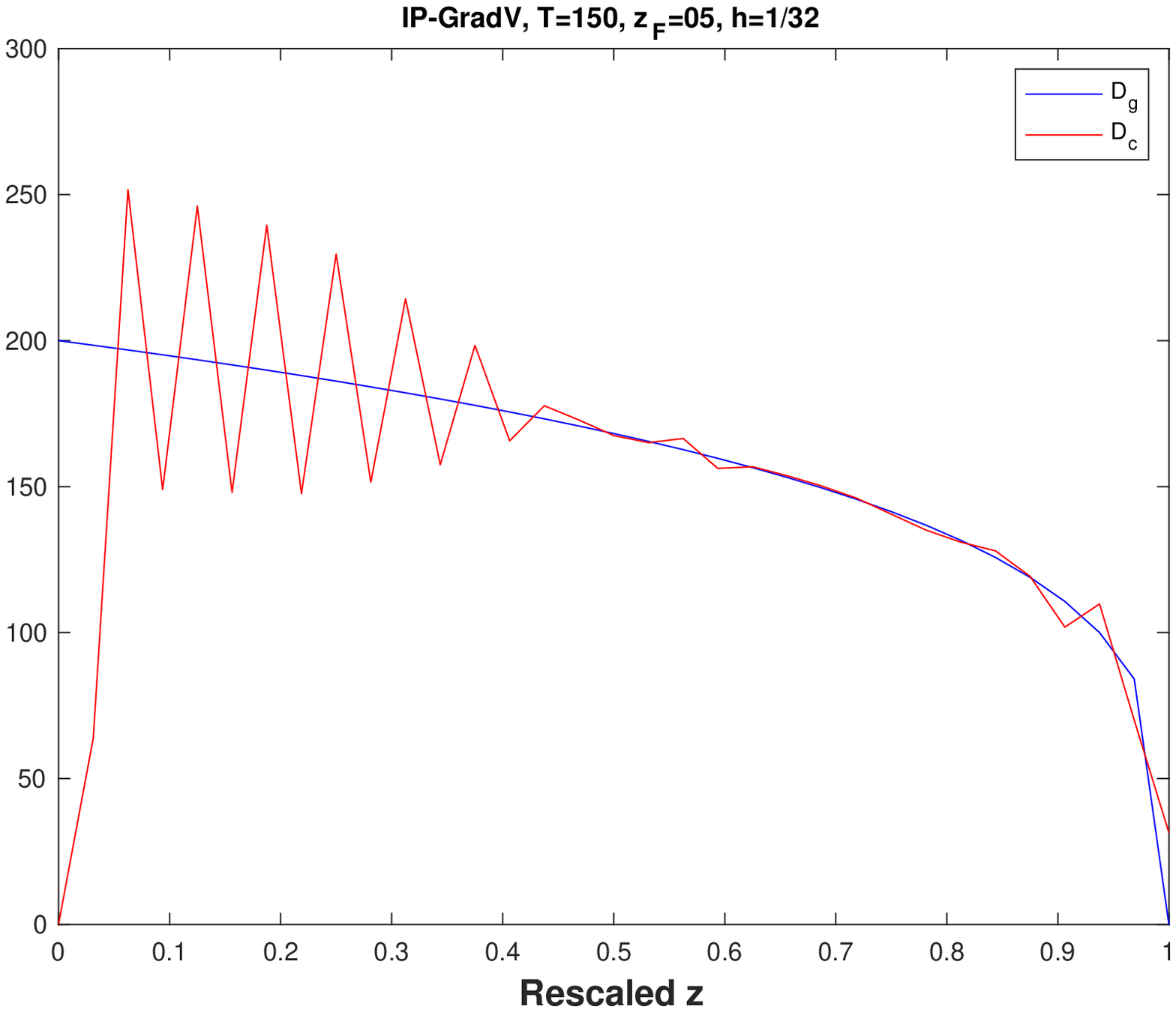}}\vspace{-5mm}\\
\subfloat{\includegraphics[scale=0.4]{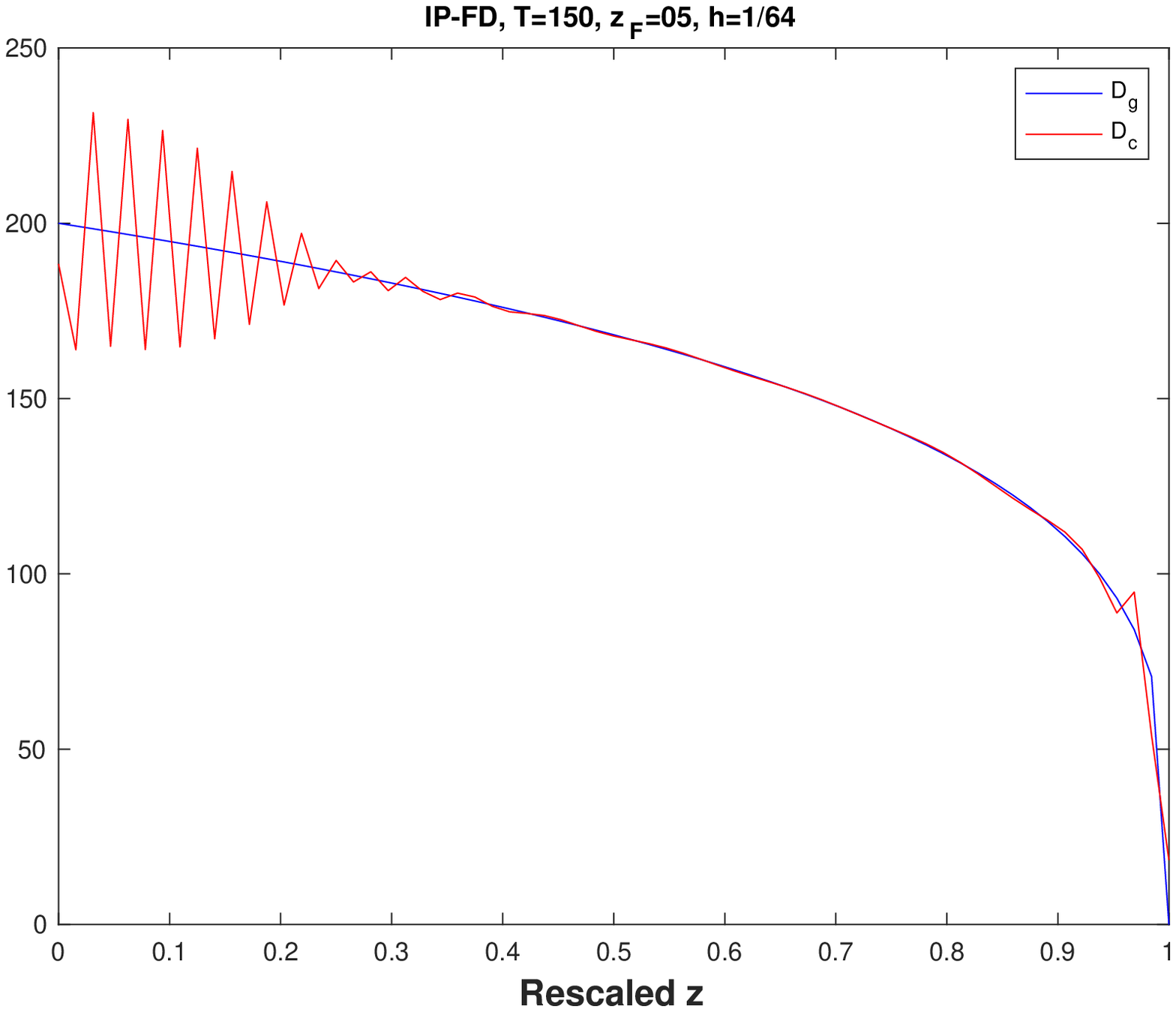}} \hspace{-10mm}&
\subfloat{\includegraphics[scale=0.4]{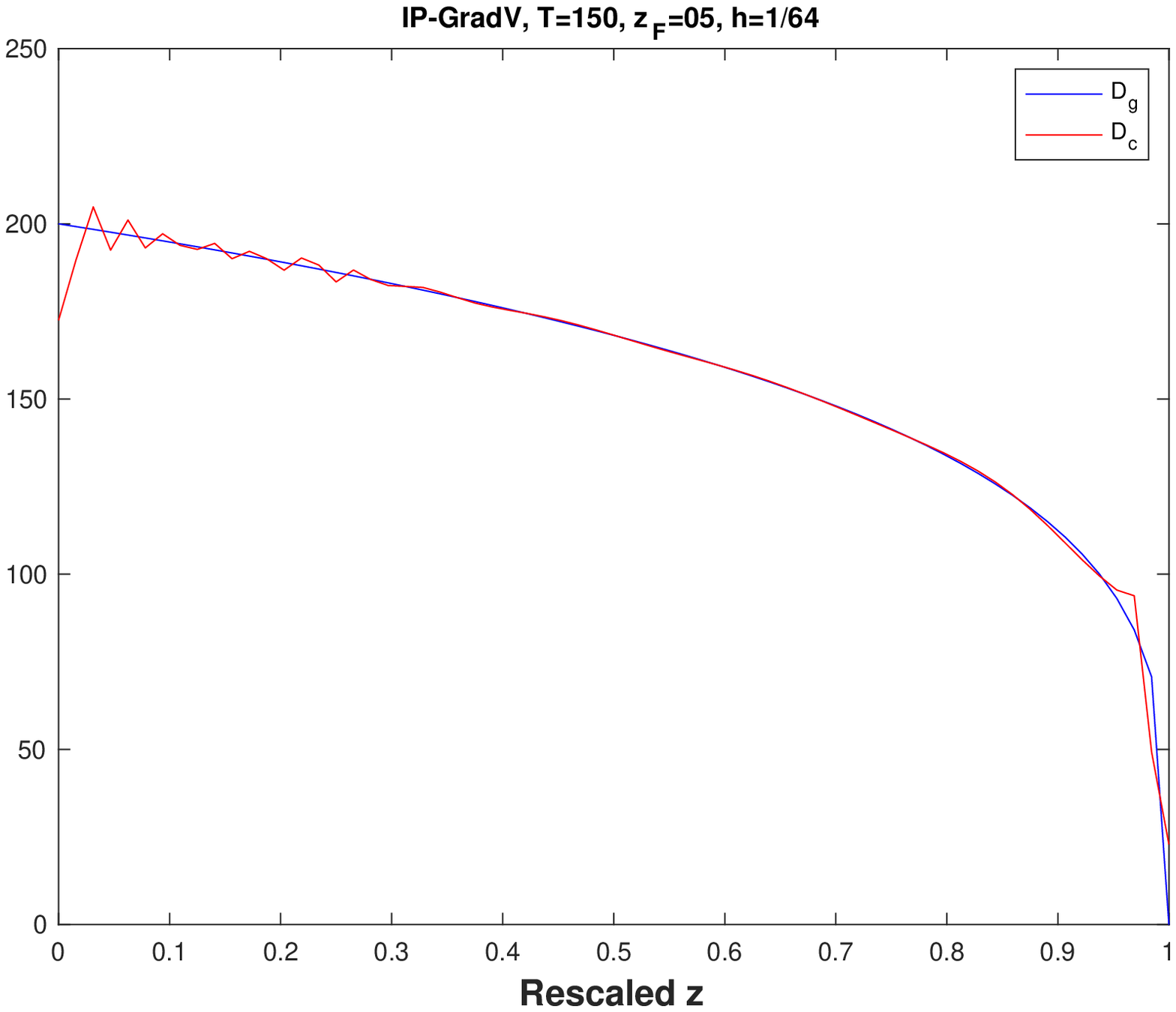}}\vspace{-5mm}\\
\subfloat{\includegraphics[scale=0.4]{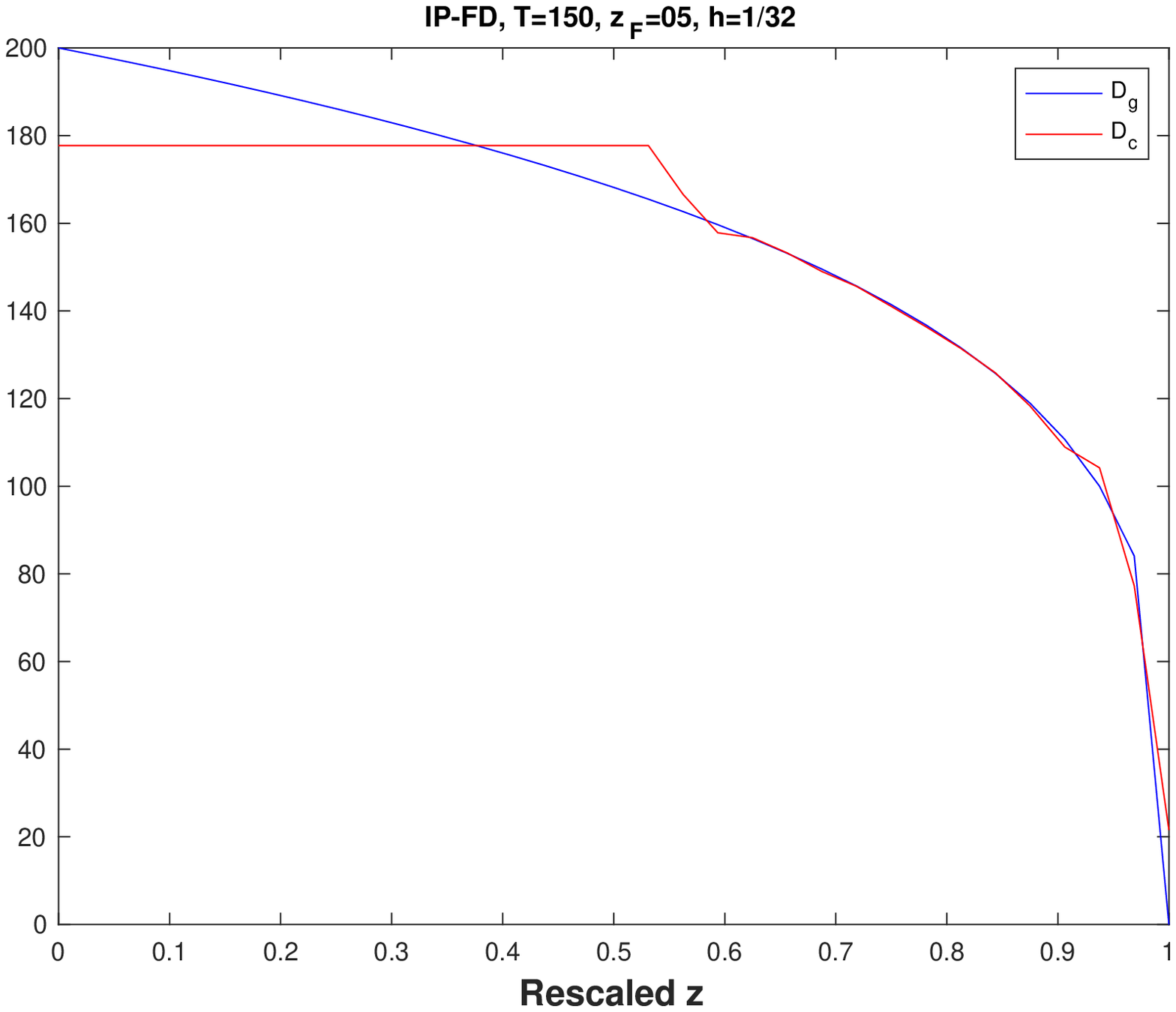}} \hspace{-10mm}&
\subfloat{\includegraphics[scale=0.4]{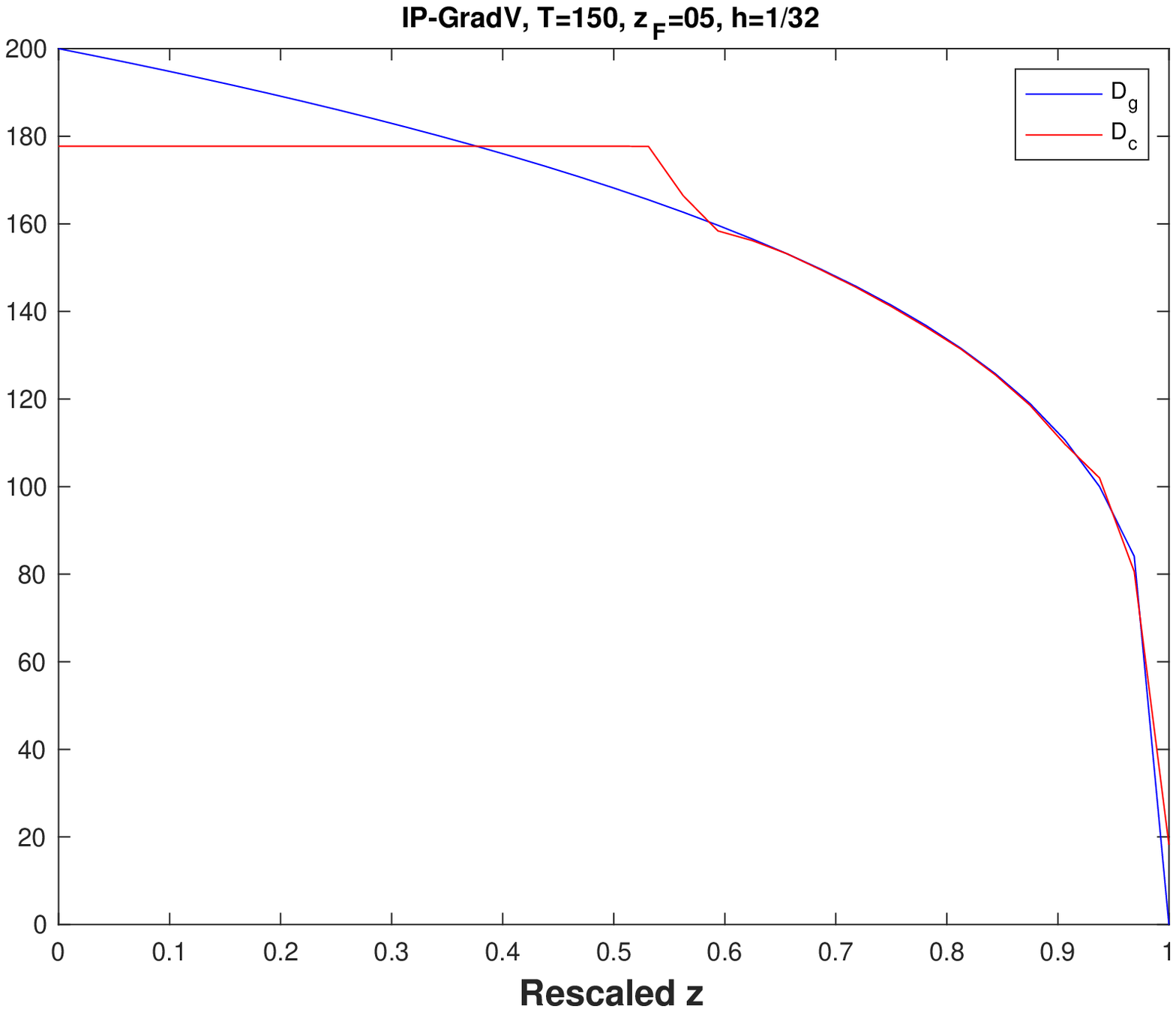}}\vspace{-5mm}\\
\subfloat{\includegraphics[scale=0.4]{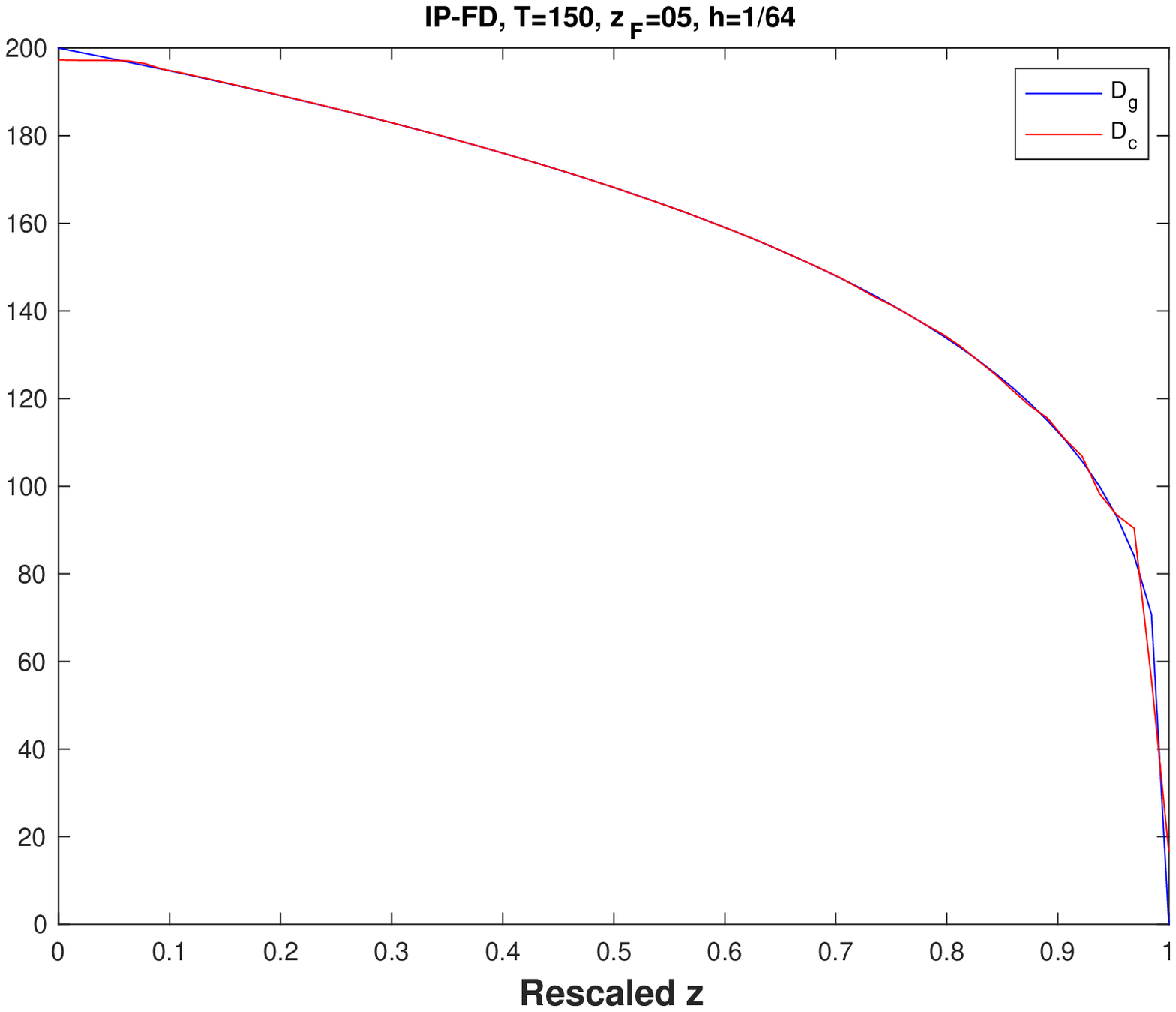}} \hspace{-10mm}&
\subfloat{\includegraphics[scale=0.4]{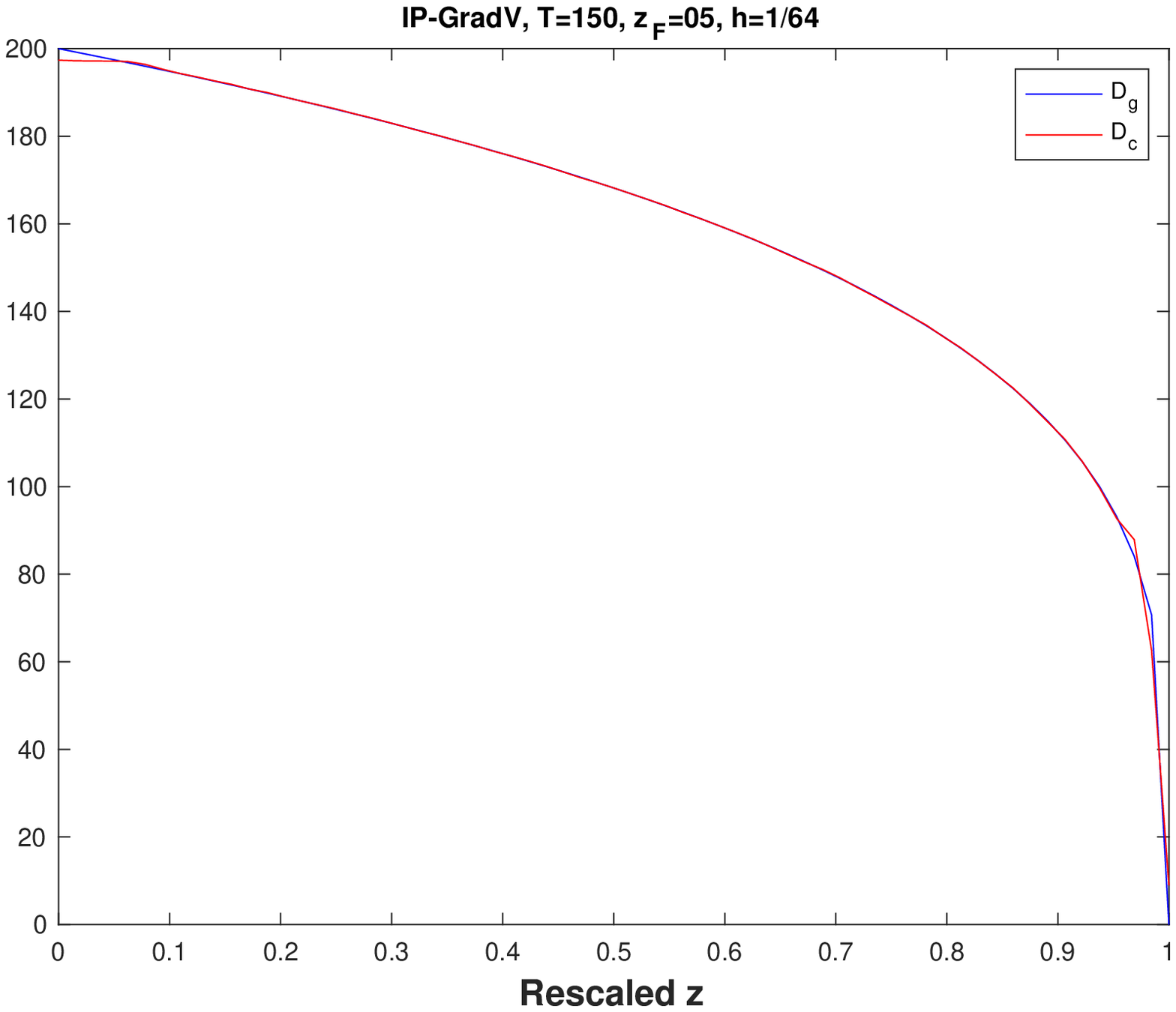}}\vspace{-5mm}\\
\end{tabular}
\centering
\caption{\it \small The $D_{CO2}$ solution for Case2d, using FminCon's IP with FD (left) and Algorithm \ref{alg:grad} (right) methods of problems \eqref{eq:min} (top) and \eqref{eq:min} with decreasing D (bottom) with initial guess $d_0 = 0$, $z_F = 5$, and $T_e = 150$, for $h = 1/32,1/64$.}\label{fig:fmin3}\vspace{-10mm}
\end{figure}

Table 
\ref{tab:inv4F} shows
the runtime of the algorithms for solving the contrained problem \eqref{eq:min} in seconds needed till convergence in (iter) iterations for stopping criteria $10^{-6}$ and the L2 relative error with the exact $D_{CO2}$. The obtained results for FminCon's SQP and IP methods for Case 
2d  
are shown, where the gradient is computed using MATLAB's built-in Finite differences or  Algorithm \ref{alg:grad}. 
Similarly to FminUnc, the corresponding versions with Algorithm \ref{alg:grad} are 10 times faster with the same order of relative errors. Moreover, the IP method has a more stable behavior than SQP as it converged for all cases with relative errors of order $10^{-1}$. The top four plots in Figure 
\ref{fig:fmin3} show the corresponding solution for IP method where the gradient is computed using MATLAB's built-in Finite differences or  Algorithm \ref{alg:grad} for Cases 
2d. 
The obtained solutions have a very similar profile with nonsmooth behavior that is reduced for a smaller mesh. 

By imposing the extra condition that $d$ is decreasing, we get better results as shown in  
Table 
\ref{tab:inv4FD}. The obtained results for FminCon's SQP and IP methods for Case 
2d 
are shown, where the gradient is computed using MATLAB's built-in Finite differences or  Algorithm \ref{alg:grad}. 
Again the corresponding versions with Algorithm \ref{alg:grad} are 10 times faster with the same order of relative errors. However, in this case  both SQP and IP methods have a stable behavior and converge for all cases with relative errors of order $10^{-1}$ and $10^{-2}$. The bottom four plots in Figure 
\ref{fig:fmin3} show the corresponding solution for IP method where the gradient is computed using MATLAB's built-in Finite differences or  Algorithm \ref{alg:grad} for Case
2d. 
It is clear that the solutions have a smoother profile than the non-decreasing case,  with smoother profile for a smaller mesh. 

\begin{table}[H]
%\centering
\small
 \hspace{-2.3cm}\begin{tabular}{|c|c|c||l|l|l|l|l|l||l|l|l|l|l|l|}
      \cline{4-15}
      \multicolumn{1}{c}{}&\multicolumn{1}{c}{}&\multicolumn{1}{c|}{}&\multicolumn{6}{c||}{FminCon  }&\multicolumn{6}{c|}{FminCon with Alg \ref{alg:grad}}\\
     \cline{4-15}
         \multicolumn{1}{c}{}& \multicolumn{1}{c}{}& \multicolumn{1}{c|}{} & \multicolumn{3}{c|}{SQP} & \multicolumn{3}{c||}{IP} & \multicolumn{3}{c|}{SQP} & \multicolumn{3}{c|}{IP}   \\ \hline
        \multicolumn{1}{|c}{$z_F$}&\multicolumn{1}{|c}{$T_e$}&\multicolumn{1}{|c||}{$h$} & Time/s & iter & L2 err & Time/s & iter & L2 err & Time/s & iter & L2 err& Time/s & iter & L2 err   \\ \hline
       \multirow{9}{*}{5} &\multirow{3}{*}{50} & 1/16 & 1.47E+0 & 291 & 4.74E-1 & 2.48E+0 & 491 & 4.82E-1 & 7.01E-1 & 258 & 4.71E-1 & 6.48E-1 & 276 & 4.81E-1  \\  \cline{3-15}
        ~ & ~ & 1/32 & 5.48E+0 & 240 & 2.92E-1 & 1.52E+1 & 650 & 2.99E-1 & 1.80E+0 & 291 & 2.90E-1 & 1.16E+0 & 288 & 2.97E-1  \\  \cline{3-15}
        ~ & ~ & 1/64 & 6.86E+1 & 413 & 3.67E+2 & 7.25E+1 & 433 & 1.08E-1 & 8.73E+1 & 1861 & 2.49E+2 & 9.16E+0 & 375 & 5.13E-2  \\  \cline{2-15}
        &\multirow{3}{*}{100} &1/16 & 8.29E-1 & 157 & 4.41E-1 & 2.10E+0 & 414 & 4.82E-1 & 6.11E-1 & 259 & 4.72E-1 & 6.60E-1 & 283 & 4.77E-1  \\  \cline{3-15}
        ~ & ~ & 1/32 & 8.01E+0 & 353 & 5.48E+1 & 1.85E+1 & 776 & 3.04E-1 & 2.13E+0 & 338 & 2.93E-1 & 1.84E+0 & 349 & 2.98E-1  \\  \cline{3-15}
        ~ & ~ & 1/64 & 1.75E+1 & 103 & 1.32E+3 & 8.49E+1 & 514 & 1.09E-1 & 7.06E+1 & 1444 & 1.06E+1 & 1.19E+1 & 480 & 1.05E-1  \\  \cline{2-15}
         &\multirow{3}{*}{150}&1/16 & 1.20E+0 & 243 & 4.62E-1 & 3.39E+0 & 655 & 5.63E-1 & 6.76E-1 & 261 & 4.83E-1 & 6.26E-1 & 262 & 4.79E-1  \\  \cline{3-15}
        ~ & ~ & 1/32 & 1.53E+1 & 648 & 1.25E+0 & 1.48E+1 & 641 & 3.02E-1 & 2.16E+0 & 347 & 3.00E-1 & 2.07E+0 & 366 & 3.03E-1  \\  \cline{3-15}
        ~ & ~ & 1/64 & 1.87E+2 & 1123 & 6.15E+0 & 6.72E+1 & 408 & 8.14E-2 & 1.38E+2 & 2644 & 2.25E+0 & 1.16E+1 & 378 & 3.33E-2  \\ \hline\hline
        %~ & ~ &
       \multirow{9}{*}{10} &\multirow{3}{*}{50} & 1/16 & 1.04E+0 & 210 & 5.90E-1 & 1.41E+0 & 272 & 6.06E-1 & 3.76E-1 & 163 & 5.71E-1 & 2.54E-1 & 163 & 5.69E-1  \\  \cline{3-15}
        ~ & ~ & 1/32 & 6.10E+0 & 266 & 3.73E-1 & 1.04E+1 & 445 & 3.89E-1 & 1.56E+0 & 257 & 3.79E-1 & 1.26E+0 & 259 & 3.73E-1  \\  \cline{3-15}
        ~ & ~ & 1/64 & 4.08E+1 & 229 & 7.19E-2 & 5.45E+1 & 327 & 9.04E-2 & 1.16E+1 & 275 & 1.04E-1 & 7.83E+0 & 342 & 8.63E-2  \\  \cline{2-15}
        %~ & ~ & 
        &\multirow{3}{*}{100} &1/16 & 1.01E+0 & 213 & 5.99E-1 & 1.24E+0 & 248 & 6.05E-1 & 5.73E-1 & 251 & 6.05E-1 & 3.14E-1 & 194 & 5.75E-1  \\  \cline{3-15}
        ~ & ~ & 1/32 & 6.03E+0 & 268 & 3.69E-1 & 1.37E+1 & 588 & 3.88E-1 & 1.45E+0 & 236 & 3.78E-1 & 1.66E+0 & 317 & 3.88E-1  \\  \cline{3-15}
        ~ & ~ & 1/64 & 9.78E+1 & 579 & 3.18E+0 & 5.69E+1 & 333 & 7.34E-2 & 1.64E+1 & 293 & 7.70E-2 & 8.23E+0 & 342 & 8.27E-2  \\ \cline{2-15}
        %~ & ~ & 
        &\multirow{3}{*}{150} &1/16 & 8.17E-1 & 169 & 5.88E-1 & 1.26E+0 & 246 & 6.05E-1 & 4.92E-1 & 215 & 5.85E-1 & 3.00E-1 & 194 & 5.71E-1  \\  \cline{3-15}
        ~ & ~ & 1/32 & 6.49E+0 & 290 & 3.82E-1 & 1.73E+1 & 733 & 4.06E-1 & 1.78E+0 & 275 & 3.88E-1 & 1.41E+0 & 315 & 3.87E-1  \\  \cline{3-15}
        ~ & ~ & 1/64 & 5.54E+1 & 334 & 9.85E-2 & 5.62E+1 & 339 & 6.77E-2 & 1.11E+1 & 263 & 6.25E-2 & 7.53E+0 & 325 & 7.56E-2 \\ \hline   
      \end{tabular}
     \caption{Convergence results of MATLAB's FminCon SQP and IP methods with built-in Finite Differences or Algorithm \ref{alg:grad} for solving Inverse problem \eqref{eq:min},  with initial guess $d_0 = 0$, for Case 2d with different $z_F, T_e,$ and $h$ values.
      }\label{tab:inv4F}
\end{table}

\begin{table}[H]
\small
 \hspace{-2.3cm}\begin{tabular}{|c|c|c||l|l|l|l|l|l||l|l|l|l|l|l|}
      \cline{4-15}
       \multicolumn{1}{c}{}&\multicolumn{1}{c}{}&\multicolumn{1}{c|}{}&\multicolumn{6}{c||}{FminCon  }&\multicolumn{6}{c|}{FminCon with Alg \ref{alg:grad}}\\
     \cline{4-15}
         \multicolumn{1}{c}{}& \multicolumn{1}{c}{}& \multicolumn{1}{c|}{} & \multicolumn{3}{c|}{SQP} & \multicolumn{3}{c||}{IP} & \multicolumn{3}{c|}{SQP} & \multicolumn{3}{c|}{IP}   \\ \hline
        \multicolumn{1}{|c}{$z_F$}&\multicolumn{1}{|c}{$T_e$}&\multicolumn{1}{|c||}{$h$} & Time/s & iter & L2 err & Time/s & iter & L2 err & Time/s & iter & L2 err& Time/s & iter & L2 err   \\ \hline
       \multirow{9}{*}{5}&\multirow{3}{*}{50} & 1/16 & 4.63E-1 & 93 & 1.58E-1 & 9.84E-1 & 199 & 1.82E-1 & 3.56E-1 & 129 & 1.57E-1 & 4.73E-1 & 183 & 1.80E-1  \\ \cline{3-15}
        ~ & ~ & 1/32 & 5.53E+0 & 252 & 5.83E-2 & 9.10E+0 & 403 & 5.90E-2 & 1.56E+0 & 214 & 6.60E-2 & 1.49E+0 & 288 & 4.89E-1  \\ \cline{3-15}
        ~ & ~ & 1/64 & 1.73E+2 & 964 & 5.96E-3 & 1.51E+2 & 840 & 1.13E-2 & 4.07E+1 & 1056 & 1.98E-2 & 2.57E+1 & 1006 & 1.07E-2  \\ \cline{2-15}
        &\multirow{3}{*}{100} & 1/16 & 6.52E-1 & 136 & 1.87E-1 & 8.36E-1 & 170 & 1.84E-1 & 4.01E-1 & 148 & 1.55E-1 & 4.23E-1 & 182 & 1.60E-1  \\ \cline{3-15}
        ~ & ~ & 1/32 & 6.25E+0 & 284 & 5.88E-2 & 8.03E+0 & 356 & 5.74E-2 & 1.61E+0 & 223 & 6.55E-2 & 1.89E+0 & 424 & 5.99E-2  \\ \cline{3-15}
        ~ & ~ & 1/64 & 1.74E+2 & 976 & 6.76E-3 & 1.47E+2 & 823 & 7.63E-3 & 4.24E+1 & 1103 & 4.63E-3 & 2.00E+1 & 948 & 1.07E-2  \\\cline{2-15}
        &\multirow{3}{*}{150} & 1/16 & 6.32E-1 & 132 & 1.86E-1 & 8.24E-1 & 166 & 1.83E-1 & 2.76E-1 & 100 & 1.57E-1 & 4.95E-1 & 202 & 1.85E-1  \\ \cline{3-15}
        ~ & ~ & 1/32 & 7.18E+0 & 328 & 6.61E-2 & 8.50E+0 & 376 & 6.09E-2 & 1.68E+0 & 232 & 6.59E-2 & 1.79E+0 & 357 & 6.00E-2  \\ \cline{3-15}
        ~ & ~ & 1/64 & 1.72E+2 & 962 & 7.11E-3 & 1.30E+2 & 722 & 1.54E-2 & 4.32E+1 & 1112 & 9.66E-3 & 2.07E+1 & 947 & 9.68E-3  \\ \hline \hline
        \multirow{9}{*}{10}&\multirow{3}{*}{50} & 1/16 & 6.74E-1 & 139 & 2.84E-1 & 1.11E+0 & 219 & 2.47E-1 & 3.64E-1 & 128 & 3.50E-1 & 3.82E-1 & 190 & 2.42E-1  \\ \cline{3-15}
        ~ & ~ & 1/32 & 1.84E+0 & 83 & 2.22E-1 & 7.61E+0 & 335 & 1.43E-1 & 1.30E+0 & 179 & 2.89E-1 & 1.86E+0 & 388 & 1.39E-1  \\ \cline{3-15}
        ~ & ~ & 1/64 & 1.11E+2 & 613 & 3.90E-2 & 1.69E+2 & 939 & 5.81E-2 & 2.36E+1 & 622 & 9.34E-2 & 2.25E+1 & 836 & 2.80E-2  \\ \cline{2-15}
        &\multirow{3}{*}{100} & 1/16 & 5.35E-1 & 109 & 3.21E-1 & 1.03E+0 & 206 & 2.29E-1 & 3.31E-1 & 121 & 3.45E-1 & 5.07E-1 & 220 & 2.53E-1  \\ \cline{3-15}
        ~ & ~ & 1/32 & 1.36E+0 & 61 & 1.68E-1 & 9.08E+0 & 399 & 1.41E-1 & 7.13E-1 & 99 & 2.82E-1 & 1.90E+0 & 392 & 1.52E-1  \\ \cline{3-15}
        ~ & ~ & 1/64 & 1.51E+2 & 838 & 4.05E-2 & 1.53E+2 & 854 & 1.45E-2 & 2.97E+1 & 790 & 3.15E-2 & 2.17E+1 & 887 & 3.33E-2  \\ \cline{2-15}
        &\multirow{3}{*}{150} & 1/16 & 6.04E-1 & 125 & 2.72E-1 & 1.01E+0 & 203 & 2.48E-1 & 3.71E-1 & 125 & 3.40E-1 & 5.45E-1 & 188 & 2.46E-1  \\ \cline{3-15}
        ~ & ~ & 1/32 & 3.38E+0 & 154 & 1.42E-1 & 1.04E+1 & 458 & 1.31E-1 & 1.64E+0 & 199 & 1.58E-1 & 1.82E+0 & 336 & 1.02E-1  \\ \cline{3-15}
         ~& ~& 1/64 & 1.54E+2 & 865 & 4.33E-2 & 1.58E+2 & 884 & 3.73E-2 & 3.24E+1 & 799 & 3.05E-2 & 2.70E+1 & 1036 & 2.65E-2 \\ \hline
         \end{tabular} 
         \caption{Convergence results of MATLAB's FminCon SQP and IP methods with built-in Finite Differences or Algorithm \ref{alg:grad} for solving Inverse problem  with decreasing $D$,  for Case 2d with different $z_F, T_e,$ and $h$ values and  initial guess $d_0 = 0$. }\label{tab:inv4FD}
\end{table}

\vspace{-4mm} 

\section{Conclusion}\label{sec:conc}
In this paper, we study the theoretical and computational aspects of the Firn direct problem. Then, the inverse problem is formulated by defining the objective function in section, where its gradient is computed using directional derivatives in a computationally efficient way. Testings on the efficiency of the computed gradient using MATLAB's FminUnc and FminCon functions and Nonlinear Conjugate Gradient method validate that it gives similar results to the case where the gradient is computed using Finite difference but in less runtime (up to 10 times faster).

This work lays the ground for constructing a robust inverse problem algorithm that should extract past history diffusion coefficients of different gases of interest for understanding climate changes, and which will be the content of a  forthcoming work . 

\newpage
\section*{Appendix: Existence and Uniqueness of Solution to the Discrete Problem \eqref{eq:mat2}} 
\noindent{
{To prove that the Discrete Problem \eqref{eq:mat2}, with the approximated matrices $A(D)$ and $S(D)$, has a unique solution,  we seek the invertibility of the matrix $M+T_e\Delta t C$, looking for conditions that would make it  positive definite.\\
Let $v\not=0\in \mathbb{R}^{n-1}$. Starting with the identity:
\begin{eqnarray}
    v^T(M+T_e\,\Delta t\, C)v =v^T\left(M+T_e\,\Delta t\, \left( \dfrac{\mathcal{G}}{f} M+\dfrac{1}{z_F^2f}S(D_\alpha)-\dfrac{\mathcal{M}_\alpha}{z_Ff}A(D_\alpha)+\dfrac{1}{z_F}Q\right)\right)v
\end{eqnarray}

\noindent one obtains:
\begin{equation}\label{eq:v}
 v^T(M+T_e\,\Delta t\, C)v \;{\color{black}\geq} \;\left(1+T_e\,\Delta t\; \dfrac{\mathcal{G}}{f}\right) \Mnorm{v}^2 +\dfrac{T_e}{z_F f}\,\Delta t\, v^T G\,v   
\end{equation}
{\color{black} since $\dfrac{1}{z_F}v^TQv \geq 0$ by lemma \ref{lemma:27}, } with $G=\dfrac{1}{z_F}S(D_\alpha)-{\mathcal{M}_\alpha}A(D_\alpha)\in\mathbb{R}^{n-1}$.

\noindent Moreover, associated with the mass matrix $M$,  the norm:
\begin{equation}\label{eq:M_norm}
\Mnorm{v}=\{v^TMv\}^{1/2}, 
\end{equation}
which is equivalent to the weighted $l^2$ norm:
\begin{equation}
\lhnorm{v} =\{\sum_{i=1}^{n-1}{h|v_i|^2}\}^{1/2},
\end{equation}
as shown in the following Lemma.
\begin{lemma}
The $M$-norm, $\Mnorm{.}$, and the weighted $l^2$ norm, $\lhnorm{.}$, satisfy the inequalities:
\begin{equation}\label{eq:norm_eq}
   \dfrac{1}{6}\lhnorm{v}^2\le  v^TMv\le \lhnorm{v}^2
\end{equation}
\end{lemma}
\begin{proof} From Lemma \ref{M:spd} one has:
$$        v^TMv = 
\dfrac{h}{3} \sum\limits_{i=1}^{n-2}(v_i^2+v_iv_{i
+1}+v_{i+1}^2). $$
By Geometric inequality ($ ab \leq  \frac{a^2}{2} + \frac{b^2}{2}$) and definition of weighted $l^2$ norm we get
    \begin{eqnarray}
    v^TMv &=& \dfrac{h}{3} \sum\limits_{i=1}^{n-2} (v_i^2+v_iv_{i+1}+v_{i+1}^2) \nonumber\\
    &\leq&  \dfrac{h}{3} \sum\limits_{i=1}^{n-2} (v_i^2+0.5v_i^2+0.5v_{i+1}^2+v_{i+1}^2) =  \dfrac{h}{2} \sum\limits_{i=1}^{n-2} (v_i^2+v_{i+1}^2)  \;\leq\;  \lhnorm{v}^2 \nonumber\\
     v^TMv &=&\dfrac{h}{6}  \sum\limits_{i=1}^{n-2}(v_i+v_{i+1})^2 +  \dfrac{h}{6}\sum\limits_{i=1}^{n-2}v_i^2 + \dfrac{h}{6} \sum\limits_{i=1}^{n-2}v_{i+1}^2 \nonumber\\
     &=& \dfrac{h}{6}  \sum\limits_{i=1}^{n-2}(v_i+v_{i+1})^2 +  \dfrac{h}{6}\sum\limits_{i=2}^{n-2}v_i^2 + \dfrac{1}{6} \lhnorm{v}^2 \;\geq \;  \dfrac{1}{6} \lhnorm{v}^2
    \end{eqnarray}
\end{proof}}
\noindent We prove now the following lemma.
\begin{lemma}\label{lem:vGv} $v^TGv \,=\, \mbox{Term}_1 - \mbox{Term}_2, $
where\\ $\displaystyle \mbox{Term}_1 =\dfrac{1}{2hz_F}\left(\sum_{i=2}^{n-1}{(D_i+D_{i+1})(v_i-v_{i-1})^2+(D_1+D_2)v_1^2}\right)$, \\$\displaystyle 
\mbox{Term}_2=\dfrac{\mathcal{M}_\alpha}{4}\left(\sum_{i=1}^{n-2}{}(D_i-D_{i+2})v_i^2+(D_n+D_{n-1})v_{n-1}^2\right).$
\end{lemma}
\begin{proof}
    This result is straightforwardly obtained from the expressions of the matrices $S(D_\alpha)$ and $A(D_\alpha)$.
\end{proof}
\noindent We focus now on Term$_1$.
For that purpose, let $D_{i+1/2}=\dfrac{1}{2}(D_i+D_{i+1})$ and introduce $v_0=0$, then $\mbox{Term}_1$ can be rewritten as:
$$\mbox{Term}_1=\dfrac{1}{z_F}\sum_{i=1}^{n-1}{hD_{i+1/2}\left(\dfrac{v_i-v_{i-1}}{h}\right)^2}.$$
\noindent We derive now two inequalities that can provide lower bounds for $\mbox{Term}_1$.
\begin{lemma}\label{lemma211} Under the continuity and positivity assumption of $D_\alpha(z)$, \eqref{eq:Dpos},  there exists an $h_0<1$ ($n_0>1$), such that for $h<h_0$ ($n>n_0$)    one has
    \begin{equation}\label{eq:estv}
    |v_i|^2\leq \,c_D\,I(D_\alpha)\,\mbox{Term}_1, \quad \forall i=1,...,n-1\,
\end{equation}
leading to \begin{equation}\label{eq:estv2}
\linftynorm{v}^2=\max_{i}{|v_i|^2}\le c_D\,I(D_\alpha)\,\mbox{Term}_1,\end{equation}
 and
 a generalized Poincaré's inequality
\begin{equation}\label{eq:genpoinc}
  \dfrac{1}{c_D\,I(D_\alpha)}\lhnorm{v}^2\le \mbox{Term}_1,  \vspace{-3mm}
\end{equation}
where $\displaystyle I(D_\alpha)=z_F\int\limits_0^{1}\dfrac{1}{D_\alpha}$ and $c_D$  a positive constant independent from $h$.
\end{lemma}
\begin{proof}
    Let $v_0=0$, then based on the identity 
    $$v_i=\sum_{k=1}^i{(v_k-v_{k-1})}=\sum_{k=1}^i{(hD_{k+1/2})^{1/2}\left(\dfrac{v_k-v_{k-1}}{h}\right)\dfrac{h^{1/2}}{D^{1/2}_{k+1/2}}},$$
    and using Schwarz inequality, one deduces that for $i=1,..,n-1$ :
    $$|v_i|^2\leq\sum_{k=1}^i{hD_{k+1/2}\left(\dfrac{v_k-v_{k-1}}{h}\right)^2} \sum_{k=1}^i{\dfrac{h}{D_{k+1/2}}}\,.$$
    Consequently, one has for $,\,i=1,..,n-1$:
    $$|v_i|^2\leq z_F\mbox{Term}_1\sum_{k=1}^{n-1}{\dfrac{h}{D_{k+1/2}}} = \mbox{Term}_1 \,\mathcal{S}_n,$$ 
    where $\displaystyle \mathcal{S}_n = z_F \,\sum_{k=1}^{n-1}{\dfrac{h}{D_{k+1/2}}}$.\\
    As the sequence $\mathcal{S}_n$ converges to $I(D _\alpha)$, i.e 
    $$ \lim\limits_{{h\to 0} \\{(n\to\infty)}} \mathcal{S}_n = \lim\limits_{{h\to 0} \\{(n\to\infty)}} z_F{\sum_{k=1}^{n-1}{\dfrac{h}{D_{k+1/2}}}}=I(D_\alpha)=z_F\int_0^{1}{\dfrac{1}{D_\alpha(z)}} \;<\; \infty$$
    then there exists an $h_0<1$ ($n_0>1$), such that for $h<h_0$ ($n>n_0$)    one has,
    \begin{equation}
        \mathcal{S}_n \leq c_D \, I(D_\alpha),
 \end{equation}
 where $0< c_D < 2$. 
    Then, for $\,i=1,..,n-1$:
    $$|v_i|^2\leq c_D\,I(D_\alpha)\,\mbox{Term}_1,$$
    and consequently:
    $$\lhnorm{v}^2=\sum_{i=1}^{n-1}{h|v_i|^2}\leq c_D\,I(D_\alpha)\,\mbox{Term}_1.$$
\end{proof}
\noindent We now consider Term$_2$.
We prove the following result.
\begin{lemma}\label{lemma212} Under the continuity, positivity \eqref{eq:Dpos}, and Lipschitz continuity \eqref{eq:Lipschitz} assumptions on $D_\alpha(z)$,  there exists an $h_0<1$ ($n_0>1$), and $D_{n-2}< \epsilon < D_{n-3}$,  such that for $h<h_0$ ($n>n_0$)    one has
$$|\mbox{Term}_2| \leq \dfrac{\mathcal{M}_\alpha}{2}\left(L_\delta \,\lhnorm{v}^2+ 2 \,\epsilon \,\linftynorm{v}^2\right)$$
and consequently from estimate \eqref{eq:estv2},
$$|\mbox{Term}_2| \leq\dfrac{\mathcal{M}_\alpha}{2}\left(L_\delta\, \lhnorm{v}^2+ 2\,\epsilon\,c_D\,I(D_\alpha)\,\mbox{Term}_1\right).$$
\end{lemma}
 \begin{proof} 
Let $\epsilon >0$, the there exists $\delta := \delta(\epsilon)$ such that \begin{equation}\label{eq:maxeps}
\max\limits_{z_i\in [z_F-\delta,z_F]} D(z_i) < \epsilon\end{equation}
given that $\lim\limits_{z \rightarrow z_f} D(z) = 0$ and $D(.)$ is continuous. \\
Thus, Term$_2$ can be expressed as  
\begin{eqnarray}
\mbox{Term}_2 &=&\dfrac{\mathcal{M}_\alpha}{4}\sum\limits_{i \,\notin \, \mathcal{
O}} (D_i-D_{i+2})v_i^2+\dfrac{\mathcal{M}_\alpha}{4}\sum\limits_
{i\, \in \,\mathcal{
O}} (D_i-D_{i+2})v_i^2+\dfrac{\mathcal{M}_\alpha}{2}D_{n-1/2}v_{n-1}^2 \nonumber\\
&=&\dfrac{\mathcal{M}_\alpha}{4}\sum\limits_{i=1}^{i_o-1} (D_i-D_{i+2})v_i^2+\dfrac{\mathcal{M}_\alpha}{4}\sum\limits_
{i=i_o}^{n-2} (D_i-D_{i+2})v_i^2+\dfrac{\mathcal{M}_\alpha}{2}D_{n-1/2}v_{n-1}^2 \nonumber
\end{eqnarray} 
where \vspace{-8mm}\begin{eqnarray}\mathcal{
O} &=& \{i \,|\, z_i \in [z_F - \delta, z_F] \; \&  \; z_{i+2} \in [z_F - \delta, z_F], \;\; for \; i = 1,2, \cdots , n-2\}\nonumber\\
&=& \{i_{{o}}, i_{{o}} + 1, \cdots, n-2\}.\nonumber
\end{eqnarray}
Then, using  the triangle inequality, \eqref{eq:maxeps}, and Lipschitz continuity where $|D_i - D_{i+2}| \leq 2hL_\delta$, we get 
\begin{eqnarray}
   |\mbox{Term}_2| &\leq&
\dfrac{\mathcal{M}_\alpha}{2}\sum\limits_{i=1}^{i_o-1} \,h\,L_\delta |v_i|^2+\dfrac{\mathcal{M}_\alpha}{2}\sum\limits_
{i=i_o}^{n-2} \epsilon \, |v_i|^2+\dfrac{\mathcal{M}_\alpha}{2} \epsilon \, |v_{n-1}|^2 
\label{eq:deltazf}  \\
    &\leq& \dfrac{\mathcal{M}_\alpha}{2}   L_\delta \, \lhnorm{v}^2 +\dfrac{\mathcal{M}_\alpha}{2}\sum\limits_
{i=i_o}^{n-1} \epsilon \, \linftynorm{v}^2 
\end{eqnarray}
Note that $\epsilon$ could be chosen very small, such that $\mathcal{O} = \{n-2\}$, i.e. $D_{n-2}< \epsilon < D_{n-3}$.\\ In this case we get,
\begin{eqnarray}
   |\mbox{Term}_2| 
    &\leq& \dfrac{\mathcal{M}_\alpha}{2}   L_\delta \, \lhnorm{v}^2 +{\mathcal{M}_\alpha} \epsilon \, \linftynorm{v}^2 
\end{eqnarray}

\end{proof}
We can now deduce a lower bound estimate on $v^TGv$.
\begin{lemma} Under the continuity, positivity \eqref{eq:Dpos}, and Lipschitz continuity \eqref{eq:Lipschitz} assumptions on $D_\alpha(z)$,  there exists an $h_0<1$ ($n_0>1$), $D_{n-2}< \epsilon < D_{n-3}$, and there exist a constant $K_G$,  such that for $h<h_0$ ($n>n_0$)    one has
$$v^TGv\ge -|K_G|. \Mnorm{v}^2.$$
\end{lemma}
\begin{proof}
$$v^TGv=\mbox{Term}_1-\mbox{Term}_2\ge (1-\epsilon\,\mathcal{M}_\alpha\, \,c_D\,I(D_\alpha))\,\mbox{Term}_1-\dfrac{\mathcal{M}_\alpha}{2} \,L_\delta\,\lhnorm{v}^2, $$
Selecting $h_0$ to be sufficiently small ($n_0$ sufficiently large) so that 
$$\epsilon\,\mathcal{M}_\alpha\, \,c_D\,I(D_\alpha) < 0.5$$ 
and using \eqref{eq:genpoinc},  and the equivalences between $\lhnorm{v}$ and $\Mnorm{v}$ {\eqref{eq:norm_eq}}, one has:
\begin{eqnarray}
    v^TGv &\ge& \dfrac{\lhnorm{v}^2}{2\,c_D\,I(D_\alpha)}-\dfrac{\mathcal{M}_\alpha}{2}L_\delta\lhnorm{v}^2 \;\geq\; -|K_G| \,\lhnorm{v}^2 \\
     &\ge&  -|K_G| \, \Mnorm{v}^2 \end{eqnarray}
where $K_G =   \dfrac{1}{2\,c_D\,I(D_\alpha)}-\dfrac{\mathcal{M}_\alpha}{2}L_\delta $ can be either positive or negative.
\end{proof}
\noindent We are now ready to prove the invertibility of $M+\Delta T_eC$, the matrix of the discrete system. 
\begin{theorem}%[{\bf Invertibility and Stability}] 
The matrix $M+T_e\,\Delta t\, C$ of the discrete system is invertible if $h$ is chosen to be sufficiently small and if \,$0 < \Delta t <  \dfrac{z_Ff}{T_e |\,\,|K_G| - {\,z_F} {\mathcal{G}}\,|}$.
    \end{theorem}
\begin{proof}
    Since from \eqref{eq:v}, one has:
    $$ v^T(M+T_e\,\Delta t\, C)v \;{\color{black} \geq} \;\left(1+T_e\,\Delta t\; \dfrac{\mathcal{G}}{f}\right) \Mnorm{v}^2 +\dfrac{T_e}{z_F f}\,\Delta t\, v^T G\,v,  $$
    then using the result of the previous lemma, we deduce for $h$ sufficiently small that:
    $$ v^T(M+T_e\,\Delta t\, C)v\ge \left(1+T_e\,\Delta t\; \dfrac{\mathcal{G}}{f}\right) \Mnorm{v}^2 -|K_G|\dfrac{T_e}{z_F f}\,\Delta t \Mnorm{v}^2,  $$
    i.e.,
$$ v^T(M+T_e\,\Delta t\, C)v\ge \left(1+\dfrac{T_e}{f}\,\Delta t\; (\mathcal{G}  -|K_G|\dfrac{1}{z_F})\right) \Mnorm{v}^2$$  
{\color{black} $M+T_e\,\Delta t\, C$ is positive definite 
if  
\begin{eqnarray} 1+\dfrac{T_e}{f}\,\Delta t\; \left(\mathcal{G}  -|K_G|\dfrac{1}{z_F}\right) &>& 0\\
\iff \Delta t\; \left({z_F}\mathcal{G}  -|K_G|\right) &>& -  z_F\dfrac{f}{T_e} \label{eqdt}\\
\Delta t &<&  \dfrac{z_Ff}{T_e (|K_G| - {\,z_F} {\mathcal{G}})} \label{eqdt2}
\end{eqnarray} Assuming  ${z_F}\, {\mathcal{G}}-{ |K_G|}< 0$ in \eqref{eqdt} leads to condition \eqref{eqdt2} on $\Delta t$.\\ If ${z_F} {\mathcal{G}}-{ |K_G|}\geq 0$, then it is sufficient to chose $\Delta t  >0$. 
Thus, setting $$0 < \Delta t <  \dfrac{z_Ff}{T_e |\,\,|K_G| - {\,z_F} {\mathcal{G}}\,|}$$ guarantees the existence of a unique solution to the Discrete Problem \eqref{eq:mat2} irrespective of the sign of ${\,z_F}\, {\mathcal{G}}-{ |K_G|}$.}
\end{proof}
}

\end{document}